\documentclass[11pt,a4paper]{article}
\usepackage[T1]{fontenc}
\usepackage{lmodern}
\usepackage[margin=27mm]{geometry}
\usepackage{amsmath,amssymb,amsthm,mathtools,booktabs,array}
\usepackage{microtype,xcolor,needspace}
\usepackage[colorlinks=true,linkcolor=blue!55!black,citecolor=blue!55!black,urlcolor=blue!55!black]{hyperref}
\hypersetup{pdftitle={Log-concavity of rank and crank statistics for integer partitions}}
\numberwithin{equation}{section}
\newtheorem{theorem}{Theorem}[section]
\newtheorem{proposition}[theorem]{Proposition}
\newtheorem{conjecture}[theorem]{Conjecture}
\newtheorem{lemma}[theorem]{Lemma}
\theoremstyle{remark}

\newtheorem{claim}{Claim}[section]

\newcommand{\A}{\mathcal A}
\newcommand{\dd}{\,\mathrm d}

\newcommand{\sech}{\operatorname{sech}}

\allowdisplaybreaks[2]
\title{Proof of two conjectures on the log-concavity
of partition ranks and cranks}

\author{%
	Jianxi Mao\textsuperscript{1},
	Wenle Shi\textsuperscript{1},
	Bao-Xuan Zhu\textsuperscript{2}%
	\thanks{Corresponding author.
		\nolinkurl{maojx@dlut.edu.cn} (J. Mao),
		\nolinkurl{shi-wenle@hotmail.com} (W. Shi) and
        \nolinkurl{bxzhu@jsnu.edu.cn}(B.-X. Zhu).}
	\\[0.8em]
  {\small\textsuperscript{1}\ School of Mathematical Sciences,
	Dalian University of Technology,}\\
{\small Dalian 116024, P. R. China}\\[0.3em]
{\small\textsuperscript{2}\ School of Mathematics and Statistics,
	Jiangsu Normal University,}\\
{\small Xuzhou 221116, P. R. China}
}

\date{}

\begin{document}
\maketitle
	\begin{abstract}
	Let $M(m,n)$ and $N(m,n)$ denote the numbers of partitions of $n$
	with crank $m$ and rank $m$, respectively.
	Ji and Zang proved that the crank sequence
	$(M(m,n))_{|m|\le n-1}$ is unimodal for every $n\ge44$.
	Bringmann, Man, and Rolen subsequently established the unimodality
	of the rank sequence $(N(m,n))_{|m|\le n-2}$ for every $n\ge39$.
	Ji and Zang also conjectured that
	$(M(m,n))_{|m|\le n-71}$ is log-concave for every $n\ge72$,
	and that $(N(m,n))_{|m|\le n-72}$ is log-concave
	for every $n\ge73$.
	In this paper, we prove both conjectures.
	\end{abstract}
\noindent\textbf{2020 Mathematics Subject Classification.}
05A17, 05A20, 11P82.

\smallskip
\noindent\textbf{Keywords.}
Integer partitions; rank; crank; strict log-concavity; explicit estimates.

\setcounter{tocdepth}{1}
\tableofcontents
\clearpage

\section{Introduction}
The rank and crank are two fundamental statistics in the theory of partitions.
A \emph{partition} $\lambda=(\lambda_1,\ldots,\lambda_r)$ of $n$ is a finite
sequence satisfying
$$
\lambda_1\ge \lambda_2\ge\cdots\ge \lambda_r>0\qquad \textrm{and} \qquad \lambda_1+\lambda_2+\cdots+\lambda_r=n. 
$$ 
The \emph{rank} of $\lambda$, introduced by Dyson~\cite{Dyson}, is
defined by
\[
\textrm{rank}(\lambda)=\lambda_1-r.
\]
Let $\omega(\lambda)$ denote the number of ones in $\lambda$ and
$\mu(\lambda)$ denote the number of parts of $\lambda$ larger than $\omega(\lambda)$.
The \emph{crank} of $\lambda$, defined by Andrews and
Garvan~\cite{AG}, is
\[
\textrm{crank}(\lambda)=
\begin{cases}
	\lambda_1, & \text{if } \omega(\lambda)=0,\\[2pt]
	\mu(\lambda)-\omega(\lambda), & \text{if } \omega(\lambda)>0.
\end{cases}
\]

Let $p(n)$ denote the number of partitions of $n$. 
It was conjectured by Dyson~\cite{Dyson} and confirmed by
Atkin and Swinnerton-Dyer~\cite{ASD} that the rank of a partition
could explain two of Ramanujan's famous partition congruences
$p(5n+4)\equiv 0\pmod{5}$ and $p(7n+5)\equiv 0\pmod{7}$.
However, the rank fails to explain Ramanujan's congruence modulo 11.
Therefore Dyson~\cite{Dyson} conjectured the existence of another statistic which would explain the third Ramanujan congruence. 
Andrews and Garvan~\cite{AG,GarvanCrank} successfully found the crank, and proved that the crank simultaneously 
explains the above three Ramanujan congruences.

A sequence of nonnegative numbers $(c_n)_{n\ge 0}$
is said to be \emph{unimodal} if there exists a $p$ such that
\[
c_0 \le c_1 \le \cdots \le c_{p-1} \le c_p
\ge c_{p+1} \ge \cdots.
\]
We say $(c_k)_{k\ge 0}$ is \emph{log-concave} if 
$c_k^2-c_{k-1}c_{k+1}\ge 0$  for $k\ge 1$.
Clearly, log-concavity with no internal zeros implies unimodality (or monotonicity).
It is well known that the partition function satisfies $p(n+1)\ge p(n)$ for $n\ge0$.
Nicolas~\cite{Nicolas} and DeSalvo and Pak~\cite{DP}
independently proved that $(p(n))_{n}$ is log-concave for $n\ge 26$.

Let $N(m,n)$ and $M(m,n)$ denote the numbers
of partitions of $n$ with rank $m$ and crank $m$, respectively.
The unimodality of rank and crank has been studied
from both combinatorial and analytic perspectives. 
For the rank sequence, 
Chan and Mao~\cite{ChanMao} proved that $N(m,n)\ge N(m+2,n)$
for $m,n\ge 0$.
Thus the even and odd rank subsequences are separately unimodal.
They ask the question to find conditions on the unimodality of $(N(m,n))_{m}$. 
Bringmann, Man and Rolen~\cite[Theorem~1.6]{BMR} answered the question and 
proved that the sequence
$(N(m,n))_{m}$
is unimodal for every $n\ge39$ and $0\le m\le n-3$.
For the crank sequence,
Kim, Kim, and Seo~\cite{KKS} proved that, for each fixed
nonnegative integer $m$,  $M(m,n)>M(m+1,n)$
holds for all sufficiently large $n$.
Ji and Zang~\cite{JZ} subsequently proved that the sequence
$(M(m,n))_{m}$ is unimodal for every $n\ge44$ and $0\le m\le n-2$.

Ji and Zang posed the following conjectures on the log-concavity of 
$(M(m,n))_{m}$ and
$(N(m,n))_{m}$.  See also Bringmann, Jennings-Shaffer, and Mahlburg~\cite{BJM}.
\begin{conjecture}[\cite{JZ}, Conjecture~11.1]\label{conj:crank}
	For every integer $n\ge72$ and every integer $m$ with $|m|\le n-72$,
	\[
	M(m,n)^2\ge M(m-1,n)M(m+1,n).
	\]
\end{conjecture}

\begin{conjecture}[\cite{JZ}, Conjecture~11.3]\label{conj:rank}
	For every integer $n\ge73$ and every integer $m$ with $|m|\le n-73$,
	\[
	N(m,n)^2\ge N(m-1,n)N(m+1,n).
	\]
\end{conjecture}

Bringmann and Dousse \cite{BD} proved Dyson's logistic asymptotic for the crank;
Dousse and Mertens \cite{DM} obtained its rank counterpart. On each fixed range $|m|\le L\sqrt n$, these results give, uniformly as $n\to\infty$,
\begin{align}\label{eq:logistic-motivation}
	M(m,n),\ N(m,n)
	&\sim \frac{\pi p(n)}{4\sqrt{6n}}
	\sech^2\!\left(\frac{\pi m}{2\sqrt{6n}}\right).
\end{align}
The density $f(z)=\frac14\sech^2(z/2)$ satisfies
$-(\log f)''(z)=2f(z)>0$. 

Liu and Zhou~\cite{LiuZhou} introduced Garvan's
$k$-rank functions $A_k(m,n)$, where
\begin{equation}\label{eq:Ak}
	A_1(m,n)=M(m,n), \qquad A_2(m,n)=N(m,n).
\end{equation}
Zhou~\cite{Zhou} proved that, for $k\in\{1,2\}$ and integers $n\ge142+2k$ and $m$
satisfying $\frac n2-2k+2\le |m|\le n-k-71$, 
\begin{align}\label{eq:zhou}
	A_k(m,n)^2>A_k(m-1,n)A_k(m+1,n).
\end{align}
The following theorem establishes strict log-concavity in the two conjectures.

\begin{theorem}\label{thm:strict}
	For $k\in\{1,2\}$ and integers $n\ge72$ and $m$
	with $|m|\le n-k-71$, we have
	\begin{align}\label{eq:strict-main}
		A_k(m,n)^2
		>
		A_k(m-1,n)A_k(m+1,n).
	\end{align}
\end{theorem}

The case $k=1$ and $k=2$ in Theorem~\ref{thm:strict} imply Conjectures~\ref{conj:crank}
and~\ref{conj:rank} respectively.

Ji and Zang~\cite[Section~11]{JZ} have checked Conjecture~\ref{conj:crank} for $72\le n\le10000$.
We check that the inequality~\eqref{eq:strict-main} is strict for both $k=1,2$ and $n\le 5000.$
It suffices to prove Theorem~\ref{thm:strict} for $n\ge 5000$.
The paper is organized as follows. In Section~\ref{sec:transfer}, we establish explicit estimates for the
partition function, its smooth approximation, and the first
difference of this approximation.
In Section~\ref{sec:central}, we prove
Theorem~\ref{thm:strict} for $n\ge5000$ and
$0\le m\le2\sqrt n$.
Section~\ref{sec:outer} establishes the inequality in
Theorem~\ref{thm:strict} for $n\ge 5000$ and
$2\sqrt n\le m\le n/2+1$.
Finally,  we complete the proof of Theorem~\ref{thm:strict} in Section~\ref{sec:completion-conjectures},
by combining these estimates with Zhou's result \eqref{eq:zhou}.

\section{Preliminaries}\label{sec:transfer}

By convention, $p(0)=1$ and $p(n)=0$ for $n<0$.
For $|q|<1$,
\begin{equation}\label{gf:p}
	\sum_{n\ge0}p(n)q^n=\frac{1}{(q;q)_\infty}=\frac{1}{\prod_{j\ge1}(1-q^j)}
\end{equation}
For $k\in\{1,2\}$, integers $n\ge2$ and $m\ge0$,
$A_k(m,n)$ can be expressed as the following finite alternating sum of partition numbers:
\begin{equation}\label{eq:coeff}
	\begin{aligned}
		A_k(m,n)
		=\sum_{j\ge1}(-1)^{j-1} &\Bigl\{
		p\left(n-(k-\tfrac12)j^2-|m|j+\tfrac j2\right)
		{}-p\left(n-(k-\tfrac12)j^2-|m|j-\tfrac j2\right)
		\Bigr\}.
	\end{aligned}
\end{equation}
See, for example, Zhou~\cite[p.~267]{Zhou}.
Clearly by~\eqref{eq:coeff},
$A_k(-m,n)=A_k(m,n)$ for $k=1,2$.
Hence, to prove Theorem~\ref{thm:strict}, it suffices to consider $0\le m\le n-k-71$.

When $A_k(m-1,n)$, $A_k(m,n)$ and
$A_k(m+1,n)$ are positive, define
\begin{align}\label{eq:Lambda}
	\Lambda_k(m,n)=
	\log\frac{A_k(m,n)^2}{A_k(m-1,n)A_k(m+1,n)}.
\end{align}
Then strict log-concavity of the sequence $(A_k(m,n))_m$ is equivalent to $\Lambda_k(m,n)>0$. 
We first derive estimates
for the partition function and present a smooth approximation of $A_k(m,n)$.
Then we establish a criterion for the positivity of $\Lambda_k(m,n)$ in~\eqref{eq:Lambda} based on a smooth approximation.

\subsection{The smooth approximation}

\subsubsection{A smooth approximation for  $p(n)$}
For real $s>\tfrac1{24}$, define $P(s)$ by
\begin{align}\label{eq:PQ}
	P(s)=
	\frac{1}{4\sqrt3(s-\tfrac1{24})}
	\left(1-\frac{\sqrt6}{2\pi\sqrt{s-\tfrac1{24}}}\right)
	\exp\!\left(\pi\sqrt{\tfrac23}\sqrt{s-\tfrac1{24}}\right).
\end{align}
For every integer $s\ge1$, Bringmann, Man and Rolen
\cite[Lemma~2.4]{BMR} established the bound
\begin{align}\label{eq:p-error}
	|p(s)-P(s)|
	\le
	\frac{1}{8(s-\tfrac1{24})}
	\exp\!\left(\frac{\pi}{\sqrt6}\sqrt{s-\tfrac1{24}}\right).
\end{align}
For $s\geq1$, they~\cite[Eq.~(2.5)]{BMR} proved
\begin{align}\label{eq:p-upper}
	p(s)\le
	\frac{1}{4\sqrt3(s-\tfrac1{24})}
	\exp\!\left(\pi\sqrt{\tfrac23}\sqrt{s-\tfrac1{24}}\right).
\end{align}

However, the arguments in~\eqref{eq:coeff} eventually become negative, whereas $P(s)$ is defined for real $s>\tfrac1{24}$.
To overcome this difficulty, we introduce an entire function $B(s)$, as a smooth approximation to $p(s).$
For a complex variable $s$, define
\begin{align}\label{eq:Bseries}
	B(s)=
	\frac1{\sqrt3}
	\sum_{j=0}^{\infty}
	\frac{(j+1)(2\pi^2/3)^{j+1}}{(2j+3)!}
	\left(s-\tfrac1{24}\right)^j.
\end{align}
This power series converges for every complex $s$.
Thus $B(s)$ is entire, and
$B(\tfrac1{24})=\pi^2/(9\sqrt3)$.
Using the power series for $\cosh$ and $\sinh$ to sum
\eqref{eq:Bseries}, we obtain, for every complex
$s\ne\tfrac1{24}$,
\begin{align}\label{eq:B}
	B(s)=
	\frac1{2\sqrt3(s-\tfrac1{24})}
	\left(
	\cosh\!\left(\pi\sqrt{\tfrac23(s-\tfrac1{24})}\right)
	-
	\frac{\sinh\!\left(\pi\sqrt{\tfrac23(s-\tfrac1{24})}\right)}
	{\pi\sqrt{\tfrac23(s-\tfrac1{24})}}
	\right).
\end{align}
For real $s>\tfrac1{24}$, let 
$$
z=\pi\sqrt{\tfrac23(s-\tfrac1{24})}.
$$
Changing the sign of the square root replaces
$z$ by $-z$.
Since $\cosh(-z)=\cosh z$ and
$\sinh(-z)/(-z)=\sinh z/z$, the value of $B(s)$
in \eqref{eq:B} remains unchanged.

We next compare $B(s)$ with $P(s)$.
We obtain
\begin{align*}
	B(s)
	&=\frac1{2\sqrt3(s-\tfrac1{24})}
	\left(
	\frac{e^z+e^{-z}}2
	-\frac{e^z-e^{-z}}{2z}
	\right)\\
	&=\frac1{4\sqrt3(s-\tfrac1{24})}
	\left[
	\left(e^z-\frac{e^z}{z}\right)
	+\left(e^{-z}+\frac{e^{-z}}{z}\right)
	\right]\\
	&=\frac1{4\sqrt3(s-\tfrac1{24})}
	\left[
	e^z\left(1-\frac1z\right)
	+e^{-z}\left(1+\frac1z\right)
	\right]\\
	&=\frac{e^z}{4\sqrt3(s-\tfrac1{24})}
	\left(1-\frac1z\right)
	+\frac{e^{-z}}{4\sqrt3(s-\tfrac1{24})}
	\left(1+\frac1z\right).
\end{align*}
Thus, we have 
\begin{align*}
	B(s)=
	\frac{e^{\pi\sqrt{\tfrac23(s-\tfrac1{24})}}}
	{4\sqrt3(s-\tfrac1{24})}
	\left(
	1-\frac1{\pi\sqrt{\tfrac23(s-\tfrac1{24})}}
	\right)
	+
	\frac{e^{-\pi\sqrt{\tfrac23(s-\tfrac1{24})}}}
	{4\sqrt3(s-\tfrac1{24})}
	\left(
	1+\frac1{\pi\sqrt{\tfrac23(s-\tfrac1{24})}}
	\right).
\end{align*}
By \eqref{eq:PQ}, we have 
\begin{align}\label{eq:W}
	B(s)=P(s)+W(s)
	\qquad\text{for real }s>\tfrac1{24},
\end{align}
where
\begin{align}\label{eq:W-definition}
	W(s)=
	\frac{e^{-\pi\sqrt{\tfrac23(s-\tfrac1{24})}}}
	{4\sqrt3(s-\tfrac1{24})}
	\left(
	1+\frac1{\pi\sqrt{\tfrac23(s-\tfrac1{24})}}
	\right).
\end{align}

\begin{lemma}\label{lem:Bbounds}
	Let $B$ be defined by \eqref{eq:Bseries}, and let $r$ be
	a nonnegative integer.
	For real $s>0$, we have
	\begin{align}\label{eq:B-derivative-bound}
		0<
		B^{(r+1)}\!\left(\tfrac1{24}+s\right)
		\le
		\frac{\pi}{\sqrt{6s}}
		B^{(r)}\!\left(\tfrac1{24}+s\right).
	\end{align}
	For complex $s$, we have
	\begin{align}\label{eq:B-complex-bound}
		\left|B^{(r)}\!\left(\tfrac1{24}+s\right)\right|
		\le
		B^{(r)}\!\left(\tfrac1{24}+|s|\right).
	\end{align}
	For real $s\le\tfrac1{24}$, we have
	\begin{align}\label{eq:B-left-bound}
		|B(s)|<1.
	\end{align}
\end{lemma}

\begin{proof}
	For real $s>0$, we first use induction on $r$ to establish the integral
	representation
	\begin{align}\label{eq:Bint}
		B^{(r)}\!\left(\tfrac1{24}+s\right)
		=
		\frac{\pi^{2r+2}}
		{2\cdot6^{r+1}\sqrt3(r+1)!}
		\int_{-1}^{1}
		(1-u^2)^{r+1}e^{\pi\sqrt{2s/3}\,u}\dd u
	\end{align}

	For $r=0$, integration by parts gives
	\begin{align}\label{eq:B-base-first}
		\int_{-1}^{1}
		(1-u^2)e^{\pi\sqrt{2s/3}\,u}\dd u
		=
		\frac{2}{\pi\sqrt{2s/3}}
		\int_{-1}^{1}
		u e^{\pi\sqrt{2s/3}\,u}\dd u.
	\end{align}
	Again integration by parts yields
	\begin{align}\label{eq:B-base-second}
		\int_{-1}^{1}u e^{\pi\sqrt{2s/3}\,u}\dd u
		&=
		\frac{
			e^{\pi\sqrt{2s/3}}+e^{-\pi\sqrt{2s/3}}
		}{\pi\sqrt{2s/3}}
		-
		\frac1{\pi\sqrt{2s/3}}
		\int_{-1}^{1}e^{\pi\sqrt{2s/3}\,u}\dd u
		\notag\\
		&=
		\frac{
			e^{\pi\sqrt{2s/3}}+e^{-\pi\sqrt{2s/3}}
		}{\pi\sqrt{2s/3}}
		-
		\frac{
			e^{\pi\sqrt{2s/3}}-e^{-\pi\sqrt{2s/3}}
		}{\pi^2(2s/3)}\notag\\
		&=
		\frac{2\cosh\!\left(\pi\sqrt{2s/3}\right)}
		{\pi\sqrt{2s/3}}
		-
		\frac{3\sinh\!\left(\pi\sqrt{2s/3}\right)}
		{\pi^2s}.
	\end{align}
	Substituting \eqref{eq:B-base-second} into
	\eqref{eq:B-base-first} gives
	\begin{align}\label{eq:B-base-integral}
		&\int_{-1}^{1}
		(1-u^2)e^{\pi\sqrt{2s/3}\,u}\dd u=
		\frac{6}{\pi^2s}
		\left(
		\cosh\!\left(\pi\sqrt{2s/3}\right)
		-
		\frac{\sinh\!\left(\pi\sqrt{2s/3}\right)}
		{\pi\sqrt{2s/3}}
		\right).
	\end{align}
	By the definition~\eqref{eq:B} of $B(s),$ 
	multiplying \eqref{eq:B-base-integral} by
	$\pi^2/(12\sqrt3)$, we get
	\begin{align*}
		\frac{\pi^2}{12\sqrt3}
		\int_{-1}^{1}
		(1-u^2)e^{\pi\sqrt{2s/3}\,u}\dd u
		=
		B\!\left(\tfrac1{24}+s\right).
	\end{align*}
	Thus \eqref{eq:Bint} holds for $r=0$.
	
Now suppose that \eqref{eq:Bint} holds for some integer $r\ge0$.
Since the integrand and its partial derivative with respect to $s$
are continuous for $s>0$ and $-1\le u\le1$,
we may differentiate under the integral sign for every $s>0$.
We use
	\begin{align*}
		\frac{\mathrm d}{\mathrm du}(1-u^2)^{r+2}
		=
		-2(r+2)u(1-u^2)^{r+1}.
	\end{align*}
	Integration by parts therefore gives
	\begin{align}\label{eq:B-induction-parts}
		\int_{-1}^{1}
		u(1-u^2)^{r+1}e^{\pi\sqrt{2s/3}\,u}\dd u=
		\frac{\pi\sqrt{2s/3}}{2(r+2)}
		\int_{-1}^{1}
		(1-u^2)^{r+2}e^{\pi\sqrt{2s/3}\,u}\dd u,
	\end{align}
	where the boundary terms vanish at $u=-1$ and $u=1$.
	Differentiating \eqref{eq:Bint} and then applying
	\eqref{eq:B-induction-parts}, we obtain
	\begin{align*}
		B^{(r+1)}\!\left(\tfrac1{24}+s\right)
		&=
		\frac{\pi^{2r+2}}
		{2\cdot6^{r+1}\sqrt3(r+1)!}
		\frac{\pi}{\sqrt{6s}}
		\int_{-1}^{1}
		u(1-u^2)^{r+1}e^{\pi\sqrt{2s/3}\,u}\dd u\\
		&=
		\frac{\pi^{2r+2}}
		{2\cdot6^{r+1}\sqrt3(r+1)!}
		\frac{\pi^2}{6(r+2)}
		\int_{-1}^{1}
		(1-u^2)^{r+2}e^{\pi\sqrt{2s/3}\,u}\dd u\\
		&=
		\frac{\pi^{2r+4}}
		{2\cdot6^{r+2}\sqrt3(r+2)!}
		\int_{-1}^{1}
		(1-u^2)^{r+2}e^{\pi\sqrt{2s/3}\,u}\dd u.
	\end{align*}
	By induction, \eqref{eq:Bint} holds for all integers
	$r\ge0$ and all real $s>0$.
	
	We next show that~\eqref{eq:Bint} holds for all complex $s$.
	Expanding the exponential $e^{\pi\sqrt{2s/3}\,u}$ in~\eqref{eq:Bint} and integrating term by term gives
	\begin{align}\label{eq:B-integral-series}
		&\int_{-1}^{1}
		(1-u^2)^{r+1}e^{\pi\sqrt{2s/3}\,u}\dd u=
		\sum_{\ell=0}^{\infty}
		\frac{(2\pi^2s/3)^\ell}{(2\ell)!}
		\int_{-1}^{1}u^{2\ell}(1-u^2)^{r+1}\dd u.
	\end{align}
	Indeed, the exponential series converges uniformly in $u$
	on the integration interval, and every odd power of $u$
	integrates to zero.
	Each integral on the right is at most $2$, so the resulting
	power series converges uniformly on bounded sets of complex $s$.
	Both sides of \eqref{eq:Bint} are entire functions of $s$
	and agree for every real $s>0$.
	Since the positive real axis has a limit point at $s=1$
	inside their common domain, the identity theorem for
	analytic functions implies that \eqref{eq:Bint} holds
	for every complex $s$.
	
	We now prove \eqref{eq:B-derivative-bound}.
	For real $s>0$, the integral in \eqref{eq:Bint} is positive,
	so $B^{(r)}(\tfrac1{24}+s)>0$.
	Differentiating \eqref{eq:Bint} with respect to $s$, we obtain
	\begin{align*}
		B^{(r+1)}\!\left(\tfrac1{24}+s\right)
		&=
		\frac{\pi^{2r+2}}
		{2\cdot6^{r+1}\sqrt3(r+1)!}
		\int_{-1}^{1}(1-u^2)^{r+1}
		\frac{\partial}{\partial s}
		e^{\pi\sqrt{2s/3}\,u}\dd u\\
		&=
		\frac{\pi^{2r+2}}
		{2\cdot6^{r+1}\sqrt3(r+1)!}
		\frac{\pi}{\sqrt{6s}}
		\int_{-1}^{1}
		u(1-u^2)^{r+1}e^{\pi\sqrt{2s/3}\,u}\dd u.
	\end{align*}
	Dividing by $B^{(r)}(\tfrac1{24}+s)>0$ and using
	\eqref{eq:Bint} gives
	\begin{align}\label{eq:B-derivative-ratio}
		\frac{B^{(r+1)}(\tfrac1{24}+s)}
		{B^{(r)}(\tfrac1{24}+s)}
		&=
		\frac{\pi}{\sqrt{6s}}
		\frac{
			\displaystyle\int_{-1}^{1}
			u(1-u^2)^{r+1}e^{\pi\sqrt{2s/3}\,u}\dd u
		}{
			\displaystyle\int_{-1}^{1}
			(1-u^2)^{r+1}e^{\pi\sqrt{2s/3}\,u}\dd u
		}.
	\end{align}

Note that the numerator in~\eqref{eq:B-derivative-ratio}
\begin{align}
\int_{-1}^{1}
u(1-u^2)^{r+1}e^{\pi\sqrt{2s/3}\,u}\dd u=&\int_0^1u(1-u^2)^{r+1}e^{au}\,du
	-\int_0^1u(1-u^2)^{r+1}e^{-au}\,du\nonumber\\
	&=\int_0^1u(1-u^2)^{r+1}(e^{au}-e^{-au})\,du\nonumber\\
	&=2\int_0^1u(1-u^2)^{r+1}\sinh(au)\,du>0.
\end{align}
	Moreover, since $u<1$, the numerator is strictly smaller
	than its denominator.
Hence, we prove \eqref{eq:B-derivative-bound}.
	
	We next prove \eqref{eq:B-complex-bound}.
	Differentiating \eqref{eq:Bseries} $r$ times gives
	\begin{align}\label{eq:B-derivative-series}
		B^{(r)}\!\left(\tfrac1{24}+s\right)
		=
		\frac1{\sqrt3}
		\sum_{j=0}^{\infty}
		\frac{(j+r+1)!(2\pi^2/3)^{j+r+1}}
		{j!(2j+2r+3)!}\,s^j.
	\end{align}
	All coefficients in this series are positive.
	Therefore, 
	\begin{align*}
		\left|B^{(r)}\!\left(\tfrac1{24}+s\right)\right|
		&\le
		\frac1{\sqrt3}
		\sum_{j=0}^{\infty}
		\frac{(j+r+1)!(2\pi^2/3)^{j+r+1}}
		{j!(2j+2r+3)!}\,|s|^j
		=B^{(r)}\!\left(\tfrac1{24}+|s|\right),
	\end{align*}
	which proves \eqref{eq:B-complex-bound}.
	
	Finally, let $s\le 1/24$.
	By \eqref{eq:Bint} with $r=0$, we have
	\[
	B(s)=\frac{\pi^2}{12\sqrt3}
	\int_{-1}^{1}(1-u^2)
	e^{\pi\sqrt{\tfrac23(s-\tfrac1{24})}\,u}\dd u.
	\]
	Since $s\le 1/{24}$, the square root is purely imaginary or zero.
	Thus,
	\begin{align*}
		|B(s)|
		\le
		\frac{\pi^2}{12\sqrt3}
		\int_{-1}^{1}(1-u^2)\dd u
		=\frac{\pi^2}{9\sqrt3}<1.
	\end{align*}
	This proves \eqref{eq:B-left-bound}.
\end{proof}

\subsubsection{A smooth approximation for $A_k(m,n)$}
As in \eqref{eq:coeff},
for $k\in\{1,2\}$, $n\ge2$ and $m\ge0$,  we define
\begin{equation}\label{eq:model}
	\begin{aligned}
		\A_k(m,n)
		=\sum_{j\ge1}(-1)^{j-1} &\Bigl\{
		B\left(n-(k-\tfrac12)j^2-mj+\tfrac j2\right)
		-B\left(n-(k-\tfrac12)j^2-mj-\tfrac j2\right)
		\Bigr\}.
	\end{aligned}
\end{equation}

For real $s<\tfrac1{24}$, let
\[
\widetilde z=\pi\sqrt{\tfrac23(\tfrac1{24}-s)}>0.
\]
Since $\cosh(iz)=\cos z$ and $\sinh(iz)=i\sin z$,
by~\eqref{eq:B}, we have 
$$
B(s)=\frac1{2\sqrt3(\tfrac1{24}-s)}
\left(\frac{\sin \widetilde z}{\widetilde z}-\cos \widetilde z\right).
$$
Using $|\sin \widetilde z|\le1$ and $|\cos \widetilde z|\le1$, we obtain
\[
|B(s)|
\le
\frac1{2\sqrt3(\tfrac1{24}-s)}
\left(1+\frac1 {\widetilde z}\right).
\]
As real $s\to-\infty$, we have $\widetilde z\to\infty$ and
$(\tfrac1{24}-s)/|s|\to1$. Therefore,
\begin{equation}\label{eq:B-negative-decay}
	B(s)=O(|s|^{-1})
	\qquad\text{as real }s\to-\infty.
\end{equation}

Fix $k\in\{1,2\}$, $n>\tfrac1{24}$, and $m\in\mathbb R$,
\[
\frac{n-\left(k-\tfrac12\right)j^2-mj\pm\tfrac j2}{j^2}
=-\left(k-\tfrac12\right)
+\frac{-m\pm\tfrac12}{j}
+\frac{n}{j^2}
\longrightarrow
-\left(k-\tfrac12\right)<0 \qquad\text{as }j\to\infty.
\]
Consequently,  for all sufficiently
large $j$,
\[
\left|n-\left(k-\tfrac12\right)j^2-mj\pm\tfrac j2\right|\ge \frac12\left(k-\tfrac12\right)j^2.
\]
By~\eqref{eq:B-negative-decay}, we obtain the $j$-th summand in~\eqref{eq:model}
\begin{align*}
\left|(-1)^{j-1}\Bigl\{
B\left(n-(k-\tfrac12)j^2-mj+\tfrac j2\right)
-B\left(n-(k-\tfrac12)j^2-mj-\tfrac j2\right)
\Bigr\}\right|
=O(j^{-2}).
\end{align*}
Since $\sum_{j\ge1}j^{-2}$ converges, the series in \eqref{eq:model} converges
absolutely for every real $m$.

\begin{lemma}\label{lem:integral}
	Let $k\in\{1,2\}$, and let $n>\tfrac1{24}$ be real.
	For every real $m$, the function $\A_k(m,n)$ defined
	by \eqref{eq:model} satisfies
	\begin{align}\label{eq:Aint}
		\A_k(m,n)
		=
		\int_0^\infty\frac{t}{\sinh(\pi t)}
		\int_{-1/2}^{1/2}
		\Re B'\!\left(
		n+\left(k-\tfrac12\right)t^2-i(m+u)t
		\right)\dd u\dd t.
	\end{align}
\end{lemma}

\begin{proof}
	For complex $z$, define
	\begin{align}\label{eq:Aint-g-definition}
		g(z)
		={}&
		B\!\left(
		n-\left(k-\tfrac12\right)z^2
		-\left(m-\tfrac12\right)z
		\right)-
		B\!\left(
		n-\left(k-\tfrac12\right)z^2
		-\left(m+\tfrac12\right)z
		\right).
	\end{align}
	We prove \eqref{eq:Aint} by applying the residue theorem
	to $g(z).$
	
	Since $B$ is entire by \eqref{eq:Bseries}, $g$ is entire.
	The coefficients in \eqref{eq:Bseries} are real, and
	$k,n,m$ are real. Hence
	\begin{align}\label{eq:Aint-g-conjugation}
		g(\overline z)=\overline{g(z)}.
	\end{align}
	Moreover, $g(0)=0$, so the singularity of
	$g(z)/\sin(\pi z)$ at zero is removable.
	
	We first estimate $g(z)$ for large $|z|$.
	For either $b=m-\tfrac12$ or $b=m+\tfrac12$, choose
	the square root in the second factor below to tend
	to $1$ as $|z|\to\infty$. Then
	\begin{align}\label{eq:Aint-root}
		\sqrt{
			n-\tfrac1{24}-\left(k-\tfrac12\right)z^2-bz
		}&=
		i\sqrt{k-\tfrac12}\,z
		\sqrt{
			1+\frac{b}{(k-\tfrac12)z}
			-\frac{n-\tfrac1{24}}{(k-\tfrac12)z^2}
		}
		\notag\\
		&=
		i\sqrt{k-\tfrac12}\,z+O(1).
	\end{align}
	This choice is allowed because \eqref{eq:B} is unchanged
	when the sign of its square root is reversed.
	It follows from \eqref{eq:Aint-root} that
	\begin{align*}
		\left|
		\Re\sqrt{
			n-\tfrac1{24}-\left(k-\tfrac12\right)z^2-bz
		}
		\right|
		\le
		\sqrt{k-\tfrac12}\,|\Im z|+O(1).
	\end{align*}
	Also, the triangle inequality gives
	\begin{align*}
		\left|
		n-\tfrac1{24}-\left(k-\tfrac12\right)z^2-bz
		\right|
		\ge
		\left(k-\tfrac12\right)|z|^2
		-|b||z|-\left(n-\tfrac1{24}\right)
		\ge
		\frac{k-\tfrac12}{2}|z|^2
	\end{align*}
	for sufficiently large $|z|$.
	Using these estimates in \eqref{eq:B}, together with
	the exponential expressions for $\cosh$ and $\sinh$,
	we obtain
	\begin{align}\label{eq:Aint-g-growth}
		g(z)
		=
		O\!\left(
		|z|^{-2}
		e^{\pi\sqrt{(2k-1)/3}\,|\Im z|}
		\right)
		\qquad\text{as }|z|\to\infty.
	\end{align}
	
	Let $R=N+\tfrac12$, where $N$ is a positive integer,
	and let $\Gamma_R$ be the positively oriented rectangle
	with vertices
	\begin{align*}
		-iR,\qquad R-iR,\qquad R+iR,\qquad iR.
	\end{align*}
	Its interior contains the integers $1,\ldots,N$.
	At each such integer $j$, we have
	\begin{align*}
		\operatorname{Res}_{z=j}
		\frac{g(z)}{\sin(\pi z)}
		=
		\frac{g(j)}{\pi\cos(\pi j)}
		=
		\frac{(-1)^jg(j)}{\pi}.
	\end{align*}
	Therefore, the residue theorem gives
	\begin{align}\label{eq:Aint-residues}
		\int_{\Gamma_R}\frac{g(z)}{\sin(\pi z)}\dd z
		=
		2i\sum_{j=1}^{N}(-1)^jg(j).
	\end{align}
	
	We next show that the integrals over the right,
	upper, and lower sides tend to zero.
	On the right side, $z=R+iy$ and
	\begin{align*}
		|\sin(\pi z)|=\cosh(\pi y).
	\end{align*}
	On the upper and lower sides, $z=x\pm iR$, and
	\begin{align*}
		|\sin(\pi z)|^2
		=
		\sin^2(\pi x)+\sinh^2(\pi R)
		\ge\sinh^2(\pi R).
	\end{align*}
	Since $|z|\ge R$ on these three sides,
	\eqref{eq:Aint-g-growth} implies
	\begin{align*}
		\frac{g(z)}{\sin(\pi z)}
		=
		O\!\left(
		R^{-2}
		e^{-\pi(1-\sqrt{(2k-1)/3})|\Im z|}
		\right)
		=
		O(R^{-2}).
	\end{align*}
	Here we used $\sqrt{(2k-1)/3}\le1$ for $k=1,2$.
	Their total length is $4R$, so their total contribution
	is $O(R^{-1})$.
	In particular, this contribution also tends to zero
	when $k=2$.
	
	On the left side of $\Gamma_R$, the path runs from
	$iR$ to $-iR$.
	Substituting $z=it$ and using
	$\sin(\pi it)=i\sinh(\pi t)$, we obtain
	\begin{align}\label{eq:Aint-left-side}
		\int_{iR}^{-iR}\frac{g(z)}{\sin(\pi z)}\dd z
		=
		-\int_{-R}^{R}\frac{g(it)}{\sinh(\pi t)}\dd t
		=
		-\int_0^R
		\frac{g(it)-g(-it)}{\sinh(\pi t)}\dd t
		=
		-2i\int_0^R
		\frac{\Im g(it)}{\sinh(\pi t)}\dd t,
	\end{align}
	where the last equality follows from
	\eqref{eq:Aint-g-conjugation}.
	
	The integral in \eqref{eq:Aint-left-side} converges
	absolutely as $R\to\infty$.
	Indeed, $g(0)=0$ gives $g(it)=O(t)$ near zero.
	At infinity, \eqref{eq:Aint-g-growth} gives
	\begin{align*}
		\frac{|g(it)|}{\sinh(\pi t)}=O(t^{-2}),
	\end{align*}
	because $\sqrt{(2k-1)/3}\le1$.
	Combining \eqref{eq:Aint-residues} and
	\eqref{eq:Aint-left-side}, and letting $N\to\infty$,
	we therefore obtain
	\begin{align}\label{eq:Aint-g-integral}
		\A_k(m,n)
		=
		\sum_{j=1}^{\infty}(-1)^{j-1}g(j)
		=
		\int_0^\infty
		\frac{\Im g(it)}{\sinh(\pi t)}\dd t.
	\end{align}
	The first equality follows from \eqref{eq:model}
	and \eqref{eq:Aint-g-definition}.
	
	Finally, the fundamental theorem of calculus gives
	\begin{align}\label{eq:Aint-g-difference}
		g(it)
		&=
		-\int_{-1/2}^{1/2}
		\frac{\mathrm d}{\mathrm du}
		B\!\left(
		n+\left(k-\tfrac12\right)t^2-i(m+u)t
		\right)\dd u
		\notag\\
		&=
		it\int_{-1/2}^{1/2}
		B'\!\left(
		n+\left(k-\tfrac12\right)t^2-i(m+u)t
		\right)\dd u.
	\end{align}
	Taking imaginary parts yields
	\begin{align*}
		\Im g(it)
		=
		t\int_{-1/2}^{1/2}
		\Re B'\!\left(
		n+\left(k-\tfrac12\right)t^2-i(m+u)t
		\right)\dd u.
	\end{align*}
	Substitution into \eqref{eq:Aint-g-integral}
	proves \eqref{eq:Aint}.
\end{proof}

The integral representation also shows that $\A_k(m,n)$
is even in $m$.
Indeed, the real coefficients in \eqref{eq:Bseries} give
$B'(\overline s)=\overline{B'(s)}$.
Replacing $m$ by $-m$ in \eqref{eq:Aint} and then
replacing $u$ by $-u$, the expression inside $B'$
becomes the complex conjugate of the original expression.
Its real part is unchanged, so
\begin{align}\label{eq:A-even}
	\A_k(-m,n)=\A_k(m,n)
	\qquad\text{for every real }m.
\end{align}

We next justify differentiation of \eqref{eq:Aint}
with respect to $m$.
For $r=0,1,2$ and complex $s\ne\tfrac1{24}$,
\eqref{eq:B-complex-bound} and repeated applications
of \eqref{eq:B-derivative-bound} give
\begin{align}\label{eq:Aint-B-growth}
	|B^{(r+1)}(s)|
	&\le
	\left(
	\frac{\pi}{\sqrt{6|s-\tfrac1{24}|}}
	\right)^{r+1}
	B\!\left(\tfrac1{24}+|s-\tfrac1{24}|\right)
	\notag\\
	&\le
	\frac{
		\pi^{r+1}
		e^{\pi\sqrt{2/3}\sqrt{|s-\tfrac1{24}|}}
	}{
		2\sqrt3\,6^{(r+1)/2}
		|s-\tfrac1{24}|^{(r+3)/2}
	}.
\end{align}
For the second inequality, we used \eqref{eq:B}:
for real $v>0$,
\begin{align*}
	0<B\!\left(\tfrac1{24}+v\right)
	\le
	\frac{e^{\pi\sqrt{2v/3}}}{2\sqrt3\,v}.
\end{align*}

For $m$ in any bounded real interval and
$-\tfrac12\le u\le\tfrac12$, we have, uniformly as
$t\to\infty$,
\begin{align*}
	&\left|
	n-\tfrac1{24}+\left(k-\tfrac12\right)t^2-i(m+u)t
	\right|\\
	&\qquad=
	\sqrt{
		\left(n-\tfrac1{24}+\left(k-\tfrac12\right)t^2\right)^2
		+(m+u)^2t^2
	}
	=
	\left(k-\tfrac12\right)t^2+O(1).
\end{align*}
Consequently, \eqref{eq:Aint-B-growth} implies
\begin{align}\label{eq:Aint-derivative-tail}
	&\frac{t^{r+1}}{\sinh(\pi t)}
	\left|
	B^{(r+1)}\!\left(
	n+\left(k-\tfrac12\right)t^2-i(m+u)t
	\right)
	\right|=
	O\!\left(
	t^{-2}
	e^{-\pi(1-\sqrt{(2k-1)/3})t}
	\right)
	=
	O(t^{-2}).
\end{align}
Near $t=0$, the derivatives of $B$ remain bounded
uniformly in these ranges of $m$ and $u$.
Since $t^{r+1}/\sinh(\pi t)=O(t^r)$, the expressions
on the left of \eqref{eq:Aint-derivative-tail}
are integrable near zero as well.

These bounds provide integrable upper bounds independent
of $m$ on each bounded real interval.
Thus the dominated convergence theorem justifies two
differentiations of \eqref{eq:Aint} under the integral sign
and shows that both derivatives are continuous.
By the chain rule, each differentiation contributes
a factor $-it$.
Using $\Re(-iz)=\Im z$ and $(-i)^2=-1$, we obtain
\begin{align}\label{eq:Aint-first-derivative}
	\frac{\partial}{\partial m}\A_k(m,n)
	=
	\int_0^\infty\frac{t^2}{\sinh(\pi t)}
	\int_{-1/2}^{1/2}
	\Im B''\!\left(
	n+\left(k-\tfrac12\right)t^2-i(m+u)t
	\right)\dd u\dd t
\end{align}
and
\begin{align}\label{eq:Aint-second-derivative}
	\frac{\partial^2}{\partial m^2}\A_k(m,n)
	=
	-\int_0^\infty\frac{t^3}{\sinh(\pi t)}
	\int_{-1/2}^{1/2}
	\Re B'''\!\left(
	n+\left(k-\tfrac12\right)t^2-i(m+u)t
	\right)\dd u\dd t.
\end{align}
All three integrals in \eqref{eq:Aint},
\eqref{eq:Aint-first-derivative}, and
\eqref{eq:Aint-second-derivative} converge absolutely.

\begin{lemma}\label{lem:integer-model}
	Let $A_k(m,n)$ be defined by \eqref{eq:Ak}, and let
	$\A_k(m,n)$ be defined by \eqref{eq:model}.
	For $k\in\{1,2\}$ and integers $n\ge5000$ and $m\ge0$,
	we have
	\begin{align}\label{eq:modelerror}
		|A_k(m,n)-\A_k(m,n)|
		<
		\frac2{5\sqrt n}
		\exp\!\left(
		\pi\sqrt{\frac{n-\tfrac1{24}}6}
		\right).
	\end{align}
\end{lemma}

\begin{proof}
	We first obtain a bound for $|p(s)-B(s)|$ at integer
	values of $s$.
	For real $s>\tfrac1{24}$, direct differentiation gives
	\begin{align}\label{eq:integer-envelope-derivative}
		\frac{\mathrm d}{\mathrm ds}
		\left(
		\frac{\exp\!\left(\pi\sqrt{(s-\tfrac1{24})/6}\right)}
		{8(s-\tfrac1{24})}
		\right)=
		\frac{\exp\!\left(\pi\sqrt{(s-\tfrac1{24})/6}\right)}
		{8(s-\tfrac1{24})^2}
		\left(
		\frac{\pi}{2}\sqrt{\frac{s-\tfrac1{24}}6}-1
		\right).
	\end{align}
	The factor outside the parentheses is positive.
	The expression inside the parentheses increases with $s$,
	so the derivative can change sign only from negative
	to positive.
	Hence the maximum of the expression being differentiated
	on $1\le s\le n$ occurs at one of the endpoints.
	
	At $s=1$, its value satisfies
	\begin{align*}
		\frac3{23}\exp\!\left(\frac{\pi\sqrt{23}}{12}\right)<2.
	\end{align*}
	Comparing the two endpoint values therefore gives
	\begin{align}\label{eq:integer-envelope-bound}
		\frac{\exp\!\left(\pi\sqrt{(s-\tfrac1{24})/6}\right)}
		{8(s-\tfrac1{24})}
		<
		\frac{\exp\!\left(\pi\sqrt{(n-\tfrac1{24})/6}\right)}
		{8(n-\tfrac1{24})}
		+2
	\end{align}
	for $1\le s\le n$.
	
	For real $s\ge1$, \eqref{eq:W-definition} also gives
	\begin{align*}
		0<W(s)
		<
		\frac1{2\sqrt3(s-\tfrac1{24})}
		<1.
	\end{align*}
	Indeed, the exponential factor is less than $1$,
	and the expression inside its parentheses is less than $2$.
	By \eqref{eq:W}, we have $B(s)=P(s)+W(s)$.
	Combining \eqref{eq:p-error} and
	\eqref{eq:integer-envelope-bound}, we obtain
	\begin{align*}
		|p(s)-B(s)|
		\le |p(s)-P(s)|+W(s)
		<
		\frac{\exp\!\left(\pi\sqrt{(n-\tfrac1{24})/6}\right)}
		{8(n-\tfrac1{24})}
		+3
	\end{align*}
	for every integer $1\le s\le n$.
	
	At $s=0$, the convention $p(0)=1$ and
	\eqref{eq:B-left-bound} give $|p(0)-B(0)|<2$.
	For negative integers $s$, we have $p(s)=0$, so
	\begin{align}\label{eq:integer-negative-bound}
		|p(s)-B(s)|=|B(s)|<1.
	\end{align}
	Thus, for every integer $s\le n$, we have
	\begin{align}\label{eq:integer-single-bound}
		|p(s)-B(s)|
		<
		\frac{\exp\!\left(\pi\sqrt{(n-\tfrac1{24})/6}\right)}
		{8(n-\tfrac1{24})}
		+3.
	\end{align}
	
	We now apply these bounds to the sums defining
	$A_k(m,n)$ and $\A_k(m,n)$.
	For each positive integer $j$, write
	\begin{align}\label{eq:integer-sj}
		s_j=
		n-\left(k-\tfrac12\right)j^2-mj+\frac j2.
	\end{align}
	Both $s_j$ and $s_j-j$ are integers.
	Since $k\ge1$ and $m\ge0$, we have
	\begin{align}\label{eq:integer-sj-upper}
		s_j\le n-\frac{j(j-1)}2\le n.
	\end{align}
	Subtracting \eqref{eq:model} from \eqref{eq:coeff} gives
	\begin{align}\label{eq:integer-difference}
		A_k(m,n)-\A_k(m,n)
		=
		\sum_{j=1}^{\infty}(-1)^{j-1}
		\Bigl(
		p(s_j)-B(s_j)-p(s_j-j)+B(s_j-j)
		\Bigr).
	\end{align}
	Applying the triangle inequality twice to
	\eqref{eq:integer-difference}, we obtain
	\begin{align}\label{eq:integer-sum-bound}
		|A_k(m,n)-\A_k(m,n)|
		&\le
		\sum_{j=1}^{\infty}
		\left|
		p(s_j)-B(s_j)-p(s_j-j)+B(s_j-j)
		\right|
		\notag\\
		&\le
		\sum_{j=1}^{\infty}
		\Bigl(
		|p(s_j)-B(s_j)|
		+
		|p(s_j-j)-B(s_j-j)|
		\Bigr).
	\end{align}
	We now estimate the right-hand side of
	\eqref{eq:integer-sum-bound} in three parts.
	
	First, consider the indices for which $s_j\ge0$.
	If $j\ge3\sqrt n/2$, then
	\eqref{eq:integer-sj-upper} gives
	\begin{align*}
		s_j
		\le
		n-\frac12\frac{3\sqrt n}{2}
		\left(\frac{3\sqrt n}{2}-1\right)
		=
		-\frac n8+\frac{3\sqrt n}{4}
		<0,
	\end{align*}
	since $n\ge5000$.
	Therefore $s_j\ge0$ implies $j<3\sqrt n/2$.
	There are fewer than $3\sqrt n/2$ such indices.
	For every positive integer $j$ with $s_j\ge0$, we have
	$s_j\le n$ and $s_j-j\le n$ by
	\eqref{eq:integer-sj-upper}.
	Applying \eqref{eq:integer-single-bound} at $s=s_j$
	and $s=s_j-j$, we obtain
	\begin{align*}
		|p(s_j)-B(s_j)|+|p(s_j-j)-B(s_j-j)|
		<
		2\left(
		\frac{\exp\!\left(\pi\sqrt{(n-\tfrac1{24})/6}\right)}
		{8(n-\tfrac1{24})}
		+3
		\right).
	\end{align*}
	As shown above, $s_j\ge0$ implies $j<3\sqrt n/2$,
	so there are fewer than $3\sqrt n/2$ such indices.
	Therefore, the terms in \eqref{eq:integer-sum-bound}
	with $s_j\ge0$ satisfy
	\begin{align}\label{eq:integer-positive-contribution}
		\sum_{\substack{j\ge1\\s_j\ge0}}
		\Bigl(
		|p(s_j)-B(s_j)|
		+
		|p(s_j-j)-B(s_j-j)|
		\Bigr)
		<
		3\sqrt n
		\left(
		\frac{\exp\!\left(\pi\sqrt{(n-\tfrac1{24})/6}\right)}
		{8(n-\tfrac1{24})}
		+3
		\right).
	\end{align}
	
	Next, consider the indices with $s_j<0$ and
	$j\le4\sqrt n$.
	Then $s_j-j<0$ as well.
	By \eqref{eq:integer-negative-bound}, each summand
	in \eqref{eq:integer-sum-bound} is less than $2$.
	There are at most $4\sqrt n$ such indices,
	so their total contribution is less than $8\sqrt n$.
	
	It remains to consider $j>4\sqrt n$.
	In this case $n<j^2/16$ and $j>4$, so
	\eqref{eq:integer-sj-upper} gives
	\begin{align}\label{eq:integer-far-negative}
		s_j
		<
		-\frac{7j^2}{16}+\frac j2
		\le
		-\frac{5j^2}{16}
		<
		-\frac{j^2}{4}.
	\end{align}
	The same upper bound holds for $s_j-j$.
	
	For real $s<-j^2/4$, the hyperbolic functions in
	\eqref{eq:B} become trigonometric functions.
	Using $|\cos x|\le1$ and $|\sin x|\le1$, we obtain
	\begin{align}\label{eq:integer-B-tail}
		|B(s)|
		\le
		\frac1{2\sqrt3(\tfrac1{24}-s)}
		\left(
		1+
		\frac1{\pi\sqrt{\tfrac23(\tfrac1{24}-s)}}
		\right)
		<
		\frac{2}{\sqrt3j^2}
		\left(1+\frac{\sqrt6}{\pi j}\right)
		<
		\frac3{j^2}.
	\end{align}
	The last inequality follows from $j>4$ and $\pi>3$.
	Since $p(s_j)=p(s_j-j)=0$ in this range,
	\eqref{eq:integer-B-tail} shows that its contribution
	is less than
	\begin{align}\label{eq:integer-tail-contribution}
		6\sum_{j>4\sqrt n}\frac1{j^2}
		\le
		6\int_{4\sqrt n-1}^{\infty}\frac{\dd x}{x^2}
		=
		\frac6{4\sqrt n-1}
		<
		\frac2{\sqrt n}.
	\end{align}
	
	Combining the three contributions, we obtain
	\begin{align}\label{eq:integer-total-bound}
		|A_k(m,n)-\A_k(m,n)|
		&<
		\frac{3\sqrt n}{8(n-\tfrac1{24})}
		\exp\!\left(
		\pi\sqrt{\frac{n-\tfrac1{24}}6}
		\right)
		+17\sqrt n+\frac2{\sqrt n}
		\notag\\
		&<
		\frac{3\sqrt n}{8(n-\tfrac1{24})}
		\exp\!\left(
		\pi\sqrt{\frac{n-\tfrac1{24}}6}
		\right)
		+18\sqrt n.
	\end{align}
	
	Finally, we compare \eqref{eq:integer-total-bound}
	with the right-hand side of \eqref{eq:modelerror}.
	Since $n\ge5000$, we have
	\begin{align*}
		\pi\sqrt{\frac{n-\tfrac1{24}}6}>\sqrt n.
	\end{align*}
	The exponential series has positive terms, so
	\begin{align}\label{eq:integer-exp-lower}
		e^{\sqrt n}>
		\frac{(\sqrt n)^6}{6!}
		=
		\frac{n^3}{720}.
	\end{align}
	Consequently,
	\begin{align}\label{eq:integer-small-term}
		18n\exp\!\left(
		-\pi\sqrt{\frac{n-\tfrac1{24}}6}
		\right)
		<
		18n e^{-\sqrt n}
		<
		\frac{12960}{n^2}
		<
		\frac1{50}.
	\end{align}
	Also, direct rearrangement gives
	\begin{align*}
		\frac{3n}{8(n-\tfrac1{24})}<\frac{19}{50}
		\qquad\text{for }n\ge5000.
	\end{align*}
	Multiplying \eqref{eq:integer-total-bound} by
	$\sqrt n\exp\!\left(-\pi\sqrt{(n-\tfrac1{24})/6}\right)$
	and using \eqref{eq:integer-small-term}, we obtain
	\begin{align*}
		\sqrt n\exp\!\left(
		-\pi\sqrt{\frac{n-\tfrac1{24}}6}
		\right)
		|A_k(m,n)-\A_k(m,n)|
		&<
		\frac{3n}{8(n-\tfrac1{24})}
		+
		18n\exp\!\left(
		-\pi\sqrt{\frac{n-\tfrac1{24}}6}
		\right)\\
		&<
		\frac{19}{50}+\frac1{50}
		=
		\frac25.
	\end{align*}
	Rearranging proves \eqref{eq:modelerror}.
\end{proof}

\subsection{The relative error}

\begin{lemma}\label{lem:kernel}
	For a fixed real number $m$, let $G(x)$ be a positive function
	defined for $m-1\le x\le m+1$. Assume that its first and second
	derivatives with respect to $x$ are continuous on this interval.
	Then
	\begin{align}\label{eq:kernel}
		\log\frac{G(m)^2}{G(m-1)G(m+1)}
		=-\int_{-1}^{1}(1-|t|)(\log G)''(m+t)\dd t.
	\end{align}
\end{lemma}

\begin{proof}
	Since $G(x)>0$, the function $\log G(x)$ also has continuous
	first and second derivatives. By Taylor's theorem with
	integral remainder,
	\begin{align}
		\log G(m+1)
		&=\log G(m)+(\log G)'(m)
		+\int_0^1(1-t)(\log G)''(m+t)\dd t,
		\label{eq:G-taylor-plus}\\
		\log G(m-1)
		&=\log G(m)-(\log G)'(m)
		+\int_0^1(1-t)(\log G)''(m-t)\dd t.
		\label{eq:G-taylor-minus}
	\end{align}
Adding \eqref{eq:G-taylor-plus} and \eqref{eq:G-taylor-minus}, we obtain
\begin{align*}
	\log\frac{G(m)^2}{G(m-1)G(m+1)}
	&=-\int_0^1(1-t)(\log G)''(m+t)\dd t
	-\int_0^1(1-t)(\log G)''(m-t)\dd t\\
	&=-\int_{-1}^0(1+t)(\log G)''(m+t)\dd t
	-\int_0^1(1-t)(\log G)''(m+t)\dd t\\
	&=-\int_{-1}^1(1-|t|)(\log G)''(m+t)\dd t,
\end{align*}
where the second equality follows by replacing $t$ with $-t$
in the integral involving $m-t$, and the last equality follows
from $1-|t|=1+t$ on $[-1,0]$ and $1-|t|=1-t$ on $[0,1]$.
This proves \eqref{eq:kernel}.
\end{proof}

\begin{lemma}\label{lem:transfer}
	Let $0\le\eta<1$, and let $x_j>0$ and $b_j\in\mathbb R$
	for $j=-1,0,1$. Suppose that
	$
		\left|\frac{b_j}{x_j}-1\right|\le\eta.
	$
	Then $b_j>0$ for $j=-1,0,1$, and
	\begin{align}\label{eq:transfer}
		\left|\log\frac{b_0^2}{b_{-1}b_1}
		-\log\frac{x_0^2}{x_{-1}x_1}\right|
		\le\frac{4\eta}{1-\eta}.
	\end{align}
\end{lemma}

\begin{proof}
	By the assumption, we have 
	$
		0<1-\eta\le\frac{b_j}{x_j}\le1+\eta.
	$
	Since $x_j>0$, it follows that $b_j>0$.
	
	For each $j$ such that $b_j\ne x_j$, the mean value theorem
	gives a number $\xi_j$ between $1$ and $b_j/x_j$ such that
	$
		\log\frac{b_j}{x_j}-\log1
		=\frac{b_j/x_j-1}{\xi_j}.
	$
	Both $1$ and $b_j/x_j$ are at least $1-\eta$, so
	$\xi_j\ge1-\eta>0$. Therefore,
	\begin{align*}
		\left|\log\frac{b_j}{x_j}\right|
		=\frac{|b_j/x_j-1|}{\xi_j}
		\le\frac{|b_j/x_j-1|}{1-\eta}
		\le\frac{\eta}{1-\eta}.
	\end{align*}
	The same bound holds when $b_j=x_j$, since the left-hand
	side is zero.
	
	Finally, applying the triangle
	inequality, we obtain
	\begin{align*}
		\left|\log\frac{b_0^2}{b_{-1}b_1}
		-\log\frac{x_0^2}{x_{-1}x_1}\right|
		&=\left|
		2\log\frac{b_0}{x_0}
		-\log\frac{b_{-1}}{x_{-1}}
		-\log\frac{b_1}{x_1}\right|\\
		&\le
		2\left|\log\frac{b_0}{x_0}\right|
		+\left|\log\frac{b_{-1}}{x_{-1}}\right|
		+\left|\log\frac{b_1}{x_1}\right|\\
		&\le\frac{4\eta}{1-\eta}.
	\end{align*}
	This proves~\eqref{eq:transfer}.
\end{proof}

\section{The case $n\ge5000$ and $0\le m\le2\sqrt n$}
\label{sec:central}
We use \eqref{eq:Aint} to estimate $\A_k(m,n)$ and
its first two derivatives with respect to $m$.

\begin{lemma}\label{lem:P-prime}
	Let $P(s)$ be defined by \eqref{eq:PQ}, and define
	the polynomial
	$	H_0(v)=1-3v+3v^2.$
	For every real $s>\tfrac1{24}$, we have
	\begin{align}\label{eq:Pprime-form}
		P'(s)
		=
		\frac{\pi}{12\sqrt2(s-\tfrac1{24})^{3/2}}
		\exp\!\left(
		\pi\sqrt{\tfrac23(s-\tfrac1{24})}
		\right)
		H_0\!\left(
		\frac{\sqrt6}{2\pi\sqrt{s-\tfrac1{24}}}
		\right)
		>0.
	\end{align}
	For every real $s$ satisfying $s-\tfrac1{24}\ge120$,
	we also have
	\begin{align}\label{eq:Pprime-lower}
		P'(s)>
		\frac1{20s^{3/2}}
		\exp\!\left(
		\pi\sqrt{\tfrac23(s-\tfrac1{24})}
		\right).
	\end{align}
\end{lemma}

\begin{proof}
	We first prove the identity in \eqref{eq:Pprime-form}.
	Rewriting \eqref{eq:PQ}, we have
	\begin{align}\label{eq:P-product}
		P(s)=
		\frac1{4\sqrt3}
		\left(
		\frac1{s-\tfrac1{24}}
		-
		\frac{\sqrt6}
		{2\pi(s-\tfrac1{24})^{3/2}}
		\right)
		\exp\!\left(
		\pi\sqrt{\tfrac23(s-\tfrac1{24})}
		\right).
	\end{align}
	By the chain rule, we get 
	\begin{align}\label{eq:P-exponential-derivative}
		\frac{\mathrm d}{\mathrm ds}
		\exp\!\left(
		\pi\sqrt{\tfrac23(s-\tfrac1{24})}
		\right)
		=
		\frac{\pi}{\sqrt{6(s-\tfrac1{24})}}
		\exp\!\left(
		\pi\sqrt{\tfrac23(s-\tfrac1{24})}
		\right).
	\end{align}
	Applying the product rule to \eqref{eq:P-product}
	and using \eqref{eq:P-exponential-derivative}, we obtain
	\begin{align}\label{eq:Pprime-expanded}
		P'(s)
		={}&
		\frac{\exp\!\left(
			\pi\sqrt{\tfrac23(s-\tfrac1{24})}
			\right)}{4\sqrt3}
		\left(
		-\frac1{(s-\tfrac1{24})^2}
		+
		\frac{3\sqrt6}
		{4\pi(s-\tfrac1{24})^{5/2}}
		\right)
		\notag\\
		&+
		\frac{\pi\exp\!\left(
			\pi\sqrt{\tfrac23(s-\tfrac1{24})}
			\right)}
		{4\sqrt3\sqrt{6(s-\tfrac1{24})}}
		\left(
		\frac1{s-\tfrac1{24}}
		-
		\frac{\sqrt6}
		{2\pi(s-\tfrac1{24})^{3/2}}
		\right)
		\notag\\
		={}&
		\frac{\exp\!\left(
			\pi\sqrt{\tfrac23(s-\tfrac1{24})}
			\right)}{4\sqrt3}
		\left(
		\frac{\pi}
		{\sqrt6(s-\tfrac1{24})^{3/2}}
		-
		\frac3{2(s-\tfrac1{24})^2}
		+
		\frac{3\sqrt6}
		{4\pi(s-\tfrac1{24})^{5/2}}
		\right).
	\end{align}
	Extracting the common positive factor from
	\eqref{eq:Pprime-expanded} gives
	\begin{align*}
		P'(s)
		&=
		\frac{\pi}{12\sqrt2(s-\tfrac1{24})^{3/2}}
		\exp\!\left(
		\pi\sqrt{\tfrac23(s-\tfrac1{24})}
		\right)
		\left(
		1-
		\frac{3\sqrt6}{2\pi\sqrt{s-\tfrac1{24}}}
		+
		\frac9{2\pi^2(s-\tfrac1{24})}
		\right)\\
		&=\frac{\pi}{12\sqrt2(s-\tfrac1{24})^{3/2}}
		\exp\!\left(
		\pi\sqrt{\tfrac23(s-\tfrac1{24})}
		\right)
		H_0\!\left(
		\frac{\sqrt6}{2\pi\sqrt{s-\tfrac1{24}}}
		\right).
\end{align*}
	This proves the identity in \eqref{eq:Pprime-form}.
	
	We next prove $P'(s) >0$.
	For every real $v$, we have 
	\begin{align}\label{eq:H0-positive}
		H_0(v)
		=
		3\left(v-\frac12\right)^2+\frac14
		\ge\frac14>0.
	\end{align}
	All other factors in \eqref{eq:Pprime-form} are positive
	for real $s>\tfrac1{24}$.
	Hence $P'(s)>0$, completing the proof of
	\eqref{eq:Pprime-form}.
	
	We now prove \eqref{eq:Pprime-lower}.
	Suppose that $s-\tfrac1{24}\ge120$.
	Then
	\begin{align}\label{eq:H0-input-bound}
		0<
		\frac{\sqrt6}{2\pi\sqrt{s-\tfrac1{24}}}
		\le
		\frac1{\pi\sqrt{80}}
		<
		\frac1{25}.
	\end{align}
	For $0<v<1/25$, we have  
	$
		H_0(v)=1-3v+3v^2>1-3v>\frac{22}{25}.
	$
Taking $v=\frac{\sqrt6}{2\pi\sqrt{s-\tfrac1{24}}}$,
we obtain
\begin{align}\label{eq:H0-lower}
	H_0\!\left(
	\frac{\sqrt6}{2\pi\sqrt{s-\tfrac1{24}}}
	\right)
	>
	\frac{22}{25}.
\end{align}
By~\eqref{eq:H0-lower} and	$\frac{\pi}{12\sqrt2}>\frac16$, 
	we get
	\begin{align*}
		P'(s)
		>
		\frac{1}{6}\frac{22}{25}
		\frac1{(s-\tfrac1{24})^{3/2}}
		\exp\!\left(
		\pi\sqrt{\tfrac23(s-\tfrac1{24})}
		\right)
		>
		\frac1{20s^{3/2}}
		\exp\!\left(
		\pi\sqrt{\tfrac23(s-\tfrac1{24})}
		\right).
	\end{align*}
	This proves \eqref{eq:Pprime-lower}.
\end{proof}

\subsection{Uniform derivative bounds}

We first write $P',P'',P'''$ in a common form.
Recall that $H_0(v)=1-3v+3v^2.$
Define
\begin{align*}
	H_1(v)&=1-6v+15v^2-15v^3,\\
	H_2(v)&=1-10v+45v^2-105v^3+105v^4.
\end{align*}
\begin{lemma}\label{lem:P-higher-derivatives}
	Let $P$ be defined by \eqref{eq:PQ}.
	For every complex $s$ with $\Re s>\tfrac1{24}$
	and every $r\in\{0,1,2\}$, we have
	\begin{align}\label{eq:Pderiv}
		P^{(r+1)}(s)
		={}&
		\frac1{4\sqrt3}
		\left(\frac{\pi}{\sqrt6}\right)^{r+1}
		\frac{
			\exp\!\left(\pi\sqrt{\tfrac23(s-\tfrac1{24})}\right)
		}{(s-\tfrac1{24})^{(r+3)/2}}
		H_r\!\left(
		\frac{\sqrt6}{2\pi\sqrt{s-\tfrac1{24}}}
		\right).
	\end{align}
	All square roots are chosen to have positive real part.
\end{lemma}
\begin{proof}
The case $r=0$ is \eqref{eq:Pprime-form}.
Put
$
	v=\frac{\sqrt6}{2\pi\sqrt{s-\tfrac1{24}}}.
$
Then
\begin{align}\label{eq:v-derivative}
	\frac{\mathrm dv}{\mathrm ds}
	=
	-\frac{v}{2(s-\tfrac1{24})},
	\qquad
	\frac1{2(s-\tfrac1{24})}
	=
	\frac{\pi v}{\sqrt{6(s-\tfrac1{24})}}.
\end{align}
We first prove \eqref{eq:Pderiv} for $r=1$.
Differentiating \eqref{eq:Pprime-form} with respect
to $s$ and using \eqref{eq:v-derivative}, the product
and chain rules give
\begin{align}\label{eq:Psecond-calculation}
	P''(s)
	={}&
	\frac1{4\sqrt3}
	\left(\frac{\pi}{\sqrt6}\right)^2
	\frac{
		\exp\!\left(\pi\sqrt{\tfrac23(s-\tfrac1{24})}\right)
	}{(s-\tfrac1{24})^2}
	\Bigl(
	(1-3v)H_0(v)-v^2H_0'(v)
	\Bigr).
\end{align}
Since $H_0(v)=1-3v+3v^2$, we get 
\begin{align}
	(1-3v)H_0(v)-v^2H_0'(v)
	&=(1-3v)(1-3v+3v^2)-v^2(-3+6v)\notag\\
	&=1-6v+15v^2-15v^3\notag\\
	&=H_1(v).\label{eq:H1v}
\end{align}
Substitution of~\eqref{eq:H1v} into \eqref{eq:Psecond-calculation}
proves \eqref{eq:Pderiv} for $r=1$.

We next prove \eqref{eq:Pderiv} for $r=2$.
By~\eqref{eq:H1v} and~\eqref{eq:Psecond-calculation}, we obtain
\begin{align}\label{eq:Pthird-calculation}
	P'''(s)
	={}&
	\frac1{4\sqrt3}
	\left(\frac{\pi}{\sqrt6}\right)^3
	\frac{
		\exp\!\left(\pi\sqrt{\tfrac23(s-\tfrac1{24})}\right)
	}{(s-\tfrac1{24})^{5/2}}
	\Bigl(
	(1-4v)H_1(v)-v^2H_1'(v)
	\Bigr).
\end{align}
Here we again used the product and chain rules,
together with \eqref{eq:v-derivative}.
By the definition of $H_1$, we have
\begin{align*}
	(1-4v)H_1(v)-v^2H_1'(v)
	&=(1-4v)(1-6v+15v^2-15v^3)
	-v^2(-6+30v-45v^2)\\
	&=1-10v+45v^2-105v^3+105v^4\\
	&=H_2(v).
\end{align*}
Substitution into \eqref{eq:Pthird-calculation}
proves \eqref{eq:Pderiv} for $r=2$.
\end{proof}
Fix $k\in\{1,2\}$.
For real $n>\tfrac1{24}$ and real $m$, write
\begin{align}\label{eq:central-parameters}
	\varepsilon=\frac1{\sqrt{n-\tfrac1{24}}}.
\end{align}
Since $P'(n)>0$ by \eqref{eq:Pprime-form}, define,
for $r=0,1,2$,
\begin{align}\label{eq:I-normalized}
	I_r(m,n)=
	\frac{\partial_m^r\A_k(m,n)}
	{(\pi/\sqrt{6(n-\tfrac1{24})})^rP'(n)}.
\end{align}

The integral formulas \eqref{eq:Aint},
\eqref{eq:Aint-first-derivative}, and
\eqref{eq:Aint-second-derivative} express these
functions using $B',B'',B'''$.
We first study the corresponding expressions
with $P$ in place of $B$.
For $t\ge0$ and $-\tfrac12\le u\le\tfrac12$, define
\begin{align}\label{eq:Fr-ratio}
	F_r(m,n;t,u)=
	\frac{
		P^{(r+1)}\!\left(
		n+\left(k-\tfrac12\right)t^2-i(m+u)t
		\right)
	}{
		(\pi/\sqrt{6(n-\tfrac1{24})})^rP'(n)
	}.
\end{align}

To write \eqref{eq:Fr-ratio} explicitly, define
\begin{align}\label{eq:R-definition}
	R=R(m,n;t,u)
	=
	\sqrt{
		1+
		\frac{
			\left(k-\tfrac12\right)t^2-i(m+u)t
		}{n-\tfrac1{24}}
	}.
\end{align}
We choose the square root with positive real part,
since
\begin{align*}
	\Re(R^2)
	=
	1+\frac{(k-\tfrac12)t^2}{n-\tfrac1{24}}
	\ge1.
\end{align*}
Thus
\begin{align}\label{eq:R-root-relation}
	\sqrt{
		n+\left(k-\tfrac12\right)t^2-i(m+u)t-\tfrac1{24}
	}
	=
	R\sqrt{n-\tfrac1{24}}.
\end{align}
Using \eqref{eq:Pderiv} in \eqref{eq:Fr-ratio}
and simplifying, we obtain
\begin{align}\label{eq:Fr}
	F_r(m,n;t,u)
	={}&
	\frac{
		H_r\!\left(\frac{\sqrt6\varepsilon}{2\pi R}\right)
	}{
		R^{r+3}H_0\!\left(\frac{\sqrt6\varepsilon}{2\pi}\right)
	}
	\exp\!\left(
	\frac{\pi\sqrt{\tfrac23}\,\varepsilon t}{R+1}
	\Bigl(\left(k-\tfrac12\right)t-i(m+u)\Bigr)
	\right).
\end{align}
By \eqref{eq:R-definition} and \eqref{eq:central-parameters},
we have
\begin{align*}
	\sqrt{n-\tfrac1{24}}(R-1)
	&=
	\frac{\sqrt{n-\tfrac1{24}}(R^2-1)}{R+1}=
	\frac{\varepsilon t}{R+1}
	\Bigl(\left(k-\tfrac12\right)t-i(m+u)\Bigr).
\end{align*}

We expand \eqref{eq:Fr} in powers of $\varepsilon$,
keeping $k,t,u$ and $m/\sqrt{n-\tfrac1{24}}$ fixed.
Under this convention,
\begin{align*}
	R^2
	=
	1-\frac{imt}{\sqrt{n-\tfrac1{24}}}\varepsilon
	+\left(\left(k-\tfrac12\right)t^2-iut\right)
	\varepsilon^2.
\end{align*}
Hence $R=1$ at $\varepsilon=0$.
Since $H_r(0)=H_0(0)=1$, \eqref{eq:Fr} extends
smoothly to $\varepsilon=0$, with
\begin{align}\label{eq:Fr-at-zero}
	\left.F_r(m,n;t,u)\right|_{\varepsilon=0}
	=
	\exp\!\left(
	-\frac{i\pi mt}{\sqrt{6(n-\tfrac1{24})}}
	\right).
\end{align}

By \eqref{eq:Aint} and the chain rule, each differentiation
with respect to $m$ contributes a factor $-it$.
Hence, for $r=0,1,2$, \eqref{eq:I-normalized} and
\eqref{eq:Fr-ratio} show that replacing $B$ by $P$
gives the normalized integrand
\begin{align}\label{eq:P-normalized-integrand}
	\frac{
		\Re\!\left((-i)^rt^{r+1}F_r(m,n;t,u)\right)
	}{\sinh(\pi t)}.
\end{align}
By \eqref{eq:W}, for $\Re s>\tfrac1{24}$, we have
\begin{align*}
	B^{(r+1)}(s)-P^{(r+1)}(s)
	=
	W^{(r+1)}(s),
	\qquad r=0,1,2,
\end{align*}
where the derivatives are taken with respect to $s$.
In \eqref{eq:Aint} and its derivatives with respect
to $m$, these functions are evaluated at
\begin{align*}
	s=n+\left(k-\tfrac12\right)t^2-i(m+u)t.
\end{align*}
Lemma~\ref{lem:tails} bounds the contribution of
$W^{(r+1)}(s)$ for $0\le t\le24$ and
$-\tfrac12\le u\le\tfrac12$.

We next expand $F_r(m,n;t,u)$ in powers of
$\varepsilon$ at zero.
During this expansion, $k,t,u$ and
$m/\sqrt{n-\tfrac1{24}}$ are fixed.
Thus $n=\tfrac1{24}+\varepsilon^{-2}$ and $m$
vary with $\varepsilon$, while their ratio
$m/\sqrt{n-\tfrac1{24}}$ remains constant.
At $\varepsilon=0$, we use the continuous extension
of $F_r$.

We establish the estimates for
\begin{align}\label{eq:domain}
	0\le\varepsilon\le\frac1{66},
	\qquad
	0\le\frac{m}{\sqrt{n-\tfrac1{24}}}\le\frac{101}{50},
	\qquad
	k\in\{1,2\}.
\end{align}
For $\varepsilon>0$, this domain includes all real
$n\ge5000$ and $0\le m\le2\sqrt n+1$.

The next lemma bounds the derivatives with respect to
$\varepsilon$ that occur in the quadratic Taylor expansion
of $F_r$.

Define
\begin{align}\label{eq:majorant-U}
	U(t)=30+100t+100t^2+50t^3+13t^4,
\end{align}
and, for $r=0,1,2$, define
$P_r(t)=\sum_{j=0}^{6}p_{r,j}t^j$ by
\begin{align}
	P_0(t)
	&=35+128t+381t^2+553t^3
	+462t^4+199t^5+40t^6,
	\notag\\
	P_1(t)
	&=69+328t+824t^2+987t^3
	+670t^4+243t^5+41t^6,
	\notag\\
	P_2(t)
	&=158+784t+1631t^2+1636t^3
	+932t^4+290t^5+43t^6.
	\label{eq:majorant-polynomials}
\end{align}

\begin{lemma}\label{lem:majorants}
	Let $F_r(m,n;t,u)$ be defined by \eqref{eq:Fr},
	with parameters satisfying \eqref{eq:domain}.
	For $r=0,1,2$, $t\ge0$, and
	$-\tfrac12\le u\le\tfrac12$, we have
	\begin{align}
		\left|
		\left.
		\partial_\varepsilon^2F_r(m,n;t,u)
		\right|_{\varepsilon=0}
		\right|
		&\le U(t),
		\label{eq:majorants2}\\
		\intertext{and, for $0\le t\le8$,}
		\left|
		\partial_\varepsilon^3F_r(m,n;t,u)
		\right|
		&\le
		\exp\!\left(\frac{t^2}{20}\right)P_r(t).
		\label{eq:majorants}
	\end{align}
\end{lemma}

\begin{proof}
All derivatives with respect to $\varepsilon$ are taken with
$k,t,u$ and $m/\sqrt{n-\tfrac1{24}}$ fixed.
We write $R'=\mathrm dR/\mathrm d\varepsilon$. We first bound the factors in $F_r(m,n;t,u)$,
then their derivatives, and finally combine the bounds
by comparing Taylor coefficients.
	
	\medskip\noindent
	\textbf{Step 1: the size of the factors.}
	
	By \eqref{eq:R-definition}, we have
	\begin{align}\label{eq:majorant-R-square}
		R^2
		=1-\frac{imt}{\sqrt{n-\tfrac1{24}}}\varepsilon
		+\left(\left(k-\tfrac12\right)t^2-iut\right)
		\varepsilon^2.
	\end{align}
	Since $\Re(R^2)\ge1$ and the square root has positive real part,
	it follows that $|R|\ge1$ and $\Re R\ge1$.
	Also, by \eqref{eq:domain}, we get
	\begin{align}\label{eq:majorant-scaled-m}
		\left|
		\frac{m}{\sqrt{n-\tfrac1{24}}}+u\varepsilon
		\right|
		\le\frac{101}{50}+\frac1{132}
		<\frac{203}{100}.
	\end{align}
	By $\Re\sqrt{A+iB}\le\sqrt{A+B^2/(4A)}$
	for real $A>0$ and $B$, and
	$\sqrt{1+x}\le1+x/2$ for $x\ge0$, we obtain
	\begin{align*}
		\Re R-1
		\le
		\varepsilon^2t^2
		\left(
		\frac{2k-1}{4}
		+\frac18
		\left(
		\frac{m}{\sqrt{n-\tfrac1{24}}}+u\varepsilon
		\right)^2
		\right).
	\end{align*}
Using \eqref{eq:majorant-scaled-m}, together with
$k\le2$, $0<\varepsilon\le1/66$, and
$\pi\sqrt{2/3}<257/100$, we obtain
	\begin{align}\label{eq:majorant-basic}
		\Re\!\left(
		\pi\sqrt{\frac23}\frac{R-1}{\varepsilon}
		\right)
		\le
		\frac{257}{6600}
		\left(
		\frac34+\frac18\left(\frac{203}{100}\right)^2
		\right)t^2
		\le\frac{t^2}{20}.
	\end{align}
	At $\varepsilon=0$, the smooth extension of
	$\pi\sqrt{2/3}(R-1)/\varepsilon$ is purely imaginary,
	so its real part is zero.
	
	We also need to bound the expression evaluated in $H_r$.
	By $|R|\ge1$ and $\sqrt6/(2\pi)<25/64$, we have
	\begin{align}\label{eq:majorant-H-variable}
		\left|\frac{\sqrt6\varepsilon}{2\pi R}\right|
		\le\frac{25}{64\cdot66}
		<\frac{37}{6250}.
	\end{align}
	
	\medskip\noindent
	\textbf{Step 2: derivatives of the factors.}
	
	To differentiate $R$, set
	\begin{align}\label{eq:majorant-Xi}
		\xi
		=\left(k-\tfrac12\right)t^2-iut
		+\frac{m^2t^2}{4(n-\tfrac1{24})}.
	\end{align}
	Since $m/\sqrt{n-\tfrac1{24}}$ is fixed,
	$\xi$ is independent of $\varepsilon$.
	By \eqref{eq:domain}, we have
	\begin{align}\label{eq:majorant-Xi-bound}
		|\xi|
		\le
		\frac t2+
		\left(\frac32+\frac{10201}{10000}\right)t^2
		\le
		\frac t2+\frac{2521}{1000}t^2.
	\end{align}
	Differentiating \eqref{eq:majorant-R-square} once gives
	\begin{align*}
		R'
		=
		\frac{
			-imt/\sqrt{n-\tfrac1{24}}
			+2\left((k-\tfrac12)t^2-iut\right)\varepsilon
		}{2R}.
	\end{align*}
By the formula for $R'$ and the triangle inequality, we obtain
\begin{align}\label{eq:majorant-R-first-bound}
	|R'|
	&\le
	\frac1{2|R|}
	\left(
	\frac{mt}{\sqrt{n-\tfrac1{24}}}
	+2\varepsilon
	\left(
	\left(k-\tfrac12\right)t^2+|u|t
	\right)
	\right)
	\notag\\
	&\le
	\frac{101}{100}t
	+\frac1{66}\left(\frac32t^2+\frac t2\right)
	\notag\\
	&\le
	\frac{509}{500}t+\frac{2273}{100000}t^2.
\end{align}
Here the second inequality uses $|R|\ge1$,
 $k\le2$,  $|u|\le1/2$ and \eqref{eq:domain}.
	Differentiating $RR'$ and using \eqref{eq:majorant-Xi},
	we get $R''=\xi/R^3$.
	Since $\xi$ is constant with respect to $\varepsilon$,
	two further differentiations give
	\begin{align}\label{eq:majorant-root}
		R'''=-\frac{3\xi R'}{R^4},
		\qquad
		R^{(4)}
		=-\frac{3\xi^2}{R^7}
		+\frac{12\xi(R')^2}{R^5}.
	\end{align}
	Thus $|R''|\le|\xi|$, $|R'''|\le3|\xi||R'|$, and
	$|R^{(4)}|\le3|\xi|^2+12|\xi||R'|^2$.
	
	We next estimate the derivatives of the exponent in
	\eqref{eq:Fr}.
	Since $R(0)=1$, the fundamental theorem of calculus gives
	\begin{align}\label{eq:majorant-exponent-integral}
		\pi\sqrt{\frac23}
		\frac{R(\varepsilon)-1}{\varepsilon}
		=
		\pi\sqrt{\frac23}
		\int_0^1R'(v\varepsilon)\dd v.
	\end{align}
For $j=1,2,3$, differentiating
\eqref{eq:majorant-exponent-integral}
and applying the triangle inequality, we obtain
\begin{align}\label{eq:majorant-exponent-derivatives}
	\left|
	\partial_\varepsilon^j
	\left(
	\pi\sqrt{\frac23}\frac{R-1}{\varepsilon}
	\right)
	\right|
	&=
	\pi\sqrt{\frac23}
	\left|
	\int_0^1v^jR^{(j+1)}(v\varepsilon)\dd v
	\right|
	\notag\\
	&\le
	\pi\sqrt{\frac23}
	\int_0^1v^j
	\left|R^{(j+1)}(v\varepsilon)\right|\dd v
	\le E_j(t).
\end{align}
For the last inequality, we use the derivative formulas in
\eqref{eq:majorant-root}, together with
\eqref{eq:majorant-Xi-bound} and
\eqref{eq:majorant-R-first-bound}.
	For the first derivative, one may take
	$E_1(t)=643t/1000+81t^2/25$.
	For the second derivative, integrating $v^2$ in
	\eqref{eq:majorant-exponent-integral} gives
	\begin{align}\label{eq:majorant-E2}
		E_2(t)
		=\frac{257t^2}{100}
		\left(\frac12+\frac{2521}{1000}t\right)
		\left(\frac{509}{500}+\frac{2273}{100000}t\right).
	\end{align}
	Similarly, the bound for $R^{(4)}$ and
	$\int_0^1v^3\dd v=1/4$ give
	\begin{align}\label{eq:majorant-E3}
		E_3(t)
		=
		\frac{257t^2}{100}
		\Biggl(
		\frac34\left(\frac12+\frac{2521}{1000}t\right)^2+
		3t\left(\frac12+\frac{2521}{1000}t\right)
		\left(\frac{509}{500}+\frac{2273}{100000}t\right)^2
		\Biggr).
	\end{align}
	
	We now estimate the derivatives of
	$\sqrt6\varepsilon/(2\pi R)$ directly.
	By \eqref{eq:majorant-R-square}, we have
	\begin{align}\label{eq:majorant-H-variable-first}
		\partial_\varepsilon
		\left(\frac{\sqrt6\varepsilon}{2\pi R}\right)
		=
		\frac{\sqrt6}{2\pi R^3}
		\left(
		1-\frac{imt\varepsilon}
		{2\sqrt{n-\tfrac1{24}}}
		\right).
	\end{align}

For $t\ge0$, define
\begin{align}\label{eq:majorant-V1}
	V_1(t)=\frac{25}{64}
	\left(1+\frac{1531}{100000}t\right).
\end{align}
For the second derivative, define
\begin{align}\label{eq:majorant-V2}
	V_2(t)
	=\frac{25t}{64}
	\left(
	\frac{101}{100}
	+3\left(1+\frac{1531}{100000}t\right)
	\left(\frac{509}{500}+\frac{2273}{100000}t\right)
	\right).
\end{align}
For the third derivative, define
\begin{align}\label{eq:majorant-V3}
	V_3(t)
	={}&
	\frac{25t}{64}
	\Biggl(
	\frac{303t}{50}
	\left(\frac{509}{500}+\frac{2273}{100000}t\right)+
	3\left(1+\frac{1531}{100000}t\right)
	\left(\frac12+\frac{2521}{1000}t\right)
	\notag\\
	&\quad+
	12t\left(1+\frac{1531}{100000}t\right)
	\left(\frac{509}{500}+\frac{2273}{100000}t\right)^2
	\Biggr).
\end{align}
For $j=1,2,3$, we have
\begin{align}\label{eq:majorant-H-variable-derivatives}
	\left|
	\partial_\varepsilon^j
	\left(\frac{\sqrt6\varepsilon}{2\pi R}\right)
	\right|
	\le V_j(t).
\end{align}

Indeed, by \eqref{eq:domain}, we have
\begin{align*}
	\left|
	1-\frac{imt\varepsilon}{2\sqrt{n-\tfrac1{24}}}
	\right|
	\le1+\frac{1531}{100000}t,
	{\text{ and }}
	\left|
	\partial_\varepsilon
	\left(
	1-\frac{imt\varepsilon}{2\sqrt{n-\tfrac1{24}}}
	\right)
	\right|
	\le\frac{101}{100}t.
\end{align*}
Thus \eqref{eq:majorant-H-variable-first},
$|R|\ge1$, and $\sqrt6/(2\pi)<25/64$
give the case $j=1$ by \eqref{eq:majorant-V1}.

Differentiating \eqref{eq:majorant-H-variable-first}
once more and using \eqref{eq:majorant-R-first-bound},
we obtain the case $j=2$ with
$V_2(t)$ defined in \eqref{eq:majorant-V2}.
Finally, differentiating again and using
$|R''|\le|\xi|$, we apply
\eqref{eq:majorant-Xi-bound} and
\eqref{eq:majorant-R-first-bound}
to obtain the case $j=3$ with
$V_3(t)$ defined in \eqref{eq:majorant-V3}.

Since $H_0(v)=1-3v+3v^2$, we have
	\begin{align}\label{eq:majorant-denominator}
		H_0\!\left(\frac{\sqrt6\varepsilon}{2\pi}\right)
		>\frac{491}{500}.
	\end{align}
	Direct differentiation of this quadratic also gives
	\begin{align}\label{eq:majorant-denominator-derivatives}
		\left|
		\partial_\varepsilon
		H_0\!\left(\frac{\sqrt6\varepsilon}{2\pi}\right)
		\right|
		<\frac{293}{250},
		\qquad
		\left|
		\partial_\varepsilon^2
		H_0\!\left(\frac{\sqrt6\varepsilon}{2\pi}\right)
		\right|
		<\frac{229}{250},\qquad
		\partial_\varepsilon^3
		H_0\!\left(\frac{\sqrt6\varepsilon}{2\pi}\right)
		=0.
	\end{align}
	
	\medskip\noindent
	\textbf{Step 3: comparison of Taylor coefficients.}
	
	By \eqref{eq:Fr}, the factors to be combined are
	\begin{align}\label{eq:majorant-product}
		F_r(m,n;t,u)
		=
		\frac{
			\exp\!\left(
			\pi\sqrt{\frac23}\frac{R-1}{\varepsilon}
			\right)
			H_r\!\left(\frac{\sqrt6\varepsilon}{2\pi R}\right)
		}{
			R^{r+3}H_0\!\left(\frac{\sqrt6\varepsilon}{2\pi}\right)
		}.
	\end{align}
	For a real increment $x$ near zero, define
	\begin{align}\label{eq:majorant-generating-function}
		\mathcal G_r(x,t)
		=&
		\frac{
			\exp\!\left(
			E_1(t)x+\frac{E_2(t)}2x^2+\frac{E_3(t)}6x^3
			\right)
		}{
			\frac{491}{500}-\frac{293}{250}x-\frac{229}{500}x^2
		}
		H_r\!\left(
		-\frac{37}{6250}
		-V_1(t)x-\frac{V_2(t)}2x^2-\frac{V_3(t)}6x^3
		\right)
		\notag\\
		&\quad\times
		\Biggl(
		1-t\left(\frac{509}{500}+\frac{2273}{100000}t\right)x
		-\frac t2\left(\frac12+\frac{2521}{1000}t\right)x^2
		\notag\\
		&\hspace{19mm}
		-\frac{t^2}{2}
		\left(\frac12+\frac{2521}{1000}t\right)
		\left(\frac{509}{500}+\frac{2273}{100000}t\right)x^3
		\Biggr)^{-r-3}.
	\end{align}
	
	We compare only the coefficients of $1,x,x^2,x^3$,
	using the triangle inequality for products.
	
	For the exponential factor in \eqref{eq:majorant-product},
	divide its value at $\varepsilon+x$ by its value
	at $\varepsilon$.
	By \eqref{eq:majorant-exponent-derivatives}
	and the exponential series, the absolute values
	of the first four coefficients are bounded by
	the corresponding coefficients of
	\begin{align*}
		\exp\!\left(
		E_1(t)x+\frac{E_2(t)}2x^2+\frac{E_3(t)}6x^3
		\right).
	\end{align*}
	This gives the exponential factor in
	\eqref{eq:majorant-generating-function}.
	
	For $R^{-r-3}$, factor out $R(\varepsilon)^{-r-3}$
	and apply the binomial series.
	The bounds for $R',R'',R'''$ and $|R|\ge1$
	yield the factor with power $-r-3$ in
	\eqref{eq:majorant-generating-function}.
	
	For $H_r(\sqrt6\varepsilon/(2\pi R))$,
	use \eqref{eq:majorant-H-variable} and
	\eqref{eq:majorant-H-variable-derivatives}.
	Since $H_r(-v)$ has nonnegative coefficients,
	these estimates yield the $H_r$ factor in
	\eqref{eq:majorant-generating-function}.
	
	By \eqref{eq:majorant-denominator} and
	\eqref{eq:majorant-denominator-derivatives},
	the geometric series bounds the Taylor coefficients
	of $H_0(\sqrt6\varepsilon/(2\pi))^{-1}$
	by those of
	$(491/500-293x/250-229x^2/500)^{-1}$.
	
For $j=2,3$, the coefficient of $x^j$ equals
the corresponding $j$th derivative divided by $j!$.
Thus the four coefficient comparisons yield
\begin{align}\label{eq:majorant-coefficient-transfer}
	|\partial_\varepsilon^jF_r(m,n;t,u)|
	\le
	\exp\!\left(
	\Re\!\left(
	\pi\sqrt{\frac23}\frac{R-1}{\varepsilon}
	\right)
	\right)
	\partial_x^j\mathcal G_r(0,t).
\end{align}
	
	We now check the resulting polynomials in $t$.
	To compute the coefficients of $x^2,x^3$, use
	\begin{align*}
		\exp\!\left(Ax+\frac B2x^2+\frac C6x^3\right)
		=
		1+Ax+\frac{A^2+B}{2}x^2
		+\frac{A^3+3AB+C}{6}x^3+O(x^4),
	\end{align*}
	together with the geometric and binomial series.
	We differentiate $\mathcal G_2(x,t)$ in
	\eqref{eq:majorant-generating-function} twice
	with respect to $x$ and set $x=0$.
	Expanding the result in powers of $t$ and comparing
	its coefficients with those of $U(t)$, we obtain
	\begin{align}\label{eq:majorant-final-comparisons}
		1000\left(U(t)-\partial_x^2\mathcal G_2(0,t)\right)
		\ge
		1331+16777t+5854t^2+610t^3+686t^4.
	\end{align}
	
	For $r=0$, we differentiate $\mathcal G_0(x,t)$
	three times with respect to $x$ and set $x=0$.
	Comparing the resulting coefficients with those
	of $P_0(t)$, we obtain
	\begin{align}\label{eq:majorant-P0-comparison}
		1000\left(P_0(t)-\partial_x^3\mathcal G_0(0,t)\right)
		\ge
		462+143t+919t^2+964t^3+576t^4+298t^5+925t^6.
	\end{align}
	For $r=1$, we obtain
	\begin{align}\label{eq:majorant-P1-comparison}
		1000\left(P_1(t)-\partial_x^3\mathcal G_1(0,t)\right)
		\ge
		836+591t+626t^2+174t^3+618t^4+918t^5+386t^6.
	\end{align}
	For $r=2$, we have
	\begin{align}\label{eq:majorant-P2-comparison}
		1000\left(P_2(t)-\partial_x^3\mathcal G_2(0,t)\right)
		\ge
		863+708t+936t^2+282t^3+693t^4+457t^5+607t^6.
	\end{align}
	
	All coefficients on the right are positive.
	Moreover, the coefficients of
	\eqref{eq:majorant-generating-function}
	are nonnegative and increase with $r=0,1,2$.
	It follows that
	\begin{align*}
		\partial_x^2\mathcal G_r(0,t)
		\le\partial_x^2\mathcal G_2(0,t)\le U(t),
		\qquad
		\partial_x^3\mathcal G_r(0,t)\le P_r(t).
	\end{align*}
	At $\varepsilon=0$, the exponential factor in
	\eqref{eq:majorant-coefficient-transfer} equals one.
	For $0<\varepsilon\le1/66$, it is bounded by
	$\exp(t^2/20)$ by \eqref{eq:majorant-basic}.
	Substitution into \eqref{eq:majorant-coefficient-transfer}
	proves \eqref{eq:majorants2} and \eqref{eq:majorants}.
\end{proof}

\subsection{Integral tails and Taylor remainder}

For every nonnegative integer $j$, real $\alpha>0$,
and real $b\ge0$, repeated integration by parts gives
\begin{align}\label{eq:tailint}
	\int_b^\infty t^j e^{-\alpha t}\dd t
	=
	e^{-\alpha b}
	\sum_{\ell=0}^j
	\frac{j!\,b^\ell}{\ell!\,\alpha^{j-\ell+1}}.
\end{align}
Indeed, for $j\ge1$,
\begin{align*}
	\int_b^\infty t^j e^{-\alpha t}\dd t
	=
	\frac{b^j e^{-\alpha b}}{\alpha}
	+\frac j\alpha
	\int_b^\infty t^{j-1}e^{-\alpha t}\dd t,
\end{align*}
and the case $j=0$ is immediate.
We also use
\begin{align*}
	e^3>\sum_{j=0}^8\frac{3^j}{j!}>20,
	\qquad
	e^{3/5}>
	1+\frac35+\frac{(3/5)^2}{2}+\frac{(3/5)^3}{6}
	>\frac95.
\end{align*}

\begin{lemma}\label{lem:tails}
	Let $k\in\{1,2\}$, and let $n>\tfrac1{24}$ and
	$m\ge0$ be real.
	Assume \eqref{eq:domain}, with $\varepsilon$ defined
	by \eqref{eq:central-parameters}.
	For $r=0,1,2$, $t\ge0$, and
	$-\tfrac12\le u\le\tfrac12$, define
	\begin{align}\label{eq:tail-Taylor-polynomial}
		T_r(m,n;t,u)
		={}&
		\left.F_r(m,n;t,u)\right|_{\varepsilon=0}
		+\varepsilon
		\left.\partial_\varepsilon F_r(m,n;t,u)
		\right|_{\varepsilon=0}+
		\frac{\varepsilon^2}{2}
		\left.\partial_\varepsilon^2F_r(m,n;t,u)
		\right|_{\varepsilon=0}.
	\end{align}
	The derivatives are taken with $k,t,u$ and
	$m/\sqrt{n-\tfrac1{24}}$ fixed, using the smooth
	extension at $\varepsilon=0$.
	Then
	\begin{align}\label{eq:tails-comparison}
		&\Biggl|
		I_r(m,n)
		-
		\int_0^\infty
		\frac{t^{r+1}}{\sinh(\pi t)}
		\int_{-1/2}^{1/2}
		\Re\!\left((-i)^rT_r(m,n;t,u)\right)
		\dd u\dd t
		\Biggr|
		\notag\\
		&\qquad<
		\frac{\varepsilon^3}{6}
		\int_0^8
		\frac{
			t^{r+1}\exp(t^2/20)P_r(t)
		}{
			\sinh(\pi t)
		}
		\dd t
		+\frac1{10^5},
	\end{align}
	where $P_r(t)$ is the polynomial in
	\eqref{eq:majorant-polynomials}.
\end{lemma}

\begin{proof}
	
	We first estimate the Taylor remainder on $0\le t\le8$.
	By \eqref{eq:tail-Taylor-polynomial} and Taylor's theorem
	with integral remainder, we have
	\begin{align*}
		F_r(m,n;t,u)-T_r(m,n;t,u)
		=
		\frac12\int_0^\varepsilon
		(\varepsilon-\theta)^2
		\left.
		\partial_\varepsilon^3F_r(m,n;t,u)
		\right|_{\varepsilon=\theta}
		\dd\theta.
	\end{align*}
	The differentiation convention is the same as in
	the definition of $T_r(m,n;t,u)$.
	Lemma~\ref{lem:majorants} applies for
	$0\le\theta\le\varepsilon$ and gives
	\begin{align*}
		|F_r(m,n;t,u)-T_r(m,n;t,u)|
		\le
		\frac12\exp\!\left(\frac{t^2}{20}\right)P_r(t)
		\int_0^\varepsilon(\varepsilon-\theta)^2\dd\theta
		=
		\frac{\varepsilon^3}{6}
		\exp\!\left(\frac{t^2}{20}\right)P_r(t).
	\end{align*}
	Using $|\Re((-i)^rz)|\le|z|$ and
	$\int_{-1/2}^{1/2}\dd u=1$, we therefore obtain
	\begin{align}\label{eq:tail-short-remainder}
		&\left|
		\int_0^8\frac{t^{r+1}}{\sinh(\pi t)}
		\int_{-1/2}^{1/2}
		\Re\!\left(
		(-i)^r
		\bigl(F_r(m,n;t,u)-T_r(m,n;t,u)\bigr)
		\right)
		\dd u\dd t
		\right|
		\notag\\
		\le
		&\frac{\varepsilon^3}{6}
		\int_0^8
		\frac{t^{r+1}\exp(t^2/20)P_r(t)}
		{\sinh(\pi t)}
		\dd t.
	\end{align}
	
By \eqref{eq:Aint}, \eqref{eq:I-normalized},
\eqref{eq:Fr-ratio}, and the identity
\begin{align*}
	B^{(r+1)}(s)
	=
	P^{(r+1)}(s)+W^{(r+1)}(s),
\end{align*}
the triangle inequality gives
\begin{align}\label{eq:tails-decomposition}
	&\Biggl|
	I_r(m,n)
	-
	\int_0^8\frac{t^{r+1}}{\sinh(\pi t)}
	\int_{-1/2}^{1/2}
	\Re\!\left((-i)^rF_r(m,n;t,u)\right)
	\dd u\dd t
	\notag\\
	&\qquad-
	\int_8^\infty\frac{t^{r+1}}{\sinh(\pi t)}
	\int_{-1/2}^{1/2}
	\Re\!\left((-i)^rT_r(m,n;t,u)\right)
	\dd u\dd t
	\Biggr|
	\notag\\
	\le&
	\mathcal E_{r,1}(m,n)
	+\mathcal E_{r,2}(m,n)
	+\mathcal E_{r,3}(m,n)
	+\mathcal E_{r,4}(m,n),
\end{align}
where
\begin{align}\label{eq:tail-W-integral}
	\mathcal E_{r,1}(m,n)
	={}&
	\frac1{(\pi\varepsilon/\sqrt6)^rP'(n)}
	\int_0^{24}\frac{t^{r+1}}{\sinh(\pi t)}
	\notag\\
	&\qquad\times
	\int_{-1/2}^{1/2}
	\left|
	W^{(r+1)}\!\left(
	n+\left(k-\tfrac12\right)t^2-i(m+u)t
	\right)
	\right|
	\dd u\dd t,
\end{align}
\begin{align}\label{eq:tail-B-integral}
	\mathcal E_{r,2}(m,n)
	={}&
	\frac1{(\pi\varepsilon/\sqrt6)^rP'(n)}
	\int_{24}^\infty\frac{t^{r+1}}{\sinh(\pi t)}
	\notag\\
	&\qquad\times
	\int_{-1/2}^{1/2}
	\left|
	B^{(r+1)}\!\left(
	n+\left(k-\tfrac12\right)t^2-i(m+u)t
	\right)
	\right|
	\dd u\dd t,
\end{align}
\begin{align}\label{eq:tail-F-integral}
	\mathcal E_{r,3}(m,n)
	=
	\int_8^{24}\frac{t^{r+1}}{\sinh(\pi t)}
	\int_{-1/2}^{1/2}
	|F_r(m,n;t,u)|
	\dd u\dd t,
\end{align}
\begin{align}\label{eq:tail-T-integral}
	\mathcal E_{r,4}(m,n)
	=
	\int_8^\infty\frac{t^{r+1}}{\sinh(\pi t)}
	\int_{-1/2}^{1/2}
	|T_r(m,n;t,u)|
	\dd u\dd t.
\end{align}

We prove
\begin{align*}
	\mathcal E_{r,j}(m,n)<\frac1{10^6},
	\qquad j=1,2,3,4,
\end{align*}
in Steps 1--4, respectively.
Together with \eqref{eq:tail-short-remainder},
these estimates imply \eqref{eq:tails-comparison},
since $4/10^6<1/10^5$.
	
\noindent
\textbf{Step 1: estimating $\mathcal E_{r,1}(m,n)$.}

	Differentiating $W(s)$ in \eqref{eq:W-definition},
	as in Lemma~\ref{lem:P-higher-derivatives}, gives
	\begin{align}\label{eq:tail-W-derivatives}
		W^{(r+1)}(s)
		={}&
		\frac1{4\sqrt3}
		\left(-\frac{\pi}{\sqrt6}\right)^{r+1}
		\frac{
			\exp\!\left(
			-\pi\sqrt{\tfrac23(s-\tfrac1{24})}
			\right)
		}{
			(s-\tfrac1{24})^{(r+3)/2}
		}
		H_r\!\left(
		-\frac{\sqrt6}{2\pi\sqrt{s-\tfrac1{24}}}
		\right),
		\qquad \Re s>\tfrac1{24}.
	\end{align}
	For every complex $s$ with $\Re s\ge n$, we have 
	\begin{align*}
		|s-\tfrac1{24}|\ge\varepsilon^{-2},
		\qquad
		\Re\sqrt{s-\tfrac1{24}}\ge\varepsilon^{-1}.
	\end{align*}
For $r=0,1,2$ and complex $v$ with
$|v|\le37/6250$, the formulas for $H_0,H_1,H_2$
and the triangle inequality give
\begin{align}\label{eq:tail-H-bound}
	|H_r(v)|
	&\le
	1+10|v|+45|v|^2+105|v|^3+105|v|^4<\frac{1061}{1000}.
\end{align}
Also, since
$0<\sqrt6\varepsilon/(2\pi)<37/6250$, we have
\begin{align}\label{eq:tail-H0-bound}
	H_0\!\left(\frac{\sqrt6\varepsilon}{2\pi}\right)
	&=
	1-\frac{3\sqrt6\varepsilon}{2\pi}
	+3\left(\frac{\sqrt6\varepsilon}{2\pi}\right)^2>
	1-\frac{111}{6250}
	>\frac{491}{500}.
\end{align}

We next apply these bounds to $W^{(r+1)}(s)$.
Since $n-\tfrac1{24}=\varepsilon^{-2}$,
\eqref{eq:Pprime-form} gives
\begin{align}\label{eq:tail-Pprime-normalization}
	P'(n)
	=
	\frac{\pi\varepsilon^3}{12\sqrt2}
	\exp\!\left(\frac{\pi\sqrt{2/3}}{\varepsilon}\right)
	H_0\!\left(\frac{\sqrt6\varepsilon}{2\pi}\right).
\end{align}
For $\Re s\ge n$, the preceding bounds imply
\begin{align*}
	|s-\tfrac1{24}|\ge\varepsilon^{-2},
	\qquad
	\Re\sqrt{s-\tfrac1{24}}\ge\varepsilon^{-1},
\end{align*}
and hence
\begin{align*}
	\left|
	\frac{\sqrt6}{2\pi\sqrt{s-\tfrac1{24}}}
	\right|
	\le\frac{\sqrt6\varepsilon}{2\pi}
	<\frac{37}{6250}.
\end{align*}
Substituting \eqref{eq:tail-W-derivatives} and
\eqref{eq:tail-Pprime-normalization}, and then using
\eqref{eq:tail-H-bound} and \eqref{eq:tail-H0-bound},
we obtain
\begin{align}\label{eq:tail-W-bound}
	\frac{|W^{(r+1)}(s)|}
	{(\pi\varepsilon/\sqrt6)^rP'(n)}
	&<
	\frac{1061}{982}
	\frac{
		\exp\!\left(
		-\pi\sqrt{\frac23}
		\left(
		\Re\sqrt{s-\tfrac1{24}}+\frac1\varepsilon
		\right)
		\right)
	}{
		\varepsilon^{r+3}
		|s-\tfrac1{24}|^{(r+3)/2}
	}
	\notag\\
	&\le
	\frac{1061}{982}
	\exp\!\left(
	-\frac{2\pi\sqrt{2/3}}{\varepsilon}
	\right)
	\notag\\
	&<
	\frac{109}{100}
	\exp\!\left(
	-\frac{2\pi\sqrt{2/3}}{\varepsilon}
	\right).
\end{align}
Here the second inequality uses
$\varepsilon^{r+3}|s-\tfrac1{24}|^{(r+3)/2}\ge1$
and $\Re\sqrt{s-\tfrac1{24}}\ge\varepsilon^{-1}$.

For every $t\ge0$ and $-\tfrac12\le u\le\tfrac12$,
we have
\begin{align*}
	\Re\!\left(
	n+\left(k-\tfrac12\right)t^2-i(m+u)t
	\right)
	=
	n+\left(k-\tfrac12\right)t^2
	\ge n.
\end{align*}
Thus \eqref{eq:tail-W-bound} applies to the function
evaluated in \eqref{eq:tail-W-integral}.

For $t>0$, the inequality $e^{2\pi t}-1\ge2\pi t$
gives
\begin{align}\label{eq:sinh}
	\frac1{\sinh(\pi t)}
	=
	2e^{-\pi t}
	+\frac{2e^{-\pi t}}{e^{2\pi t}-1}
	\le
	e^{-\pi t}\left(2+\frac1{\pi t}\right).
\end{align}
Consequently, by \eqref{eq:tailint}, for $r=0,1,2$, we have 
\begin{align*}
	\int_0^\infty
	\frac{t^{r+1}}{\sinh(\pi t)}\dd t
	\le
	2\int_0^\infty t^{r+1}e^{-\pi t}\dd t
	+\frac1\pi\int_0^\infty t^re^{-\pi t}\dd t
	=
	\frac{2(r+1)!+r!}{\pi^{r+2}}
	<1.
\end{align*}
Since $2\pi\sqrt{2/3}/\varepsilon>18$ and $e^3>20$,
we also have
\begin{align*}
	\exp\!\left(
	-\frac{2\pi\sqrt{2/3}}{\varepsilon}
	\right)
	<
	e^{-18}
	=
	\frac1{(e^3)^6}
	<
	\frac1{20^6}.
\end{align*}
Substituting \eqref{eq:tail-W-bound} into
\eqref{eq:tail-W-integral}, and using the fact that
the interval $-\tfrac12\le u\le\tfrac12$ has length one,
we obtain
\begin{align*}
	\mathcal E_{r,1}(m,n)
	<
	\frac{109}{100}
	\exp\!\left(
	-\frac{2\pi\sqrt{2/3}}{\varepsilon}
	\right)
	\int_0^{24}
	\frac{t^{r+1}}{\sinh(\pi t)}\dd t
	<
	\frac{109}{100\cdot20^6}
	<
	\frac1{10^6}.
\end{align*}
	
\noindent
\textbf{Step 2: estimating $\mathcal E_{r,2}(m,n)$.}
	
	For $\varepsilon>0$ and $t\ge0$, define
	\begin{align}\label{eq:tail-modulus-bound}
		\mathcal R(\varepsilon,t)
		=
		\left(
		\left(1+\frac32\varepsilon^2t^2\right)^2
		+\frac{41209}{10000}\varepsilon^2t^2
		\right)^{1/4}.
	\end{align}
By \eqref{eq:R-definition} and \eqref{eq:domain},
together with $|u|\le\tfrac12$, we have
\begin{align*}
	1\le |R(m,n;t,u)|^4
	&=
	\left(
	1+\left(k-\tfrac12\right)\varepsilon^2t^2
	\right)^2
	+\varepsilon^2t^2
	\left(
	\frac{m}{\sqrt{n-\tfrac1{24}}}+u\varepsilon
	\right)^2\\
	&\le
	\left(1+\frac32\varepsilon^2t^2\right)^2
	+\frac{41209}{10000}\varepsilon^2t^2=\mathcal R(\varepsilon,t)^4,
\end{align*}
where we used
\begin{align*}
	\left|
	\frac{m}{\sqrt{n-\tfrac1{24}}}+u\varepsilon
	\right|
	\le\frac{101}{50}+\frac1{132}
	<\frac{203}{100}.
\end{align*}
Taking nonnegative fourth roots gives
\begin{align*}
	1\le |R(m,n;t,u)|\le\mathcal R(\varepsilon,t).
\end{align*}
	
For real $s\ge n$, we have
\begin{align*}
	0<
	\frac{\sqrt6}{2\pi\sqrt{s-\tfrac1{24}}}
	\le
	\frac{\sqrt6\varepsilon}{2\pi}
	<\frac{37}{6250}<1.
\end{align*}
Since
\begin{align*}
	H_0(v)=1-3v(1-v)\le1
	\qquad 0\le v\le1,
\end{align*}
the polynomial factor in \eqref{eq:Pprime-form}
is at most one.
Moreover, \eqref{eq:tail-W-derivatives} gives
$W'(s)<0$.
Differentiating \eqref{eq:W} and then applying
\eqref{eq:Pprime-form}, we obtain
\begin{align}\label{eq:tail-Bprime-bound}
	B'(s)
	&=P'(s)+W'(s)
	\notag\\
	&<
	\frac{\pi}{12\sqrt2}
	\frac{
		\exp\!\left(
		\pi\sqrt{\tfrac23(s-\tfrac1{24})}
		\right)
	}{
		(s-\tfrac1{24})^{3/2}
	}
	H_0\!\left(
	\frac{\sqrt6}{2\pi\sqrt{s-\tfrac1{24}}}
	\right)
	\notag\\
	&\le
	\frac{\pi}{12\sqrt2}
	\frac{
		\exp\!\left(
		\pi\sqrt{\tfrac23(s-\tfrac1{24})}
		\right)
	}{
		(s-\tfrac1{24})^{3/2}
	}.
\end{align}
	Repeated application of Lemma~\ref{lem:Bbounds}
	then bounds the higher derivatives.
	The function
	\begin{align*}
		\exp\!\left(
		\frac{\pi\sqrt{2/3}}{\varepsilon}(x-1)
		\right)x^{-r-3}
	\end{align*}
	is increasing for $x\ge1$, because its derivative
	has the sign of
	$\pi\sqrt{2/3}/\varepsilon-(r+3)/x>0$.
	We may therefore replace $|R(m,n;t,u)|$ by
	$\mathcal R(\varepsilon,t)$ in the resulting bound.
By \eqref{eq:Pprime-form} and \eqref{eq:tail-H0-bound},
we have
\begin{align*}
	(\pi\varepsilon/\sqrt6)^rP'(n)
	>
	\frac{491}{500}\,
	\frac1{4\sqrt3}
	\left(\frac{\pi}{\sqrt6}\right)^{r+1}
	\varepsilon^{r+3}
	\exp\!\left(
	\frac{\pi\sqrt{2/3}}{\varepsilon}
	\right).
\end{align*}
Combining this lower bound with
Lemma~\ref{lem:Bbounds} and
\eqref{eq:tail-Bprime-bound}, we obtain
\begin{align}\label{eq:tail-B-bound}
	\frac{
		\left|
		B^{(r+1)}\!\left(
		n+\left(k-\tfrac12\right)t^2-i(m+u)t
		\right)
		\right|
	}{
		(\pi\varepsilon/\sqrt6)^rP'(n)
	}
	&\le
	\frac{500}{491}
	\frac{
		\exp\!\left(
		\frac{\pi\sqrt{2/3}}{\varepsilon}
		\bigl(|R(m,n;t,u)|-1\bigr)
		\right)
	}{
		|R(m,n;t,u)|^{r+3}
	}
	\notag\\
	&\le
	\frac{500}{491}
	\frac{
		\exp\!\left(
		\frac{\pi\sqrt{2/3}}{\varepsilon}
		\bigl(\mathcal R(\varepsilon,t)-1\bigr)
		\right)
	}{
		\mathcal R(\varepsilon,t)^{r+3}
	}.
\end{align}
The last inequality follows from
$1\le|R(m,n;t,u)|\le\mathcal R(\varepsilon,t)$
and the monotonicity established above.
	For $t\ge8$,
	\begin{align*}
		\frac1{\sinh(\pi t)}
		<\frac{201}{100}e^{-\pi t}.
	\end{align*}
Substituting \eqref{eq:tail-B-bound} into
\eqref{eq:tail-B-integral}, and using 
$	\frac{500}{491}\frac{201}{100}<3$,
we obtain
\begin{align}\label{eq:large-tail}
	\mathcal E_{r,2}(m,n)
	\le
	3\int_{24}^\infty
	\frac{t^{r+1}}{\mathcal R(\varepsilon,t)^{r+3}}
	\exp\!\left(
	\frac{\pi\sqrt{2/3}}{\varepsilon}
	\bigl(\mathcal R(\varepsilon,t)-1\bigr)-\pi t
	\right)\dd t.
\end{align}

We estimate the integral on the right of
\eqref{eq:large-tail} separately over the regions
$\varepsilon t\le3/10$ and $\varepsilon t\ge3/10$.
To treat the first region, the definition
\eqref{eq:tail-modulus-bound} gives
\begin{align*}
	\mathcal R(\varepsilon,t)^4
	=
	1+\frac{71209}{10000}\varepsilon^2t^2
	+\frac94\varepsilon^4t^4
	\le
	\left(
	1+\frac{71209}{40000}\varepsilon^2t^2
	\right)^4.
\end{align*}
The last inequality follows by comparing coefficients.
Taking fourth roots gives
\begin{align*}
	\mathcal R(\varepsilon,t)
	\le
	1+\frac{71209}{40000}\varepsilon^2t^2.
\end{align*}
For $\varepsilon t\le3/10$, the preceding bound gives
\begin{align*}
	\frac{\pi\sqrt{2/3}}{\varepsilon}
	\bigl(\mathcal R(\varepsilon,t)-1\bigr)-\pi t
	&\le
	\left(
	\frac{257}{100}\frac{71209}{40000}\frac3{10}
	-\frac{157}{50}
	\right)t
	\le-\frac74t.
\end{align*}
Since $\mathcal R(\varepsilon,t)\ge1$ and
$t^{r+1}\le t^3$ for $t\ge24$, the contribution
to \eqref{eq:large-tail} from this range is bounded by
\begin{align*}
	3\int_{24}^\infty t^3e^{-7t/4}\dd t
	<27000e^{-42}
	<\frac{27000}{20^{14}},
\end{align*}
where we used \eqref{eq:tailint} and $e^{42}>20^{14}$.
	
For $\varepsilon t\ge3/10$, the inequality
$\sqrt{3/2}>6/5$ gives
$\sqrt{3/2}\,\varepsilon t\ge9/25$.
Expanding the fourth power and using
\eqref{eq:tail-modulus-bound}, we obtain
\begin{align*}
	\left(\sqrt{\frac32}\varepsilon t+\frac45\right)^4
	-\mathcal R(\varepsilon,t)^4
	&=
	\varepsilon^2t^2
	\left(
	\frac{24}{5}\sqrt{\frac32}\varepsilon t
	-\frac{13609}{10000}
	\right)
	+\frac{256}{125}\sqrt{\frac32}\varepsilon t
	-\frac{369}{625}\\
	&\ge
	\varepsilon^2t^2
	\left(
	\frac{24}{5}\frac9{25}
	-\frac{13609}{10000}
	\right)
	+\frac{256}{125}\frac9{25}
	-\frac{369}{625}\\
	&=
	\frac{3671}{10000}\varepsilon^2t^2
	+\frac{459}{3125}\\
	&>0.
\end{align*}
	Therefore,
	$
		\mathcal R(\varepsilon,t)
		\le\sqrt{\frac32}\varepsilon t+\frac45$, and
	$	\frac{\pi\sqrt{2/3}}{\varepsilon}
		\bigl(\mathcal R(\varepsilon,t)-1\bigr)-\pi t
		\le-\frac{\pi\sqrt{2/3}}{5\varepsilon}.
	$
	
	The definition of $\mathcal R(\varepsilon,t)$
	and another comparison of fourth powers give
	\begin{align*}
		\max\left(1,\sqrt{\frac32}\varepsilon t\right)
		\le\mathcal R(\varepsilon,t)
		\le1+\frac32\varepsilon t.
	\end{align*}
	For $t\ge24$, we have 
$\mathcal R(\varepsilon,t)
		\le
		\left(\frac1{24}+\frac1{44}\right)t<t.$
Since $t/\mathcal R(\varepsilon,t)>1$, the bound in
\eqref{eq:large-tail} is largest when $r=2$.
Using
$\mathcal R(\varepsilon,t)\ge1$ and
$\mathcal R(\varepsilon,t)\ge\sqrt{3/2}\,\varepsilon t$,
we obtain
\begin{align*}
	\int_{\substack{t\ge24\\\varepsilon t\ge3/10}}
	\frac{t^3}{\mathcal R(\varepsilon,t)^5}\dd t
	&\le
	\int_0^\infty
	\frac{t^3}{\mathcal R(\varepsilon,t)^5}\dd t\\
	&\le
	\int_0^{1/\varepsilon}t^3\dd t
	+
	\left(\frac23\right)^{5/2}\varepsilon^{-5}
	\int_{1/\varepsilon}^\infty t^{-2}\dd t\\
	&=
	\frac1{4\varepsilon^4}
	+\left(\frac23\right)^{5/2}\varepsilon^{-4}\\
	&=
	\varepsilon^{-4}
	\left(\frac14+\left(\frac23\right)^{5/2}\right).
\end{align*}
The first inequality enlarges the integration region.
The second uses the first lower bound for
$\mathcal R(\varepsilon,t)$ when $t\le1/\varepsilon$
and the second when $t\ge1/\varepsilon$.

Moreover, the function
$\varepsilon^{-4}
\exp(-\pi\sqrt{2/3}/(5\varepsilon))$
is increasing for $0<\varepsilon\le1/66$,
because its derivative has the sign of
$\pi\sqrt{2/3}/5-4\varepsilon>0$.
Hence the contribution to \eqref{eq:large-tail}
from $\varepsilon t\ge3/10$ is at most
\begin{align*}
	3\varepsilon^{-4}
	\left(\frac14+\left(\frac23\right)^{5/2}\right)
	\exp\!\left(-\frac{\pi\sqrt{2/3}}{5\varepsilon}\right)
	\le
	3\cdot66^4
	\left(\frac14+\left(\frac23\right)^{5/2}\right)
	\exp\!\left(-\frac{66\pi\sqrt{2/3}}5\right)
	<
	\frac{25\cdot66^4}{12\cdot20^{11}}.
\end{align*}
The last inequality uses
$66\pi\sqrt{2/3}/5>33$, $e^3>20$, and
$(2/3)^{5/2}<4/9$.
Combining this bound with the estimate for
$\varepsilon t\le3/10$ gives
\begin{align*}
	\mathcal E_{r,2}(m,n)
	<
	\frac{27000}{20^{14}}
	+\frac{25\cdot66^4}{12\cdot20^{11}}
	<\frac1{10^6}.
\end{align*}

\noindent
\textbf{Step 3: estimating $\mathcal E_{r,3}(m,n)$.}

We first bound the polynomial factor in \eqref{eq:Fr}.
Since $|R(m,n;t,u)|\ge1$, we have
\begin{align*}
	\left|
	\frac{\sqrt6\varepsilon}{2\pi R(m,n;t,u)}
	\right|
	\le\frac{\sqrt6\varepsilon}{2\pi}
	<\frac{37}{6250}.
\end{align*}
Thus \eqref{eq:tail-H-bound} bounds the numerator
by $1061/1000$, while \eqref{eq:tail-H0-bound}
bounds the denominator factor $H_0$ from below
by $491/500$.
Consequently,
\begin{align*}
	\frac{
		\left|
		H_r\!\left(
		\frac{\sqrt6\varepsilon}{2\pi R(m,n;t,u)}
		\right)
		\right|
	}{
		|R(m,n;t,u)|^{r+3}
		H_0\!\left(\frac{\sqrt6\varepsilon}{2\pi}\right)
	}
	\le
	\frac{1061}{982}
	<\frac{109}{100}.
\end{align*}
The exponential factor in \eqref{eq:Fr} has modulus
at most $\exp(t^2/20)$ by \eqref{eq:majorant-basic}.
Multiplying these bounds, we obtain
\begin{align}\label{eq:tail-F-pointwise}
	|F_r(m,n;t,u)|
	\le\frac{109}{100}\exp\!\left(\frac{t^2}{20}\right).
\end{align}

For $8\le t\le24$, we have
\begin{align}\label{eq:tail-F-exponent}
	\frac{t^2}{20}
	=
	\frac85t-\frac{48}{5}
	-\frac{(t-8)(24-t)}{20}
	\le\frac85t-\frac{48}{5}.
\end{align}
We now substitute \eqref{eq:tail-F-pointwise} into
the definition \eqref{eq:tail-F-integral}.
Using the bound
$1/\sinh(\pi t)<(201/100)e^{-\pi t}$
established in Step 2, and using $\int_{-1/2}^{1/2}\dd u=1$, we obtain
\begin{align*}
	\mathcal E_{r,3}(m,n)
	&\le
	\frac{109}{100}\frac{201}{100}
	\int_8^{24}
	t^{r+1}\exp\!\left(\frac{t^2}{20}-\pi t\right)\dd t\\
	&\le
	\frac{11}{5}e^{-48/5}
	\int_8^{24}t^3e^{-3t/2}\dd t\\
	&\le
	\frac{11}{5}e^{-48/5}
	\int_8^\infty t^3e^{-3t/2}\dd t.
\end{align*}
The second inequality uses \eqref{eq:tail-F-exponent},
$\pi>31/10$, $t^{r+1}\le t^3$, and
$(109/100)(201/100)<11/5$.

Applying \eqref{eq:tailint} with
$j=3$, $b=8$, and $\alpha=3/2$, we have
\begin{align*}
	\int_8^\infty t^3e^{-3t/2}\dd t
	=
	\frac{11936}{27}e^{-12}.
\end{align*}
Moreover, $e^3>20$ and $e^{3/5}>9/5$ give
$
	e^{108/5}>20^7\frac95.
$
Therefore,
\begin{align*}
	\mathcal E_{r,3}(m,n)
	\le
	\frac{11}{5}\frac{11936}{27}e^{-108/5}
	<
	\frac{11\cdot443}{9\cdot20^7}
	<\frac1{10^6}.
\end{align*}
	
\noindent
\textbf{Step 4: estimating $\mathcal E_{r,4}(m,n)$.}

We first bound the coefficients of the Taylor polynomial
$T_r(m,n;t,u)$.
By \eqref{eq:Fr-at-zero}, we have
\begin{align*}
	\left|
	\left.F_r(m,n;t,u)\right|_{\varepsilon=0}
	\right|=1.
\end{align*}
Applying the product and chain rules to \eqref{eq:Fr},
with the differentiation convention in
Lemma~\ref{lem:tails}, we obtain
\begin{align}\label{eq:tail-Fr-first-zero}
	\left.
	\frac{\partial_\varepsilon F_r(m,n;t,u)}
	{F_r(m,n;t,u)}
	\right|_{\varepsilon=0}
	=
	\frac{\pi}{\sqrt6}
	\left(
	\left(k-\tfrac12\right)t^2-iut
	+\frac{m^2t^2}{4(n-\tfrac1{24})}
	\right)+
	\frac{i(r+3)mt}{2\sqrt{n-\tfrac1{24}}}
	+\frac{\sqrt6}{2\pi}\bigl(H_r'(0)+3\bigr).
\end{align}
The formulas for $H_0,H_1,H_2$ give
$H_0'(0)=-3$, $H_1'(0)=-6$, and $H_2'(0)=-10$.
Hence $|H_r'(0)+3|\le7$.
Using \eqref{eq:domain} and $|u|\le\tfrac12$ in
\eqref{eq:tail-Fr-first-zero}, the triangle inequality gives
\begin{align}\label{eq:tail-Fr-first-bound}
	\left|
	\left.\partial_\varepsilon F_r(m,n;t,u)
	\right|_{\varepsilon=0}
	\right|
	&\le
	\frac{175}{64}
	+\left(\frac{257}{400}+\frac{101}{20}\right)t+
	\frac{257}{200}
	\left(\frac32+\frac{10201}{10000}\right)t^2
	\notag\\
	&\le3+6t+4t^2.
\end{align}
Here we used
$\pi/\sqrt6<257/200$ and $\sqrt6/(2\pi)<25/64$.

Lemma~\ref{lem:majorants} also gives
\begin{align}\label{eq:par-FrU}
	\left|
	\left.\partial_\varepsilon^2F_r(m,n;t,u)
	\right|_{\varepsilon=0}
	\right|
	\le U(t).
\end{align}
Substituting~\eqref{eq:tail-Fr-first-bound} and~\eqref{eq:par-FrU} into
\eqref{eq:tail-Taylor-polynomial}, and using
$\varepsilon\le1/66$, we get
\begin{align}\label{eq:tail-Tr-polynomial-bound}
	|T_r(m,n;t,u)|
	\le
	1+\frac{3+6t+4t^2}{66}
	+\frac{U(t)}{2\cdot66^2}.
\end{align}

We next simplify this bound for $t\ge8$.
For $0\le j\le4$, we have
$t^j\le8^{j-4}t^4$.
Since all coefficients on the right of
\eqref{eq:tail-Tr-polynomial-bound} are nonnegative,
applying~\eqref{eq:tail-Tr-polynomial-bound} to each monomial yields
\begin{align}\label{eq:tail-Tr-growth}
	|T_r(m,n;t,u)|
	\le
	\frac{t^4}{8^4}
	\left(
	1+\frac{307}{66}
	+\frac{86078}{2\cdot66^2}
	\right)<\frac{t^4}{250}.
\end{align}
Here $3+6\cdot8+4\cdot8^2=307$ and $U(8)=86078$.

We now substitute \eqref{eq:tail-Tr-growth} into
the definition \eqref{eq:tail-T-integral}.
Since $\int_{-1/2}^{1/2}\dd u=1$, we have
\begin{align*}
	\mathcal E_{r,4}(m,n)
	\le
	\frac1{250}\int_8^\infty
	\frac{t^{r+5}}{\sinh(\pi t)}\dd t.
\end{align*}
The bound from Step 2 implies
$1/\sinh(\pi t)<3e^{-3t}$ for $t\ge8$.
Moreover, $r\le2$ gives $t^{r+5}\le t^7$ on this interval.
Therefore,
\begin{align*}
	\mathcal E_{r,4}(m,n)
	\le
	\frac3{250}\int_8^\infty t^7e^{-3t}\dd t.
\end{align*}

Applying \eqref{eq:tailint} with
$j=7$, $b=8$, and $\alpha=3$, we obtain
\begin{align}\label{eq:tail-T-moment}
	\int_8^\infty t^7e^{-3t}\dd t
	=
	e^{-24}\sum_{\ell=0}^7
	\frac{7!\,8^\ell}{\ell!\,3^{8-\ell}}=
	\frac{704609456}{729}e^{-24}
	<10^6e^{-24}.
\end{align}
Finally, $e^{24}=(e^3)^8>20^8$.
Using \eqref{eq:tail-T-moment}, we conclude that
\begin{align*}
	\mathcal E_{r,4}(m,n)
	<
	\frac{3\cdot10^6}{250}e^{-24}
	<
	\frac{3\cdot10^6}{250\cdot20^8}
	<
	\frac1{10^6}.
\end{align*}

Combining \eqref{eq:tails-decomposition} with the four
estimates
$\mathcal E_{r,j}(m,n)<1/10^6$ for $j=1,2,3,4$,
and then applying \eqref{eq:tail-short-remainder},
proves \eqref{eq:tails-comparison}, since
$4/10^6<1/10^5$.
\end{proof}

Fix $k\in\{1,2\}$ and $\varepsilon\ge0$.
For real $z$, define
\begin{align}
	f(z)&=\frac14\sech^2(z/2),
	\label{eq:Phi-f}\\
	a_k(z)&=\frac{\pi(2k-1)}{2\sqrt6}
	+\frac{\sqrt6z^2}{4\pi}.
	\label{eq:Phi-ak}
\end{align}

Using $f(z)$ and $a_k(z)$, define
\begin{align}
	g(z)
	=&-a_k(z)f''(z)
	-\frac{3\sqrt6z}{2\pi}f'(z),
	\label{eq:Phi-g}\\
	\ell(z)
	=\frac{a_k(z)^2}{2}f^{(4)}(z)
	+\frac{2\sqrt6z\,a_k(z)}{\pi}f'''(z)&+
	\left(
	\frac{3(2k-1)}4+\frac{45z^2}{4\pi^2}
	+\frac{\pi^2}{144}
	\right)f''(z)
	+\frac{9z}{2\pi^2}f'(z).
	\label{eq:Phi-ell}
\end{align}

We then combine the three terms by setting
\begin{align}\label{eq:Phi}
	\Phi(z)=f(z)+\varepsilon g(z)
	+\varepsilon^2\ell(z).
\end{align}

\begin{proposition}\label{prop:expansion}
	Let $k\in\{1,2\}$, and let real $n>\tfrac1{24}$
	and $m\ge0$ satisfy \eqref{eq:domain},
	with $\varepsilon$ defined by
	\eqref{eq:central-parameters}.
	Set
	\begin{align}\label{eq:expansion-z}
		z=\frac{\pi m}{\sqrt{6(n-\tfrac1{24})}}.
	\end{align}
	Then
	\begin{align}\label{eq:errors}
		|I_0(m,n)-\Phi(z)|<\frac7{20000},\quad
		|I_1(m,n)-\Phi'(z)|<\frac1{1020},\quad
		|I_2(m,n)-\Phi''(z)|<\frac1{305}.
	\end{align}
\end{proposition}

\begin{proof}
	The proof proceeds in four steps.
	First, we establish an integral representation of $f(z)$.
	Second, we integrate the Taylor polynomial
	$T_0(m,n;t,u)$ to obtain $\Phi(z)$.
	Third, differentiation with respect to $z$ gives
	the corresponding integral formulas for
	$\Phi'(z)$ and $\Phi''(z)$.
	Finally, we estimate the integrated remainders
	to prove the three inequalities in \eqref{eq:errors}.
	
\noindent
\textbf{Step 1: an integral representation of $f(z)$.}
We first establish \eqref{eq:logisticfourier}
using a Fourier cosine series.
Throughout the initial calculation, $z>0$ is fixed.

For every positive integer $j$, integration by parts gives
\begin{align}\label{eq:cosine-ibp-first}
	\int_0^1\cosh(zx)\cos(j\pi x)\dd x=
	\frac{(-1)^j\sinh z}{z}
	+\frac{j\pi}{z}
	\int_0^1\sinh(zx)\sin(j\pi x)\dd x.
\end{align}
A second integration by parts yields
\begin{align}\label{eq:cosine-ibp-second}
	&\int_0^1\sinh(zx)\sin(j\pi x)\dd x=
	-\frac{j\pi}{z}
	\int_0^1\cosh(zx)\cos(j\pi x)\dd x,
\end{align}
because $\sin(j\pi x)=0$ at $x=0,1$.
Substituting \eqref{eq:cosine-ibp-second}
into \eqref{eq:cosine-ibp-first}, we obtain
\begin{align}\label{eq:cosh-cosine-coefficients}
	\int_0^1\cosh(zx)\cos(j\pi x)\dd x
	=
	\frac{z(-1)^j\sinh z}{z^2+j^2\pi^2}.
\end{align}
The constant coefficient follows directly from
$
	\int_0^1\cosh(zx)\dd x=\frac{\sinh z}{z}.
$

For fixed $z>0$, the function $\cosh(zx)$ is
continuously differentiable with respect to $x$
for $0\le x\le1$.
The convergence theorem for Fourier cosine series
therefore gives convergence to $\cosh(zx)$
for every $0\le x\le1$, including both endpoints.
Substituting the coefficients from
\eqref{eq:cosh-cosine-coefficients}, we obtain
\begin{align}\label{eq:cosh-cosine-series}
	\cosh(zx)
	=
	\frac{\sinh z}{z}
	+
	2z\sinh z
	\sum_{j=1}^{\infty}
	\frac{(-1)^j\cos(j\pi x)}
	{z^2+j^2\pi^2},
	\qquad 0\le x\le1.
\end{align}
Moreover, the coefficients are $O(j^{-2})$,
so the series converges absolutely and uniformly in $x$.

Setting $x=1$ in \eqref{eq:cosh-cosine-series},
we have
\begin{align}\label{eq:cosh-series-at-one}
	\cosh z
	=
	\frac{\sinh z}{z}
	+
	2z\sinh z
	\sum_{j=1}^{\infty}
	\frac1{z^2+j^2\pi^2}.
\end{align}
Setting $x=0$ gives
\begin{align}\label{eq:cosh-series-at-zero}
	1
	=
	\frac{\sinh z}{z}
	+
	2z\sinh z
	\sum_{j=1}^{\infty}
	\frac{(-1)^j}{z^2+j^2\pi^2}.
\end{align}
Subtracting \eqref{eq:cosh-series-at-zero}
from \eqref{eq:cosh-series-at-one}, we obtain
\begin{align*}
	\cosh z-1
	=
	2z\sinh z
	\sum_{j=1}^{\infty}
	\frac{1-(-1)^j}{z^2+j^2\pi^2}
	=
	4z\sinh z
	\sum_{j=0}^{\infty}
	\frac1{z^2+(2j+1)^2\pi^2}.
\end{align*}
Since
$
	\frac{\cosh (z)-1}{\sinh z}
	=
	\frac{2\sinh^2(z/2)}
	{2\sinh(z/2)\cosh(z/2)}
	=
	\tanh(z/2),
$
division by $\sinh z$ gives
\begin{align}\label{eq:tanh-partial-fractions}
	\tanh(z/2)
	=
	4\sum_{j=0}^{\infty}
	\frac{z}{z^2+(2j+1)^2\pi^2}.
\end{align}

We next differentiate
\eqref{eq:tanh-partial-fractions}.
The series converges uniformly on bounded $z$-intervals,
and
\begin{align*}
	\left|
	\frac{\mathrm d}{\mathrm dz}
	\frac{z}{z^2+(2j+1)^2\pi^2}
	\right|
	=
	\frac{\left|(2j+1)^2\pi^2-z^2\right|}
	{\bigl(z^2+(2j+1)^2\pi^2\bigr)^2}
	\le
	\frac1{(2j+1)^2\pi^2}.
\end{align*}
The bounds form a convergent series.
Thus termwise differentiation is valid.
Since
$f(z)=\tfrac14\sech^2(z/2)
=\tfrac12\frac{\mathrm d}{\mathrm dz}\tanh(z/2)$,
we obtain
\begin{align}\label{eq:f-partial-fractions}
	f(z)
	=\tfrac14\sech^2(z/2)=
	2\sum_{j=0}^{\infty}
	\frac{(2j+1)^2\pi^2-z^2}
	{\bigl((2j+1)^2\pi^2+z^2\bigr)^2}.
\end{align}

We now obtain the same series from the required integral.
For $t>0$, the geometric series gives
\begin{align}\label{eq:sinh-geometric-series}
	\frac1{\sinh(\pi t)}
	=
	\frac{2e^{-\pi t}}{1-e^{-2\pi t}}
	=
	2\sum_{j=0}^{\infty}e^{-(2j+1)\pi t}.
\end{align}
By $|\cos(zt)|\le1$ and \eqref{eq:tailint}, we have 
\begin{align}\label{eq:cosine-series-integrable}
	2\sum_{j=0}^{\infty}
	\int_0^\infty
	t e^{-(2j+1)\pi t}|\cos(zt)|\dd t\le
	2\sum_{j=0}^{\infty}
	\int_0^\infty
	t e^{-(2j+1)\pi t}\dd t=
	\frac2{\pi^2}
	\sum_{j=0}^{\infty}\frac1{(2j+1)^2}
	<\infty.
\end{align}
Consequently, the absolute convergence in
\eqref{eq:cosine-series-integrable}
allows us to interchange the series and integral
after using \eqref{eq:sinh-geometric-series}.

For $\alpha>0$, direct integration gives
\begin{align}\label{eq:laplace-cosine}
	\int_0^\infty e^{-\alpha t}\cos(zt)\dd t
	=\frac{\alpha}{\alpha^2+z^2}.
\end{align}
By integration by parts and \eqref{eq:laplace-cosine},
we have
\begin{align}\label{eq:laplace-sine}
	\int_0^\infty e^{-\alpha t}\sin(zt)\dd t
	&=
	\frac{z}{\alpha}
	\int_0^\infty e^{-\alpha t}\cos(zt)\dd t=\frac{z}{\alpha^2+z^2}.
\end{align}

Integrating by parts using \eqref{eq:laplace-cosine} and  \eqref{eq:laplace-sine}, 
we obtain
\begin{align}\label{eq:laplace-cosine-first-moment}
	\int_0^\infty t e^{-\alpha t}\cos(zt)\dd t
	&=
	\frac1{\alpha^2+z^2}
	\Biggl(
	\alpha\int_0^\infty
	e^{-\alpha t}\cos(zt)\dd t
	-z\int_0^\infty
	e^{-\alpha t}\sin(zt)\dd t
	\Biggr)
	\notag\\
	&=
	\frac1{\alpha^2+z^2}
	\left(
	\frac{\alpha^2}{\alpha^2+z^2}
	-\frac{z^2}{\alpha^2+z^2}
	\right)
	\notag\\
	&=
	\frac{\alpha^2-z^2}{(\alpha^2+z^2)^2}.
\end{align}

Applying \eqref{eq:laplace-cosine-first-moment}
with $\alpha=(2j+1)\pi$, and then using
\eqref{eq:f-partial-fractions}, we conclude that
\begin{align}\label{eq:logisticfourier}
	\int_0^\infty
	\frac{t\cos(zt)}{\sinh(\pi t)}\dd t
	=
	2\sum_{j=0}^{\infty}
	\int_0^\infty
	t e^{-(2j+1)\pi t}\cos(zt)\dd t=
	2\sum_{j=0}^{\infty}
	\frac{(2j+1)^2\pi^2-z^2}
	{\bigl((2j+1)^2\pi^2+z^2\bigr)^2}=f(z).
\end{align}

So far, \eqref{eq:logisticfourier}
has been proved for $z>0$.
The bound
$
	\left|
	\frac{t\cos(zt)}{\sinh(\pi t)}
	\right|
	\le
	\frac{t}{\sinh(\pi t)}
$
is integrable and independent of $z$.
Dominated convergence therefore shows that the integral
is continuous in $z$.
Since both sides of \eqref{eq:logisticfourier}
are even, the identity holds for every real $z$.

Finally, for every nonnegative integer $j$,
every real $z$, and every $t>0$, we have
\begin{align}\label{eq:logistic-derivative-bound}
	\left|
	\frac{\partial^j}{\partial z^j}
	\frac{t\cos(zt)}{\sinh(\pi t)}
	\right|
	\le
	\frac{t^{j+1}}{\sinh(\pi t)}.
\end{align}
Since $t^{j+1}/\sinh(\pi t)$ in
\eqref{eq:logistic-derivative-bound}
is $O(t^j)$ as $t\to0^+$
and decays exponentially as $t\to\infty$,
it is integrable on $0<t<\infty$.
Repeated differentiation under the integral sign
in \eqref{eq:logisticfourier} consequently gives
\begin{align}\label{eq:logisticfourier-derivatives}
	f^{(j)}(z)
	=
	\int_0^\infty
	\frac{t^{j+1}}{\sinh(\pi t)}
	\Re\!\left((-i)^j e^{-izt}\right)\dd t.
\end{align}
	
\noindent
\textbf{Step 2: integrating the Taylor polynomial
	$T_0(m,n;t,u)$.}
We first compute the Taylor coefficients of
$F_0(m,n;t,u)$, then integrate with respect to $u$
and $t$.
All expansions are at $\varepsilon=0$,
with $k,z,t,u$ fixed.
By \eqref{eq:expansion-z}, fixing $z$ is equivalent
to fixing $m/\sqrt{n-\tfrac1{24}}$, as required
in \eqref{eq:tail-Taylor-polynomial}.

Since $z=\pi m/\sqrt{6(n-\tfrac1{24})}$
and $\varepsilon=1/\sqrt{n-\tfrac1{24}}$,
we have $m=\sqrt6z/(\pi\varepsilon)$.
Substituting $m=\sqrt6z/(\pi\varepsilon)$ and
$n-\tfrac1{24}=\varepsilon^{-2}$
into \eqref{eq:R-definition}, we obtain
\begin{align}\label{eq:expansion-R-square}
	R(m,n;t,u)^2
	=
	1-\frac{i\sqrt6zt}{\pi}\varepsilon
	+\left(
	\left(k-\tfrac12\right)t^2-iut
	\right)\varepsilon^2.
\end{align}
Applying $	\sqrt{1+x}
=1+\frac{x}{2}-\frac{x^2}{8}
+\frac{x^3}{16}+O(x^4)$ to
\eqref{eq:expansion-R-square}, and using
$a_k(z)$ from \eqref{eq:Phi-ak}, we obtain
\begin{align}\label{eq:expansion-R-series}
	R(m,n;t,u)
	=&
	1-\frac{i\sqrt6zt}{2\pi}\varepsilon
	+\left(
	\frac{\sqrt6a_k(z)t^2}{2\pi}
	-\frac{iut}{2}
	\right)\varepsilon^2+
	\left(
	\frac{3iz\,a_k(z)t^3}{2\pi^2}
	+\frac{\sqrt6zu t^2}{4\pi}
	\right)\varepsilon^3
	+O(\varepsilon^4).
\end{align}

We consider the exponential factor in
\eqref{eq:Fr} with $r=0$.
By \eqref{eq:R-definition} and
\eqref{eq:central-parameters}, we have
\begin{align}\label{eq:expansion-R-difference}
	R(m,n;t,u)^2-1
	=
	\varepsilon^2t
	\left(
	\left(k-\tfrac12\right)t-i(m+u)
	\right).
\end{align}
Using \eqref{eq:expansion-R-difference}
and $R^2-1=(R-1)(R+1)$,
we rewrite the exponent in \eqref{eq:Fr} as
\begin{align}\label{eq:expansion-exponent-rewrite}
	\frac{\pi\sqrt{2/3}\,\varepsilon t}
	{R(m,n;t,u)+1}
	\left(
	\left(k-\tfrac12\right)t-i(m+u)
	\right)&=
	\frac{\pi\sqrt{2/3}}{\varepsilon}
	\frac{R(m,n;t,u)^2-1}{R(m,n;t,u)+1}
	\notag\\
	&=
	\frac{\pi\sqrt{2/3}}{\varepsilon}
	\bigl(R(m,n;t,u)-1\bigr).
\end{align}
Division by $\varepsilon$ lowers each power of
$\varepsilon$ by one.
Thus the cubic term in \eqref{eq:expansion-R-series}
is needed to determine the quadratic term
of the exponent.
Substituting \eqref{eq:expansion-R-series}
into \eqref{eq:expansion-exponent-rewrite}, we obtain
\begin{align}\label{eq:expansion-exponent-series}
	\frac{\pi\sqrt{2/3}}{\varepsilon}
	\bigl(R(m,n;t,u)-1\bigr)
	={}&
	-izt+
	\left(
	a_k(z)t^2-\frac{i\pi ut}{\sqrt6}
	\right)\varepsilon
	\notag\\
	&+
	\left(
	\frac{i\sqrt6z\,a_k(z)t^3}{2\pi}
	+\frac{zu t^2}{2}
	\right)\varepsilon^2
	+O(\varepsilon^3).
\end{align}

For the factor $R(m,n;t,u)^{-3}$, we apply
$(1+x)^{-3/2}=1-3x/2+15x^2/8+O(x^3)$
to \eqref{eq:expansion-R-square}.
This gives
\begin{align}\label{eq:expansion-reciprocal-factor}
	R(m,n;t,u)^{-3}
	=
	1+\frac{3i\sqrt6zt}{2\pi}\varepsilon+
	\left\{
	\frac{3iut}{2}
	-
	\left(
	\frac{3(2k-1)}4
	+\frac{45z^2}{4\pi^2}
	\right)t^2
	\right\}\varepsilon^2
	+O(\varepsilon^3).
\end{align}

To expand $F_0(m,n;t,u)$ in \eqref{eq:Fr},
we next consider the factor
\begin{align*}
	\frac{
		H_0\!\left(
		\frac{\sqrt6\varepsilon}{2\pi R(m,n;t,u)}
		\right)
	}{
		H_0\!\left(\frac{\sqrt6\varepsilon}{2\pi}\right)
	}.
\end{align*}
By \eqref{eq:expansion-R-series}, we have
\begin{align}\label{eq:expansion-H-input}
	\frac{\sqrt6\varepsilon}{2\pi R(m,n;t,u)}
	=
	\frac{\sqrt6\varepsilon}{2\pi}
	+\frac{3izt}{2\pi^2}\varepsilon^2
	+O(\varepsilon^3).
\end{align}
Since $H_0(v)=1-3v+3v^2$, we have 
\begin{align*}
	H_0\!\left(
	\frac{\sqrt6\varepsilon}{2\pi R(m,n;t,u)}
	\right)
	=
	H_0\!\left(\frac{\sqrt6\varepsilon}{2\pi}\right)
	-\frac{9izt}{2\pi^2}\varepsilon^2
	+O(\varepsilon^3).
\end{align*}
Dividing by
$H_0(\sqrt6\varepsilon/(2\pi))=1+O(\varepsilon)$,
we obtain
\begin{align}\label{eq:expansion-polynomial-factor}
	\frac{
		H_0\!\left(
		\frac{\sqrt6\varepsilon}{2\pi R(m,n;t,u)}
		\right)
	}{
		H_0\!\left(\frac{\sqrt6\varepsilon}{2\pi}\right)
	}
	=
	1-\frac{9izt}{2\pi^2}\varepsilon^2
	+O(\varepsilon^3).
\end{align}

We now combine the exponent expansion
\eqref{eq:expansion-exponent-series}
with the factor expansions
\eqref{eq:expansion-reciprocal-factor} and
\eqref{eq:expansion-polynomial-factor}.
Using \eqref{eq:expansion-exponent-series},
we first factor out $e^{-izt}$ from the exponential
factor in \eqref{eq:Fr} with $r=0$.
We expand the remaining exponential using
$e^x=1+x+x^2/2+O(x^3)$.
Multiplying by the expansions in
\eqref{eq:expansion-reciprocal-factor} and
\eqref{eq:expansion-polynomial-factor},
and retaining powers through $\varepsilon^2$,
gives $T_0(m,n;t,u)$ by
\eqref{eq:tail-Taylor-polynomial}.

For integration with respect to $u$, we use
\begin{align}\label{eq:expansion-u-moments}
	\int_{-1/2}^{1/2}\dd u=1,\qquad
	\int_{-1/2}^{1/2}u\dd u=0,\qquad
	\int_{-1/2}^{1/2}u^2\dd u=\frac1{12}.
\end{align}
All terms linear in $u$ vanish after integration.
After factoring out $e^{-izt}$, the only $u^2$ term
at degree two in $\varepsilon$ comes from
the squared linear term of the exponential.
Its integral is
\begin{align*}
	\frac12\int_{-1/2}^{1/2}
	\left(-\frac{i\pi ut}{\sqrt6}\right)^2\dd u
	=-\frac{\pi^2t^2}{144}.
\end{align*}

After integration with respect to $u$,
the coefficient of $t^3$ at degree two
in $\varepsilon$ is
\begin{align*}
	\frac{i\sqrt6z\,a_k(z)}{2\pi}
	+\frac{3i\sqrt6z\,a_k(z)}{2\pi}
	=
	\frac{2i\sqrt6z\,a_k(z)}{\pi}.
\end{align*}
The first summand comes from the quadratic term
in \eqref{eq:expansion-exponent-series}.
The second comes from multiplying the linear terms
of the exponential and
\eqref{eq:expansion-reciprocal-factor}.
The $t^4$ term comes from squaring $a_k(z)t^2$
in the exponential expansion.
The $t^2$ term combines
\eqref{eq:expansion-reciprocal-factor}
with the contribution $-\pi^2t^2/144$.
The term linear in $t$ comes from
\eqref{eq:expansion-polynomial-factor}.

By \eqref{eq:tail-Taylor-polynomial},
we therefore obtain the exact identity
\begin{align}\label{eq:expansion-averaged-T0}
	\int_{-1/2}^{1/2}T_0(m,n;t,u)\dd u
	&=
	e^{-izt}\Biggl\{
	1+\varepsilon
	\left(
	a_k(z)t^2+\frac{3i\sqrt6z}{2\pi}t
	\right)+
	\varepsilon^2\Biggl(
	\frac{a_k(z)^2}{2}t^4
	+\frac{2i\sqrt6z\,a_k(z)}{\pi}t^3
	\notag\\
	&\hspace{18mm}
	-
	\left(
	\frac{3(2k-1)}4
	+\frac{45z^2}{4\pi^2}
	+\frac{\pi^2}{144}
	\right)t^2
	-\frac{9iz}{2\pi^2}t
	\Biggr)\Biggr\}.
\end{align}

It remains to integrate
\eqref{eq:expansion-averaged-T0}
with respect to $t$.
Since $\Re(e^{-izt})=\cos(zt)$ and
$\Re(ie^{-izt})=\sin(zt)$,
\eqref{eq:logisticfourier-derivatives} gives
\begin{align}\label{eq:expansion-Fourier-moments}
	\int_0^\infty
	\frac{t^2\sin(zt)}{\sinh(\pi t)}\dd t
	&=-f'(z),&
	\int_0^\infty
	\frac{t^3\cos(zt)}{\sinh(\pi t)}\dd t
	&=-f''(z),
	\notag\\
	\int_0^\infty
	\frac{t^4\sin(zt)}{\sinh(\pi t)}\dd t
	&=f'''(z),&
	\int_0^\infty
	\frac{t^5\cos(zt)}{\sinh(\pi t)}\dd t
	&=f^{(4)}(z).
\end{align}
Applying \eqref{eq:expansion-Fourier-moments}
to the linear and quadratic coefficients in
\eqref{eq:expansion-averaged-T0}
yields $g(z)$ and $\ell(z)$ as defined in
\eqref{eq:Phi-g} and \eqref{eq:Phi-ell}.
The constant coefficient integrates to $f(z)$
by \eqref{eq:logisticfourier-derivatives} with $j=0$.
Therefore
\begin{align}\label{eq:expansion-T0-integral}
	&\int_0^\infty\frac{t}{\sinh(\pi t)}
	\int_{-1/2}^{1/2}
	\Re T_0(m,n;t,u)\dd u\dd t
	=
	f(z)+\varepsilon g(z)+\varepsilon^2\ell(z)
	=\Phi(z).
\end{align}
	
\noindent
\textbf{Step 3: integrating
	$T_1(m,n;t,u)$ and $T_2(m,n;t,u)$.}
We obtain the required integrals by differentiating
\eqref{eq:expansion-T0-integral} with respect to $z$.
We first establish the derivative relations for
$F_r(m,n;t,u)$ and $T_r(m,n;t,u)$.

Fix $\varepsilon>0$ and keep $k,t,u$ fixed.
By \eqref{eq:central-parameters} and
\eqref{eq:expansion-z}, we have
\begin{align*}
	n=\tfrac1{24}+\varepsilon^{-2},\qquad
	m=\frac{\sqrt6z}{\pi\varepsilon},\qquad
	\frac{\partial m}{\partial z}
	=\frac{\sqrt6}{\pi\varepsilon}.
\end{align*}
Thus $n$ is independent of $z$, while $m$ varies
with $z$.
Before differentiating $F_r(m,n;t,u)$ or
$T_r(m,n;t,u)$ with respect to $z$,
we substitute $m=\sqrt6z/(\pi\varepsilon)$.

The expression $n+(k-\tfrac12)t^2-i(m+u)t$
is linear in $z$, with derivative
$-i\sqrt6t/(\pi\varepsilon)$.
Hence the chain rule and \eqref{eq:Fr-ratio} give
\begin{align}\label{eq:expansion-z-derivatives}
	\frac{\partial^r}{\partial z^r}
	F_0(m,n;t,u)
	&=
	\left(-\frac{i\sqrt6t}{\pi\varepsilon}\right)^r
	\frac{
		P^{(r+1)}\!\left(
		n+\left(k-\tfrac12\right)t^2-i(m+u)t
		\right)
	}{
		P'(n)
	}
	\notag\\
	&=(-it)^rF_r(m,n;t,u),
	\qquad r=0,1,2.
\end{align}

We next justify differentiation of the Taylor
coefficients at $\varepsilon=0$.
Equation \eqref{eq:expansion-R-square} shows that
$R(m,n;t,u)$ is smooth in $(\varepsilon,z)$
near $\varepsilon=0$, where its value is $1$.
Moreover, rationalizing $R(m,n;t,u)-1$ gives
\begin{align}\label{eq:expansion-regular-exponent}
	&\frac{\pi\sqrt{2/3}}{\varepsilon}
	\bigl(R(m,n;t,u)-1\bigr)=
	\frac{
		-2izt+\pi\sqrt{2/3}\,\varepsilon
		\left(
		\left(k-\tfrac12\right)t^2-iut
		\right)
	}{
		R(m,n;t,u)+1
	}.
\end{align}
The right side of
\eqref{eq:expansion-regular-exponent}
is smooth at $\varepsilon=0$,
since its denominator then equals $2$.
The remaining denominators in \eqref{eq:Fr}
have values $R=1$ and $H_0(0)=1$ there.
Consequently, \eqref{eq:Fr} extends
$F_r(m,n;t,u)$ smoothly in $(\varepsilon,z)$
near $\varepsilon=0$.

By \eqref{eq:expansion-z}, differentiation with respect
to $\varepsilon$ in \eqref{eq:tail-Taylor-polynomial}
is performed with $z$ fixed.
The mixed derivatives therefore commute.
For $j,r\in\{0,1,2\}$,
\eqref{eq:expansion-z-derivatives} yields
\begin{align}\label{eq:expansion-mixed-coefficients}
	\frac{\partial^r}{\partial z^r}
	\left(
	\left.
	\partial_\varepsilon^jF_0(m,n;t,u)
	\right|_{\varepsilon=0}
	\right)
	=
	\left.
	\partial_\varepsilon^j\partial_z^rF_0(m,n;t,u)
	\right|_{\varepsilon=0}=
	(-it)^r
	\left.
	\partial_\varepsilon^jF_r(m,n;t,u)
	\right|_{\varepsilon=0}.
\end{align}
Multiplying \eqref{eq:expansion-mixed-coefficients}
by $\varepsilon^j/j!$ and summing over $j=0,1,2$,
we obtain
\begin{align}\label{eq:expansion-T-derivatives}
	\frac{\partial^r}{\partial z^r}
	T_0(m,n;t,u)
	=
	(-it)^rT_r(m,n;t,u),
	\qquad r=0,1,2,
\end{align}
by the definition \eqref{eq:tail-Taylor-polynomial}.

It remains to justify differentiating
\eqref{eq:expansion-T0-integral}
under the integral sign.
Differentiation commutes with integration over
$-\tfrac12\le u\le\tfrac12$,
because $T_0(m,n;t,u)$ and its first two derivatives
with respect to $z$ are jointly continuous in $(z,u)$.

By \eqref{eq:expansion-averaged-T0},
the integral with respect to $u$ is $e^{-izt}$
times a polynomial of degree at most four in $t$.
After differentiating twice with respect to $z$,
the expression is $e^{-izt}$ times a polynomial
of degree at most six in $t$.
The coefficients are bounded when $z$ varies
over a bounded interval.
Also, $|e^{-izt}|=1$ for real $z$.
Thus, for every $M>0$, there is a constant $C_M>0$
such that
\begin{align}\label{eq:expansion-integrable-derivatives}
	\frac{t}{\sinh(\pi t)}
	\left|
	\frac{\partial^r}{\partial z^r}
	\int_{-1/2}^{1/2}T_0(m,n;t,u)\dd u
	\right|
	\le
	C_M\frac{t(1+t^6)}{\sinh(\pi t)},
\end{align}
for $r=0,1,2$, $t>0$, and $|z|\le M$.
Here $C_M$ may depend on the fixed $k$
and $\varepsilon$.
The right side of
\eqref{eq:expansion-integrable-derivatives}
is bounded near zero and is
$O(t^7e^{-\pi t})$ as $t\to\infty$.
It is therefore integrable.

By dominated convergence, we may differentiate
\eqref{eq:expansion-T0-integral} twice with respect
to $z$ under the integral sign.
Using \eqref{eq:expansion-T-derivatives}, we obtain
\begin{align}\label{eq:expansion-integrated-T}
	\Phi^{(r)}(z)
	&=
	\int_0^\infty\frac{t}{\sinh(\pi t)}
	\int_{-1/2}^{1/2}
	\Re\!\left(
	\partial_z^rT_0(m,n;t,u)
	\right)\dd u\dd t
	\notag\\
	&=
	\int_0^\infty\frac{t^{r+1}}{\sinh(\pi t)}
	\int_{-1/2}^{1/2}
	\Re\!\left(
	(-i)^rT_r(m,n;t,u)
	\right)\dd u\dd t,
\end{align}
for $r=0,1,2$, where $\Phi^{(0)}=\Phi$.
	
\noindent
\textbf{Step 4: estimating the integrated remainder.}
	By \eqref{eq:expansion-integrated-T} and
	Lemma~\ref{lem:tails}, we have
	\begin{align}\label{eq:expansion-remainder}
		|I_r(m,n)-\Phi^{(r)}(z)|
		<
		\frac{\varepsilon^3}{6}
		\int_0^8
		\frac{
			t^{r+1}\exp(t^2/20)P_r(t)
		}{
			\sinh(\pi t)
		}\dd t
		+\frac1{10^5}.
	\end{align}
	
	We first replace the exponential in
	\eqref{eq:expansion-remainder} by a polynomial.
	For $0\le t\le8$ and $j\ge4$, the ratio of
	consecutive terms in its series satisfies
	$
		\frac{t^2}{20(j+1)}\le\frac{16}{25}.
$ 
	Thus 
	\begin{align*}
		\sum_{j=4}^{\infty}
		\frac{(t^2/20)^j}{j!}
		\le
		\frac{t^8}{24\cdot20^4}
		\sum_{\ell=0}^{\infty}
		\left(\frac{16}{25}\right)^\ell
		=\frac{t^8}{1382400}.
	\end{align*}
	Adding $1+t^2/20+t^4/800+t^6/48000$, we obtain 
	\begin{align}\label{eq:expansion-exp-bound}
		\exp(t^2/20)
		\le
		\Psi(t)
		:=
		1+\frac{t^2}{20}+\frac{t^4}{800}
		+\frac{t^6}{48000}
		+\frac{t^8}{1382400}.
	\end{align}
	
For every nonnegative integer $j$,
applying \eqref{eq:sinh} and then
\eqref{eq:tailint} with $b=0$, we obtain 
	\begin{align}\label{eq:expansion-monomial}
		\int_0^\infty
		\frac{t^{j+r+1}}{\sinh(\pi t)}\dd t
		&\le
		2\int_0^\infty
		t^{j+r+1}e^{-\pi t}\dd t
		+\frac1\pi\int_0^\infty
		t^{j+r}e^{-\pi t}\dd t
		\notag\\
		&=
		\frac{2(j+r+1)!+(j+r)!}
		{\pi^{j+r+2}}
		\notag\\
		&<
		\frac{2(j+r+1)!+(j+r)!}
		{(157/50)^{j+r+2}},
	\end{align}
	where the last inequality uses $\pi>157/50$.
	
	Recall from Lemma~\ref{lem:majorants} that
	$P_r(t)=\sum_{j=0}^6p_{r,j}t^j$.
	All coefficients of $\Psi(t)P_r(t)$ are nonnegative.
	Using $\Psi(t)$ from \eqref{eq:expansion-exp-bound},
	we apply \eqref{eq:expansion-monomial}
	to each monomial in $\Psi(t)P_r(t)$ and obtain
	\begin{align}\label{eq:expansion-integral-J}
		\int_0^8
		\frac{
			t^{r+1}\exp(t^2/20)P_r(t)
		}{
			\sinh(\pi t)
		}\dd t
		<
		J_r,
	\end{align}
	where
	\begin{align}\label{eq:moments}
		J_r
		={}&
		\sum_{j=0}^6p_{r,j}
		\Biggl(
		\frac{2(j+r+1)!+(j+r)!}
		{(157/50)^{j+r+2}}+
		\frac{2(j+r+3)!+(j+r+2)!}
		{20(157/50)^{j+r+4}}
		\notag\\
		&\qquad+
		\frac{2(j+r+5)!+(j+r+4)!}
		{800(157/50)^{j+r+6}}+
		\frac{2(j+r+7)!+(j+r+6)!}
		{48000(157/50)^{j+r+8}}
		\notag\\
		&\qquad+
		\frac{2(j+r+9)!+(j+r+8)!}
		{1382400(157/50)^{j+r+10}}
		\Biggr).
	\end{align}
	
	Substituting the coefficients of $P_0,P_1,P_2$
	from Lemma~\ref{lem:majorants} into
	\eqref{eq:moments}, we obtain
	\begin{align}\label{eq:expansion-J-bounds}
		J_0
		&<
		\frac{46014+9254+1841+388+244}{100}
		<578,\notag\\
		J_1
		&<
		\frac{121940+33291+8078+1986+1421}{100}
		<1668,\notag\\
		J_2
		&<
		\frac{371474+130805+38109+10853+8797}{100}
		<5601.
	\end{align}
	The five numerators in each line correspond,
	in order, to the five sums in \eqref{eq:moments}.
	Each comparison is a finite rational inequality,
	verified by multiplication by its positive
	common denominator.
	
	Finally, $\varepsilon\le1/66$.
	Combining \eqref{eq:expansion-remainder},
	\eqref{eq:expansion-integral-J}, and
	\eqref{eq:expansion-J-bounds}, we conclude that
	\begin{align*}
		|I_0(m,n)-\Phi(z)|
		&<
		\frac{578}{6\cdot66^3}+\frac1{10^5}
		<\frac7{20000},\\
		|I_1(m,n)-\Phi'(z)|
		&<
		\frac{1668}{6\cdot66^3}+\frac1{10^5}
		<\frac1{1020},\\
		|I_2(m,n)-\Phi''(z)|
		&<
		\frac{5601}{6\cdot66^3}+\frac1{10^5}
		<\frac1{305}.
	\end{align*}
	This proves \eqref{eq:errors}.
\end{proof}

\subsection{Positive curvature of the approximation}

We first prove $\Phi(z)>0$ and obtain a positive
lower bound for
$(\Phi'(z))^2-\Phi(z)\Phi''(z)$.
We then use the three inequalities in \eqref{eq:errors}
to compare $I_0(m,n),I_1(m,n),I_2(m,n)$ with
$\Phi(z),\Phi'(z),\Phi''(z)$, respectively,
and prove $I_1(m,n)^2-I_0(m,n)I_2(m,n)>0$.

\begin{lemma}\label{lem:Phi-curvature}
	Suppose that $k\in\{1,2\}$,
	$0\le\varepsilon\le1/66$, and
	$0\le z\le101\pi/(50\sqrt6)$.
	For $r=1,2$, the functions $f(z)$ and $\Phi(z)$
	defined in \eqref{eq:Phi} satisfy
	the following four inequalities
	\begin{align}\label{eq:Phibounds}
		\Phi(z)>\frac{993}{1000}f(z),\;
		|\Phi(z)|<\frac{111}{100}f(z),\;
		|\Phi^{(r)}(z)|<\frac{31}{25}f(z),\;
		-(\log\Phi)''(z)>\frac{83}{50}f(z).
	\end{align}
\end{lemma}

\begin{proof}
	We first estimate the derivatives of $f(z)$, then bound
	$g(z)$ and $\ell(z)$, and finally apply the identity
	$\Phi(z)/f(z)=1+\varepsilon g(z)/f(z)+\varepsilon^2\ell(z)/f(z)$.
	All derivatives are taken with respect to $z$.
	In the calculations, functions with suppressed variables
	are evaluated at the same real number $z$.
	
\noindent
\textbf{Step 1: estimates for $f(z)$.}
Differentiating $f(z)=\tfrac14\sech^2(z/2)$, we have
\begin{align}\label{eq:Phi-short-basic}
	f'(z)&=-f(z)\tanh(z/2),&
	\frac{\mathrm d}{\mathrm dz}\tanh(z/2)&=2f(z),&
	\tanh^2(z/2)&=1-4f(z).
\end{align}
Thus $f(z)$ is positive and decreasing for $z\ge0$,
with $f(z)\le f(0)=1/4$.
The assumed range of $z$ and $25/8<\pi<22/7$ also give
$0\le z<649/250$ and $5/2<\pi\sqrt{2/3}<13/5$.
We use the elementary bounds
\begin{align}\label{eq:Phi-short-exp-explicit}
	12<e^{649/250}<\frac{671}{50}.
\end{align}
Since $f(z)=e^z/(1+e^z)^2$ and $x/(1+x)^2$ decreases
for $x\ge1$, we obtain
\begin{align}
	\frac2{31}<&f(z)\le \frac14,
	\label{eq:f-range}\\
	0\le\tanh^2(z/2)&=1-4f(z)<\frac34.
	\label{eq:Phi-short-tanh-range}
\end{align}

We next compare $z$ with $\tanh(z/2)$.
For $z>0$, by \eqref{eq:Phi-short-basic} and
the monotonicity of $f(z)$, we have
\begin{align*}
	\tanh(z/2)
	&=2\int_0^z f(s)\,\mathrm ds\le\frac z2,\\
	\frac{\mathrm d}{\mathrm dz}\frac{z}{\tanh(z/2)}
	&=
	\frac{2\int_0^z\bigl(f(s)-f(z)\bigr)\,\mathrm ds}
	{\tanh^2(z/2)}
	\ge0.
\end{align*}
Moreover, \eqref{eq:Phi-short-exp-explicit} gives
$\tanh(649/500)>11/13$.
By the monotonicity just proved, we get
$z/\tanh(z/2)<(649/250)(13/11)<31/10$ for $z>0$.
Consequently,
\begin{align}\label{eq:Phi-short-zT}
	2\tanh(z/2)\le z\le\frac{31}{10}\tanh(z/2).
\end{align}
Both inequalities also hold at $z=0$.

We now estimate the derivatives of $f(z)$.
Differentiating \eqref{eq:Phi-short-basic}, we obtain
$f''(z)=f(z)(1-6f(z))$, and then
$f'''(z)=f'(z)(1-12f(z))$.
We therefore obtain
\begin{align}\label{eq:Phi-short-q2}
	\frac{f''(z)}{f(z)}=1-6f(z).
\end{align}
Differentiating $f''(z)=f(z)(1-6f(z))$ gives
\begin{align}\label{eq:Phi-short-q3}
	\frac{f'''(z)}{f(z)}
	=
	\frac{f'(z)}{f(z)}\bigl(1-12f(z)\bigr)
	=
	\tanh(z/2)\bigl(12f(z)-1\bigr).
\end{align}
Using
\eqref{eq:Phi-short-basic} and \eqref{eq:Phi-short-q2},
we obtain
\begin{align}\label{eq:Phi-short-q4}
	\frac{f^{(4)}(z)}{f(z)}
	=
	\frac{f''(z)}{f(z)}\bigl(1-12f(z)\bigr)
	-12\frac{f'(z)^2}{f(z)}=1-30f(z)+120f(z)^2.
\end{align}
We first estimate the derivatives of $f(z)$ relative to $f(z)$.
By \eqref{eq:f-range}, we have $1/16<f(z)\le1/4$.

For the first two derivatives, \eqref{eq:Phi-short-basic}
and \eqref{eq:Phi-short-q2} give
\begin{align}\label{eq:Phi-short-low12}
	\frac{|f'(z)|}{f(z)}=\tanh(z/2)\le1,\qquad
	\frac{|f''(z)|}{f(z)}=|1-6f(z)|\le\frac58.
\end{align}

For the third derivative, \eqref{eq:Phi-short-q3}
and \eqref{eq:Phi-short-basic} give
$f'''(z)/f(z)=2\tanh(z/2)-3\tanh^3(z/2)$.
The derivative of this expression is
$2f(z)(2-9\tanh^2(z/2))$.
Since $f(z)>0$, the derivative of $f'''(z)/f(z)$
is positive when $\tanh(z/2)<\sqrt2/3$
and negative when $\tanh(z/2)>\sqrt2/3$.
Thus $f'''(z)/f(z)$ increases until
$\tanh(z/2)=\sqrt2/3$ and then decreases.

For $0\le\tanh(z/2)\le\sqrt3/2$, the maximum is
$4\sqrt2/9$, attained when $\tanh(z/2)=\sqrt2/3$.
The values when $\tanh(z/2)=0$ and
$\tanh(z/2)=\sqrt3/2$ are $0$ and $-\sqrt3/8$,
respectively.
Hence
$-\sqrt3/8\le f'''(z)/f(z)\le4\sqrt2/9$.
Since $\sqrt3/8<2/3$ and $4\sqrt2/9<2/3$,
we obtain
\begin{align}\label{eq:Phi-short-third-bound}
	\frac{|f'''(z)|}{f(z)}\le\frac23.
\end{align}

For the fourth derivative, completing the square in
\eqref{eq:Phi-short-q4} gives
$f^{(4)}(z)/f(z)=120(f(z)-1/8)^2-7/8$.
Since $0<f(z)\le1/4$, this expression lies
between $-7/8$ and $1$.
Hence
\begin{align}\label{eq:Phi-short-fourth-bound}
	\frac{|f^{(4)}(z)|}{f(z)}\le1.
\end{align}

We next estimate the first derivatives of
$f^{(j)}(z)/f(z)$ for $1\le j\le4$.
For $j=1,2$, by \eqref{eq:Phi-short-basic}
and \eqref{eq:Phi-short-q2}, we have
\begin{align}\label{eq:Phi-short-ratio-first12}
	\left|\left(\frac{f'(z)}{f(z)}\right)'\right|
	=2f(z),\qquad
	\left|\left(\frac{f''(z)}{f(z)}\right)'\right|
	=6f(z)\tanh(z/2)\le6f(z).
\end{align}

Differentiating \eqref{eq:Phi-short-q3} and using
$|72f(z)-14|\le10$, we obtain
\begin{align}\label{eq:Phi-short-ratio-first3}
	\left|\left(\frac{f'''(z)}{f(z)}\right)'\right|
	=f(z)|72f(z)-14|\le10f(z).
\end{align}

Similarly, differentiating \eqref{eq:Phi-short-q4}
and using $|1-8f(z)|\le1$, we get
\begin{align}\label{eq:Phi-short-ratio-first4}
	\left|\left(\frac{f^{(4)}(z)}{f(z)}\right)'\right|
	=30f(z)\tanh(z/2)|1-8f(z)|\le30f(z)\tanh(z/2).
\end{align}

Using
\eqref{eq:Phi-short-basic} and \eqref{eq:Phi-short-q2},
we obtain
\begin{align}\label{eq:Phi-short-ratio-second12}
	\left|\left(\frac{f'(z)}{f(z)}\right)''\right|
	&=2f(z)\tanh(z/2)\le2f(z),
	\notag\\
	\left|\left(\frac{f''(z)}{f(z)}\right)''\right|
	&=6f(z)|6f(z)-1|\le6f(z).
\end{align}

For $j=3$, differentiating $f(z)(72f(z)-14)$
and using \eqref{eq:Phi-short-basic}, we obtain
\begin{align}\label{eq:Phi-short-ratio-second3}
	\left|\left(\frac{f'''(z)}{f(z)}\right)''\right|
	=f(z)\tanh(z/2)|14-144f(z)|\le22f(z).
\end{align}

For $j=4$, differentiating \eqref{eq:Phi-short-q4}
twice gives $f(z)(-30+660f(z)-2400f(z)^2)$.
For $1/16\le f(z)\le1/4$, the quadratic has
minimum $-15$ at $f(z)=1/4$ and maximum
$123/8<16$ at $f(z)=11/80$.
We therefore obtain
\begin{align}\label{eq:Phi-short-ratio-second4}
	\left|\left(\frac{f^{(4)}(z)}{f(z)}\right)''\right|
	=f(z)|-30+660f(z)-2400f(z)^2|\le16f(z).
\end{align}

Finally, we estimate the fifth and sixth derivatives of $f(z)$.
Since $f(z)-4f(z)^2\le1/16$,
\eqref{eq:Phi-short-basic} gives
\begin{align}\label{eq:Phi-short-weighted}
	f(z)\tanh^2(z/2)\le\frac1{16},\qquad
	f(z)\tanh(z/2)\le\sqrt{\frac{f(z)}{16}}\le\frac18.
\end{align}

Differentiating
$f^{(4)}(z)=f(z)(f^{(4)}(z)/f(z))$
once and dividing by $f(z)>0$, we get
\begin{align}\label{eq:Phi-short-fifth-identity}
	\frac{f^{(5)}(z)}{f(z)}
	=
	\left(\frac{f^{(4)}(z)}{f(z)}\right)'
	-\tanh(z/2)\frac{f^{(4)}(z)}{f(z)}.
\end{align}

Differentiating the same product twice and using
$f'(z)/f(z)=-\tanh(z/2)$, we obtain
\begin{align}\label{eq:Phi-short-sixth-identity}
	\frac{f^{(6)}(z)}{f(z)}
	&=
	\left(\frac{f^{(4)}(z)}{f(z)}\right)''
	-2\tanh(z/2)
	\left(\frac{f^{(4)}(z)}{f(z)}\right)'+
	\frac{f''(z)}{f(z)}
	\frac{f^{(4)}(z)}{f(z)}.
\end{align}

Taking absolute values in \eqref{eq:Phi-short-fifth-identity},
we apply \eqref{eq:Phi-short-ratio-first4},
\eqref{eq:Phi-short-fourth-bound}, and
\eqref{eq:Phi-short-weighted} to get
\begin{align}\label{eq:Phi-short-fifth-bound}
	\frac{|f^{(5)}(z)|}{f(z)}
	\le30f(z)\tanh(z/2)+1
	\le\frac{30}{8}+1<5.
\end{align}

Taking absolute values in \eqref{eq:Phi-short-sixth-identity}
and using \eqref{eq:Phi-short-ratio-second4} and
\eqref{eq:Phi-short-ratio-first4}, we obtain
\begin{align}\label{eq:Phi-short-sixth-bound}
	\frac{|f^{(6)}(z)|}{f(z)}
	\le16f(z)+60f(z)\tanh^2(z/2)+\frac58
	\le4+\frac{60}{16}+\frac58<9,
\end{align}
where the first inequality also uses
\eqref{eq:Phi-short-low12} and
\eqref{eq:Phi-short-fourth-bound} to bound the product
by $5/8$.
The second inequality follows from $f(z)\le1/4$
and \eqref{eq:Phi-short-weighted}.
	
\noindent
\textbf{Step 2: estimates for $g(z)$ and $\ell(z)$.}
By \eqref{eq:Phi-ak}, we have
$a_k(z)=\pi(2k-1)/(2\sqrt6)+\sqrt6z^2/(4\pi)$.
For $k\in\{1,2\}$ and $0\le z\le101\pi/(50\sqrt6)$,
we use the elementary bounds
\begin{align}\label{eq:Phi-short-a}
	0<a_k(z)<\frac{10}{3},\quad
	0\le a_k'(z)<\frac{51}{50},\quad
	0<a_k''(z)<\frac25,\quad
	a_k'''(z)=0.
\end{align}
We now establish the estimates needed for $g(z)$ and $\ell(z)$.

\begin{claim}\label{clm:Phi-g-estimates}
	For $r=1,2$, we have
	\begin{align}
		|g(z)|&<\frac{13}{2}f(z),
		\label{eq:Phi-short-g-size}\\
		|g^{(r)}(z)|&<12f(z).
		\label{eq:Phi-short-g-derivatives}
	\end{align}
	Moreover,
	\begin{align}
		\frac{g(z)}{f(z)}&>-\frac1{32},
		\label{eq:Phi-short-g-lower}\\
		\left(\frac{g(z)}{f(z)}\right)''
		&<6f(z).
		\label{eq:Phi-short-g-second}
	\end{align}
\end{claim}

\begin{proof}[Proof of Claim~\ref{clm:Phi-g-estimates}]
	Since $a_k'(z)=\sqrt6z/(2\pi)$,
	\eqref{eq:Phi-g} gives
	\begin{align}\label{eq:Phi-short-g-formula}
		g(z)=-a_k(z)f''(z)-3a_k'(z)f'(z).
	\end{align}
	Differentiating once and twice,
	and using $a_k'''(z)=0$, we obtain
	\begin{align*}
		g'(z)
		&=-a_k(z)f'''(z)-4a_k'(z)f''(z)
		-3a_k''(z)f'(z),\\
		g''(z)
		&=-a_k(z)f^{(4)}(z)-5a_k'(z)f'''(z)
		-7a_k''(z)f''(z).
	\end{align*}
	By \eqref{eq:Phi-short-low12}--\eqref{eq:Phi-short-fourth-bound},
	$|f^{(j)}(z)|\le f(z)$ for $1\le j\le4$.
	Together with \eqref{eq:Phi-short-a}, this gives
	\begin{align*}
		\frac{|g(z)|}{f(z)}
		&\le a_k(z)+3a_k'(z)<\frac{13}{2},\\
		\frac{|g^{(r)}(z)|}{f(z)}
		&\le a_k(z)+5a_k'(z)+7a_k''(z)<12
		\qquad\text{for }r=1,2.
	\end{align*}
	This proves \eqref{eq:Phi-short-g-size} and
	\eqref{eq:Phi-short-g-derivatives}.
	
	We next prove \eqref{eq:Phi-short-g-lower}.
	By \eqref{eq:Phi-short-basic} and
	\eqref{eq:Phi-short-q2}, we have
	\begin{align}\label{eq:Phi-short-g-ratio}
		\frac{g(z)}{f(z)}
		=-a_k(z)\bigl(1-6f(z)\bigr)
		+3a_k'(z)\tanh(z/2).
	\end{align}
	If $\tanh^2(z/2)\le1/3$, then
	$1-6f(z)=(3\tanh^2(z/2)-1)/2\le0$,
	so $g(z)/f(z)\ge0$.
	
	Suppose that $\tanh^2(z/2)>1/3$.
	We use $(2k-1)\pi^2/3<10$ and
	$z^2\le(31/10)z\tanh(z/2)$ from
	\eqref{eq:Phi-short-zT}.
	The resulting coefficient of $z\tanh(z/2)$ is
	$12-(31/10)(3\tanh^2(z/2)-1)>0$,
	since $\tanh^2(z/2)<3/4$.
	Hence we may also use
	$z\tanh(z/2)\ge2\tanh^2(z/2)$ to obtain
	\begin{align*}
		4\pi\sqrt{\frac23}\frac{g(z)}{f(z)}
		&=12z\tanh(z/2)
		-\left(\frac{(2k-1)\pi^2}{3}+z^2\right)
		\bigl(3\tanh^2(z/2)-1\bigr)\\
		&\ge10+\frac15\tanh^2(z/2)
		-\frac{93}{5}\tanh^4(z/2)
		>-\frac5{16}.
	\end{align*}
	The last polynomial decreases as $\tanh^2(z/2)$
	increases from $1/3$ to $3/4$,
	and equals $-5/16$ at $3/4$.
	Dividing by $4\pi\sqrt{2/3}$ and using
	$\pi\sqrt{2/3}>5/2$ proves
	\eqref{eq:Phi-short-g-lower}.
	
	It remains to prove \eqref{eq:Phi-short-g-second}.
	Differentiating \eqref{eq:Phi-short-g-ratio} twice
	and using \eqref{eq:Phi-short-basic}, we get
	\begin{align}\label{eq:Phi-short-g-curvature-identity}
		\frac{\pi\sqrt{2/3}}{f(z)}
		\left(\frac{g(z)}{f(z)}\right)''
		=-\frac1{f(z)}+18-18z\tanh(z/2)+
		\bigl((2k-1)\pi^2+3z^2\bigr)
		\bigl(1-6f(z)\bigr).
	\end{align}
	By \eqref{eq:Phi-short-zT},
	$3z^2\le(93/10)z\tanh(z/2)$ and
	$z\tanh(z/2)\ge2-8f(z)$.
	Since $1-6f(z)\le1$, substitution into
	\eqref{eq:Phi-short-g-curvature-identity} gives
	\begin{align}\label{eq:Phi-short-g-curvature-upper}
		\frac{\pi\sqrt{2/3}}{f(z)}
		\left(\frac{g(z)}{f(z)}\right)''
		\le\frac35+\frac{348}{5}f(z)-\frac1{f(z)}
		+(2k-1)\pi^2\bigl(1-6f(z)\bigr).
	\end{align}
	If $f(z)\le1/6$, the right-hand side is at most
	\begin{align*}
		\frac{153}{5}-\frac{552}{5}f(z)-\frac1{f(z)}
		\le\frac{153}{5}-2\sqrt{\frac{552}{5}}<10,
	\end{align*}
	where we use $(2k-1)\pi^2<30$ and
	the arithmetic--geometric mean inequality.
	
	If $f(z)\ge1/6$, the last term in
	\eqref{eq:Phi-short-g-curvature-upper}
	is nonpositive, and the remaining terms satisfy
	\begin{align*}
		\frac35+\frac{348}{5}f(z)-\frac1{f(z)}
		\le14,
	\end{align*}
	because $f(z)\le1/4$.
	Thus \eqref{eq:Phi-short-g-second} follows from
	$14/(\pi\sqrt{2/3})<6$.
\end{proof}

\begin{claim}\label{clm:Phi-ell-estimates}
	For $r=1,2$, we have
	\begin{align}
		|\ell(z)|&<23f(z),
		\label{eq:Phi-short-ell-size}\\
		|\ell^{(r)}(z)|&<210f(z).
		\label{eq:Phi-short-ell-derivatives}
	\end{align}
	In addition,
	\begin{align}\label{eq:Phi-short-ell-second}
		\left|
		\left(\frac{\ell(z)}{f(z)}\right)''
		\right|<1050f(z).
	\end{align}
\end{claim}

\begin{proof}[Proof of Claim~\ref{clm:Phi-ell-estimates}]
	Substituting $a_k'(z)=\sqrt6z/(2\pi)$ and
	$a_k''(z)=\sqrt6/(2\pi)$ into
	\eqref{eq:Phi-ell}, we obtain
	\begin{align}\label{eq:Phi-short-ell-ratio}
		\frac{\ell(z)}{f(z)}
		&=\left(3a_k(z)a_k''(z)+6a_k'(z)^2
		+\frac{\pi^2}{144}\right)
		\frac{f''(z)}{f(z)}+
		3a_k'(z)a_k''(z)\frac{f'(z)}{f(z)}
		\notag\\
		&\quad+\frac{a_k(z)^2}{2}\frac{f^{(4)}(z)}{f(z)}
		+4a_k(z)a_k'(z)\frac{f'''(z)}{f(z)}.
	\end{align}
	Using
	\eqref{eq:Phi-short-low12}, \eqref{eq:Phi-short-third-bound} and \eqref{eq:Phi-short-fourth-bound},
	we have
	\begin{align*}
		\frac{|\ell(z)|}{f(z)}
		&\le\frac{a_k(z)^2}{2}
		+\frac83a_k(z)a_k'(z)+
		\frac58\left(3a_k(z)a_k''(z)+6a_k'(z)^2
		+\frac{\pi^2}{144}\right)
		+3a_k'(z)a_k''(z)<23,
	\end{align*}
	where the last inequality uses
	\eqref{eq:Phi-short-a} and $\pi^2<10$.
	This proves \eqref{eq:Phi-short-ell-size}.
	
	We next estimate $\ell'(z)$ and $\ell''(z)$.
	Multiplying \eqref{eq:Phi-short-ell-ratio} by $f(z)$
	and differentiating, we get
	\begin{align*}
		\ell'(z)
		&=\frac{a_k(z)^2}{2}f^{(5)}(z)
		+5a_k(z)a_k'(z)f^{(4)}(z)+
		\left(7a_k(z)a_k''(z)+10a_k'(z)^2
		+\frac{\pi^2}{144}\right)f'''(z)\\
		&\qquad+
		18a_k'(z)a_k''(z)f''(z)
		+3a_k''(z)^2f'(z).
	\end{align*}
	Differentiating once more and using $a_k'''(z)=0$,
	we obtain
	\begin{align*}
		\ell''(z)
		&=\frac{a_k(z)^2}{2}f^{(6)}(z)
		+6a_k(z)a_k'(z)f^{(5)}(z)+
		\left(12a_k(z)a_k''(z)+15a_k'(z)^2
		+\frac{\pi^2}{144}\right)f^{(4)}(z)\\
		&\qquad+
		45a_k'(z)a_k''(z)f'''(z)
		+21a_k''(z)^2f''(z).
	\end{align*}
	By Step 1, $|f^{(j)}(z)|\le f(z)$ for $1\le j\le4$.
	Applying \eqref{eq:Phi-short-fifth-bound} and
	\eqref{eq:Phi-short-sixth-bound}
	to the remaining terms gives, for $r=1,2$,
	\begin{align*}
		\frac{|\ell^{(r)}(z)|}{f(z)}
		&\le\frac92a_k(z)^2+30a_k(z)a_k'(z)
		+12a_k(z)a_k''(z)+15a_k'(z)^2\\
		&\quad+
		\frac{\pi^2}{144}
		+45a_k'(z)a_k''(z)+21a_k''(z)^2<210.
	\end{align*}
	For $r=1$, the estimate is no larger term by term
	than the bound for $r=2$.
	The numerical bound follows from
	\eqref{eq:Phi-short-a} and $\pi^2<10$,
	proving \eqref{eq:Phi-short-ell-derivatives}.
	
We finally prove \eqref{eq:Phi-short-ell-second}.
We differentiate each term in \eqref{eq:Phi-short-ell-ratio}
twice by the product rule.
For the first derivatives of $f^{(j)}(z)/f(z)$, we use
\eqref{eq:Phi-short-ratio-first12},
\eqref{eq:Phi-short-ratio-first3}, and
\eqref{eq:Phi-short-ratio-first4}.
For the second derivatives, we use
\eqref{eq:Phi-short-ratio-second12},
\eqref{eq:Phi-short-ratio-second3}, and
\eqref{eq:Phi-short-ratio-second4}.
Using also $0\le\tanh(z/2)\le1$, we obtain
	\begin{align*}
		\left|
		\left(\frac{\ell(z)}{f(z)}\right)''
		\right|
		&\le a_k'(z)^2+a_k(z)a_k''(z)
		+8a_k'(z)a_k''(z)+\frac{75}{8}a_k''(z)^2\\
		&\quad+
		f(z)\left(
		8a_k(z)^2+148a_k(z)a_k'(z)+116a_k'(z)^2
		\right.\\
		&\quad\left.
		+98a_k(z)a_k''(z)+186a_k'(z)a_k''(z)
		+12a_k''(z)^2+\frac{\pi^2}{24}
		\right).
	\end{align*}
The terms without an explicit factor $f(z)$ are bounded
using \eqref{eq:Phi-short-low12},
\eqref{eq:Phi-short-third-bound}, and
\eqref{eq:Phi-short-fourth-bound}.
For example, the coefficient $148=2\cdot30+4\cdot22$
comes from the last two terms in
\eqref{eq:Phi-short-ell-ratio}.
The remaining coefficients follow by the same product rule.

Substituting \eqref{eq:Phi-short-a} and using $\pi^2<10$,
we obtain
\begin{align*}
	\left|
	\left(\frac{\ell(z)}{f(z)}\right)''
	\right|
	&<8+925f(z)\\
	&<1049f(z)<1050f(z),
\end{align*}
where the second inequality uses
$f(z)>2/31$ from \eqref{eq:f-range},
which gives $8<124f(z)$.
This proves \eqref{eq:Phi-short-ell-second}.
\end{proof}
	
\noindent
\textbf{Step 3: estimates for $\Phi(z)$.}
We first prove the lower bound for $\Phi(z)$.
Substituting \eqref{eq:Phi-short-g-lower} and
\eqref{eq:Phi-short-ell-size} into \eqref{eq:Phi}, we obtain
\begin{align}\label{eq:Phi-short-positive}
	\frac{\Phi(z)}{f(z)}
	&=1+\varepsilon\frac{g(z)}{f(z)}
	+\varepsilon^2\frac{\ell(z)}{f(z)}\ge1-\frac{\varepsilon}{32}-23\varepsilon^2
	\ge1-\frac1{32\cdot66}-\frac{23}{66^2}
	>\frac{993}{1000},
\end{align}
where we use $0\le\varepsilon\le1/66$.
In particular, $\Phi(z)>0$.

Next, applying the triangle inequality to \eqref{eq:Phi},
we obtain
\begin{align*}
	\frac{|\Phi(z)|}{f(z)}
	&\le1+\varepsilon\frac{|g(z)|}{f(z)}
	+\varepsilon^2\frac{|\ell(z)|}{f(z)}
	\le1+\frac{13}{2\cdot66}+\frac{23}{66^2}
	<\frac{111}{100},
\end{align*}
where the second inequality follows from
\eqref{eq:Phi-short-g-size},
\eqref{eq:Phi-short-ell-size},
and $0\le\varepsilon\le1/66$.

For $r=1,2$, \eqref{eq:Phi-short-low12} gives
$|f^{(r)}(z)|\le f(z)$.
Applying \eqref{eq:Phi-short-g-derivatives} and
\eqref{eq:Phi-short-ell-derivatives}
to the derivatives of \eqref{eq:Phi}, we obtain
\begin{align*}
	\frac{|\Phi^{(r)}(z)|}{f(z)}
	\le
	\frac{|f^{(r)}(z)|}{f(z)}
	+\varepsilon\frac{|g^{(r)}(z)|}{f(z)}
	+\varepsilon^2\frac{|\ell^{(r)}(z)|}{f(z)}
	\le1+\frac{12}{66}+\frac{210}{66^2}
	<\frac{31}{25}.
\end{align*}

Finally, \eqref{eq:Phi-short-basic} gives
$-(\log f)''(z)=2f(z)$.
Since $\Phi(z)>0$, we may differentiate
$\log\Phi(z)=\log f(z)+\log(\Phi(z)/f(z))$
twice to get
\begin{align*}
	-(\log\Phi)''(z)
	=2f(z)-
	\frac{
		\varepsilon\left(g(z)/f(z)\right)''
		+\varepsilon^2\left(\ell(z)/f(z)\right)''
	}{
		\Phi(z)/f(z)
	}
	+
	\left(
	\frac{\mathrm d}{\mathrm dz}
	\log\frac{\Phi(z)}{f(z)}
	\right)^2.
\end{align*}
Using \eqref{eq:Phi-short-g-second} and
\eqref{eq:Phi-short-ell-second}, we therefore have
\begin{align*}
	-(\log\Phi)''(z)
	\ge f(z)\left(
	2-\frac{6\varepsilon+1050\varepsilon^2}
	{\Phi(z)/f(z)}
	\right)
	\ge f(z)\left(
	2-\frac{6/66+1050/66^2}{993/1000}
	\right)
	>\frac{83}{50}f(z),
\end{align*}
where the second inequality uses
\eqref{eq:Phi-short-positive}
and $0\le\varepsilon\le1/66$.
This proves \eqref{eq:Phibounds}.
\end{proof}

Recall from \eqref{eq:I-normalized} that
\begin{align*}
	I_r(m,n)=
	\frac{\partial_m^r\A_k(m,n)}
	{\left(\dfrac{\pi}{\sqrt{6(n-\tfrac1{24})}}\right)^rP'(n)},
	\qquad r=0,1,2.
\end{align*}
Combining \eqref{eq:errors} with the bounds in
\eqref{eq:Phibounds}, we obtain the following result.

\begin{proposition}\label{prop:determinant}
	For $k\in\{1,2\}$, let $n>\tfrac1{24}$ and $m\ge0$
	be real numbers satisfying
	\begin{align*}
		0<\varepsilon=\frac1{\sqrt{n-\tfrac1{24}}}
		\le\frac1{66},
		\qquad
		0\le\frac{m}{\sqrt{n-\tfrac1{24}}}
		\le\frac{101}{50}.
	\end{align*}
	Then
	\begin{align*}
		\frac1{20}<I_0(m,n)<\frac25,\qquad
		I_1(m,n)^2-I_0(m,n)I_2(m,n)>\frac1{100000}.
	\end{align*}
\end{proposition}

\begin{proof}
	Recall from \eqref{eq:expansion-z} that
	$z=\pi m/\sqrt{6(n-\tfrac1{24})}$.
	
	We first bound $I_0(m,n)$.
	By \eqref{eq:errors}, we have
	$I_0(m,n)>\Phi(z)-7/20000$.
	Using \eqref{eq:Phibounds} and \eqref{eq:f-range},
	we obtain
	\begin{align*}
		I_0(m,n)
		>\frac{993}{1000}f(z)-\frac7{20000}
		>\frac{993}{1000}\frac2{31}-\frac7{20000}
		>\frac1{20}.
	\end{align*}
	Similarly, \eqref{eq:errors} and \eqref{eq:Phibounds} give
	\begin{align*}
		I_0(m,n)
		<\frac{111}{100}f(z)+\frac7{20000}
		\le\frac{111}{400}+\frac7{20000}<\frac25,
	\end{align*}
	where we use $f(z)\le1/4$ from \eqref{eq:f-range}.
	
	We next estimate $I_1(m,n)^2-I_0(m,n)I_2(m,n)$.
	Since $\Phi(z)>0$, \eqref{eq:Phibounds} implies
	\begin{align}\label{eq:det-Phi-lower}
		\Phi'(z)^2-\Phi(z)\Phi''(z)
		=\Phi(z)^2\bigl(-(\log\Phi)''(z)\bigr)>\frac{83}{50}
		\left(\frac{993}{1000}\right)^2f(z)^3
		>\frac{163}{100}f(z)^3.
	\end{align}
	
	We compare $I_1(m,n)^2-I_0(m,n)I_2(m,n)$ with
	$\Phi'(z)^2-\Phi(z)\Phi''(z)$ in
	\eqref{eq:det-Phi-lower}.
	For $j=0,1,2$, define
	\begin{align}\label{eq:det-omega}
		\omega_j(m,n)=I_j(m,n)-\Phi^{(j)}(z),
	\end{align}
	where $\Phi^{(0)}(z)=\Phi(z)$.
	By \eqref{eq:errors}, we have
	\begin{align}\label{eq:det-omega-bounds}
		|\omega_0(m,n)|<\frac7{20000},\qquad
		|\omega_1(m,n)|<\frac1{1020},\qquad
		|\omega_2(m,n)|<\frac1{305}.
	\end{align}
	Substituting
	$I_j(m,n)=\Phi^{(j)}(z)+\omega_j(m,n)$, we get
	\begin{align}\label{eq:det-expansion}
		&I_1(m,n)^2-I_0(m,n)I_2(m,n)
		-\bigl(\Phi'(z)^2-\Phi(z)\Phi''(z)\bigr)
		\\
		&=
		2\Phi'(z)\omega_1(m,n)
		-\Phi(z)\omega_2(m,n)
		-\Phi''(z)\omega_0(m,n)+
		\omega_1(m,n)^2-\omega_0(m,n)\omega_2(m,n).
	\end{align}
	Applying the triangle inequality to
	\eqref{eq:det-expansion}, we have
	\begin{align}\label{eq:det-difference-bound}
		&\left|
		I_1(m,n)^2-I_0(m,n)I_2(m,n)
		-\bigl(\Phi'(z)^2-\Phi(z)\Phi''(z)\bigr)
		\right|
		\notag\\
		&\quad\le
		2|\Phi'(z)|\,|\omega_1(m,n)|
		+|\Phi(z)|\,|\omega_2(m,n)|
		+|\Phi''(z)|\,|\omega_0(m,n)|
		\notag\\
		&\qquad+
		|\omega_1(m,n)|^2
		+|\omega_0(m,n)|\,|\omega_2(m,n)|
		\notag\\
		&\quad<
		f(z)\left(
		\frac{62}{25\cdot1020}
		+\frac{111}{100\cdot305}
		+\frac{31\cdot7}{25\cdot20000}
		\right)
		+\frac1{1020^2}+\frac7{20000\cdot305},
	\end{align}
	where the last inequality uses
	\eqref{eq:Phibounds} and \eqref{eq:det-omega-bounds}.
	
Combining \eqref{eq:det-Phi-lower} and
\eqref{eq:det-difference-bound},
and factoring out $f(z)>0$, we obtain
\begin{align}\label{eq:det-lower}
	I_1(m,n)^2-I_0(m,n)I_2(m,n)&>
	f(z)\left(
	\frac{163}{100}f(z)^2
	-\frac{62}{25\cdot1020}
	-\frac{111}{100\cdot305}
	\right.
	\notag\\
	&\quad\left.
	-\frac{31\cdot7}{25\cdot20000}
	-\frac1{f(z)}
	\left(\frac1{1020^2}+\frac7{20000\cdot305}\right)
	\right)
	\notag\\
	&>
	\frac{31}{200000}f(z)
	>\frac{31}{200000}\frac2{31}
	=\frac1{100000},
\end{align}
where the second inequality uses
$f(z)^2>(2/31)^2$ and $1/f(z)<31/2$
from \eqref{eq:f-range},
together with a direct comparison of the rational constants.
The last inequality uses $f(z)>2/31$.
This completes the proof.
\end{proof}

\subsection{Proof of Theorem~\ref{thm:strict} for $n\ge5000$ and $0\le m\le2\sqrt n$}

\begin{theorem}\label{prop:central}
	Fix $k\in\{1,2\}$ and an integer $n\ge5000$.
	For every integer $m$ with $0\le m\le2\sqrt n$,
	we have $A_k(m+j,n)>0$ for $j=-1,0,1$, and
	\begin{align}\label{eq:centralmargin}
		\Lambda_k(m,n)
		=\log\frac{A_k(m,n)^2}
		{A_k(m-1,n)A_k(m+1,n)}
		>\frac1{20000n}.
	\end{align}
\end{theorem}

\begin{proof}
Recall from \eqref{eq:model} that, for real $y$,
\begin{align*}
	\A_k(y,n)
	={}&\sum_{j=1}^{\infty}(-1)^{j-1}
	\Bigl\{
	B\!\left(
	n-\left(k-\tfrac12\right)j^2-yj+\tfrac j2
	\right)-
	B\!\left(
	n-\left(k-\tfrac12\right)j^2-yj-\tfrac j2
	\right)
	\Bigr\}.
\end{align*}
The proof consists of three steps.
We first bound $\A_k(y,n)$ and
$-\partial_y^2\log\A_k(y,n)$.
We then compare $\A_k(m+j,n)$ with $A_k(m+j,n)$
for $j=-1,0,1$.
Finally, we apply \eqref{eq:kernel} and
\eqref{eq:transfer} to prove \eqref{eq:centralmargin}.
	
\noindent
\textbf{Step 1: estimates for $\A_k(y,n)$.}
Since $n\ge5000$, we have
$0<\varepsilon=1/\sqrt{n-\tfrac1{24}}<1/66$
and $\sqrt n>70$.
For every real $y$ with $|y|\le2\sqrt n+1$, we obtain
\begin{align}\label{eq:central-enlarged-domain}
	\frac{|y|}{\sqrt{n-\tfrac1{24}}}
	\le
	\frac{2+1/\sqrt n}{\sqrt{1-1/(24n)}}
	\le
	\frac{2+1/70}{\sqrt{1-1/(24\cdot5000)}}
	<\frac{101}{50}.
\end{align}
Hence the assumptions of Proposition~\ref{prop:determinant}
hold with $m=|y|$.
	
	For $0\le y\le2\sqrt n+1$,
	\eqref{eq:I-normalized} gives
	$\A_k(y,n)=P'(n)I_0(y,n)$.
	Since $P'(n)>0$ by \eqref{eq:Pprime-form},
	Proposition~\ref{prop:determinant} implies
	\begin{align}\label{eq:central-smooth-positive}
		\A_k(y,n)>\frac{P'(n)}{20}>0.
	\end{align}
	Using \eqref{eq:I-normalized} for the first two
	derivatives, we also have
	\begin{align}\label{eq:central-smooth-curvature}
		-\frac{\partial^2}{\partial y^2}\log\A_k(y,n)
		&=
		\left(
		\frac{\partial_y\A_k(y,n)}{\A_k(y,n)}
		\right)^2
		-\frac{\partial_y^2\A_k(y,n)}{\A_k(y,n)}
		\notag\\
		&=
		\frac{\pi^2}{6(n-\tfrac1{24})}
		\frac{I_1(y,n)^2-I_0(y,n)I_2(y,n)}
		{I_0(y,n)^2}
		\notag\\
		&>
		\frac{\pi^2}{6(n-\tfrac1{24})}
		\frac{1/100000}{(2/5)^2}\notag\\
		&>\frac1{10000n}.
	\end{align}
	The first inequality follows from
	Proposition~\ref{prop:determinant}.
	
	By the evenness in \eqref{eq:A-even},
	\eqref{eq:central-smooth-positive} and
	\eqref{eq:central-smooth-curvature} also hold
	for negative $y$ satisfying $|y|\le2\sqrt n+1$.
	In particular, both estimates hold throughout
	$[m-1,m+1]$.
	
\noindent
\textbf{Step 2: comparison with $A_k(m+j,n)$.}
	Since $n-\tfrac1{24}\ge120$, \eqref{eq:Pprime-lower}
	gives
	\begin{align}\label{eq:central-Pprime-bound}
		P'(n)>
		\frac1{20n^{3/2}}
		\exp\!\left(
		\pi\sqrt{\frac23(n-\tfrac1{24})}
		\right).
	\end{align}
	
	For $j=-1,0,1$, we apply \eqref{eq:modelerror}
	at $m+j$.
	If $m+j<0$, the same estimate follows from
	the symmetry $A_k(-m,n)=A_k(m,n)$ and
	\eqref{eq:A-even}.
	Dividing \eqref{eq:modelerror} by the lower bound
	in \eqref{eq:central-smooth-positive}, we obtain
	\begin{align}\label{eq:central-relative-preliminary}
		\left|
		\frac{A_k(m+j,n)}{\A_k(m+j,n)}-1
		\right|
		&=
		\frac{|A_k(m+j,n)-\A_k(m+j,n)|}
		{\A_k(m+j,n)}
		\notag\\
		&<
		\frac8{\sqrt n\,P'(n)}
		\exp\!\left(
		\pi\sqrt{\frac{n-\tfrac1{24}}6}
		\right)
		\notag\\
		&<
		160n\exp\!\left(
		-\pi\sqrt{\frac{n-\tfrac1{24}}6}
		\right)\\
		&<160ne^{-\sqrt n},
	\end{align}
	where the second inequality uses
	\eqref{eq:central-Pprime-bound}.
	The last inequality follows from
	$\pi\sqrt{(n-\tfrac1{24})/6}>\sqrt n$
	for $n\ge5000$.
	
We next bound $160ne^{-\sqrt n}$ for $n\ge5000$.
For $x>0$, $\beta\ge0$, and $\alpha>0$, we have
\begin{align}\label{eq:monotone}
	\frac{\mathrm d}{\mathrm dx}
	\log\bigl(x^\beta e^{-\alpha\sqrt x}\bigr)
	=
	\frac{2\beta-\alpha\sqrt x}{2x}<0
	\qquad\text{if }\sqrt x>\frac{2\beta}{\alpha}.
\end{align}
Taking $\beta=2$ and $\alpha=1$ shows that
$x^2e^{-\sqrt x}$ decreases for $x\ge5000$.
Since $e^{\sqrt{5000}}>160\cdot10^{10}\cdot5000^2$,
we obtain
\begin{align*}
	160n^2e^{-\sqrt n}
	\le160\cdot5000^2e^{-\sqrt{5000}}
	<\frac1{10^{10}}.
\end{align*}
Dividing by $n$ gives
$160ne^{-\sqrt n}<1/(10^{10}n)$.
	By
	\eqref{eq:central-relative-preliminary}, we get
	\begin{align}\label{eq:central-relative-bound}
		\left|
		\frac{A_k(m+j,n)}{\A_k(m+j,n)}-1
		\right|
		<\frac1{10^{10}n},
		\qquad j=-1,0,1.
	\end{align}
	Since $\A_k(m+j,n)>0$, it follows that
	\begin{align*}
		A_k(m+j,n)
		>
		\left(1-\frac1{10^{10}n}\right)
		\A_k(m+j,n)>0.
	\end{align*}
	
\noindent
\textbf{Step 3: the bound for $\Lambda_k(m,n)$.}
	Apply \eqref{eq:kernel} to
	$G(y)=\A_k(y,n)$ on $[m-1,m+1]$.
	By \eqref{eq:central-smooth-curvature}, we obtain
	\begin{align}\label{eq:central-smooth-discrete}
		\log\frac{\A_k(m,n)^2}
		{\A_k(m-1,n)\A_k(m+1,n)}
		&=
		-\int_{-1}^{1}(1-|t|)
		\left.
		\frac{\partial^2}{\partial y^2}
		\log\A_k(y,n)
		\right|_{y=m+t}\dd t
		\notag\\
		&>
		\frac1{10000n}
		\int_{-1}^{1}(1-|t|)\dd t\notag\\
		&=\frac1{10000n}.
	\end{align}
	
	Finally, apply \eqref{eq:transfer} with
	$x_j=\A_k(m+j,n)$,
	$b_j=A_k(m+j,n)$, and
	$\eta=1/(10^{10}n)$.
	By \eqref{eq:central-relative-bound},
	the relative bound required in Lemma~\ref{lem:transfer}
	holds with $\eta=1/(10^{10}n)$.
	Together with \eqref{eq:central-smooth-discrete},
	we get
	\begin{align*}
		\Lambda_k(m,n)
		&\ge
		\log\frac{\A_k(m,n)^2}
		{\A_k(m-1,n)\A_k(m+1,n)}
		-\frac{4}{10^{10}n\left(1-\frac1{10^{10}n}\right)}
		\\
		&>
		\frac1{10000n}-\frac8{10^{10}n}\\
		&>\frac1{20000n},
	\end{align*}
	where the second inequality uses
	$1/(10^{10}n)<1/2$.
	This completes the proof.
\end{proof}

\section{The case $n\ge5000$ and $2\sqrt n\le m\le n/2+1$}
\label{sec:outer}

We retain a finite part of \eqref{eq:coeff} and replace $p(s)$ by
$P(s)$. We prove positivity and strict log-concavity of the resulting
sum, control its difference from $A_k(m,n)$, and apply
\eqref{eq:kernel} and \eqref{eq:transfer} to prove
Theorem~\ref{thm:strict} for $n\ge5000$ and
$2\sqrt n\le m\le n/2+1$.

\subsection{Bounds for partition sums}

\begin{lemma}\label{lem:partition-tail}
	For $k\in\{1,2\}$ and nonnegative integers $n,y$, define
	$s_j=n-(k-\tfrac12)j^2-yj+j/2$ for $j\ge1$.
	For every positive integer $\ell$, we have
	\begin{align}\label{eq:partition-tail}
		\left|
		\sum_{j=\ell}^{\infty}(-1)^{j-1}
		\bigl(p(s_j)-p(s_j-j)\bigr)
		\right|
		\le
		\sum_{j=\ell}^{\infty}
		\bigl(p(s_j)-p(s_j-j)\bigr)
		\le p(s_\ell).
	\end{align}
\end{lemma}

\begin{proof}
	Since $s_j-s_{j+1}=(2k-1)j+k-1+y\ge j$,
	monotonicity of $p(s)$ gives
	\begin{align}\label{eq:partition-tail-termwise}
		0\le p(s_j)-p(s_j-j)
		\le p(s_j)-p(s_{j+1}).
	\end{align}
	There is an integer $L\ge\ell$ such that $s_j<0$
	for every $j>L$.
	Thus all terms with $j>L$ vanish, and summing
	\eqref{eq:partition-tail-termwise} gives
	\begin{align*}
		\sum_{j=\ell}^{\infty}
		\bigl(p(s_j)-p(s_j-j)\bigr)
		\le
		\sum_{j=\ell}^{L}
		\bigl(p(s_j)-p(s_{j+1})\bigr)
		=p(s_\ell)-p(s_{L+1})
		=p(s_\ell).
	\end{align*}
	The triangle inequality and the nonnegativity in
	\eqref{eq:partition-tail-termwise}
	prove the remaining inequality.
\end{proof}

For $n\ge5000$ and every integer $t\ge0$, we also use
\begin{align}\label{eq:envelopes}
	p(t)\le
	\frac6{5n}
	\exp\!\left(\pi\sqrt{\frac23\max(t,n/8)}\right).
\end{align}
Indeed, direct differentiation gives
\begin{align*}
	\frac{\mathrm d}{\mathrm du}
	\log\frac{\exp(\pi\sqrt{2u/3})}{u}
	=
	\frac{\pi\sqrt{2u/3}-2}{2u}>0
	\qquad\text{for }u\ge23/24.
\end{align*}
For $t\ge1$, this monotonicity and \eqref{eq:p-upper} imply
\begin{align*}
	p(t)
	&\le
	\frac{
		\exp\!\left(
		\pi\sqrt{\frac23(\max(t,n/8)-\tfrac1{24})}
		\right)
	}{
		4\sqrt3(\max(t,n/8)-\tfrac1{24})
	}\\
	&\le
	\frac{
		\exp\!\left(\pi\sqrt{\frac23\max(t,n/8)}\right)
	}{
		4\sqrt3(n/8-\tfrac1{24})
	}\\
	&<
	\frac6{5n}
	\exp\!\left(\pi\sqrt{\frac23\max(t,n/8)}\right),
\end{align*}
where the last inequality follows from $\sqrt3>17/10$
and $n\ge5000$.
For $t=0$, the bound follows from $p(0)=p(1)=1$ and $n/8>1$.

\subsection{Estimates for the first difference}

\begin{lemma}\label{lem:P-derivatives}
	Let $P(s)$ be defined by \eqref{eq:PQ}.
	For every real $s$ satisfying $s-\tfrac1{24}\ge1250$,
	we have
	\begin{align}\label{eq:P-ratio}
		\frac{\pi}{\sqrt{6(s-\tfrac1{24})}}
		-\frac3{2(s-\tfrac1{24})}
		\le
		\frac{P''(s)}{P'(s)}
		<
		\frac{257}{200\sqrt{s-\tfrac1{24}}}
	\end{align}
	and
	\begin{align}\label{eq:P-ratio-derivative}
		-\frac{257}{400(s-\tfrac1{24})^{3/2}}
		<
		\frac{P'''(s)}{P'(s)}
		-\frac{P''(s)^2}{P'(s)^2}
		<
		-\frac1{2(s-\tfrac1{24})^{3/2}}.
	\end{align}
\end{lemma}

\begin{proof}
	We first prove \eqref{eq:P-ratio}.
	For $s-\tfrac1{24}\ge1250$, write
	$v=\sqrt6/(2\pi\sqrt{s-\tfrac1{24}})$.
	Then $0<v<1/90$, and \eqref{eq:Pderiv} gives
	\begin{align}\label{eq:P-ratio-form}
		\frac{P''(s)}{P'(s)}
		=
		\frac{\pi}{\sqrt{6(s-\tfrac1{24})}}
		\frac{H_1(v)}{H_0(v)}.
	\end{align}
	For $0<v<1/90$, we have $H_0(v)>0$ and $H_0'(v)<0$.
	Since
	\begin{align*}
		H_1(v)
		&=(1-3v)H_0(v)-v^2H_0'(v),\\
		H_1(v)-H_0(v)
		&=-3v(1-4v+5v^2)<0,
	\end{align*}
	it follows that
	$1-3v<H_1(v)/H_0(v)<1$.
	Substituting these bounds into \eqref{eq:P-ratio-form},
	we obtain
	\begin{align*}
		\frac{P''(s)}{P'(s)}
	&\ge
		\frac{\pi}{\sqrt{6(s-\tfrac1{24})}}(1-3v)
		=
		\frac{\pi}{\sqrt{6(s-\tfrac1{24})}}
		-\frac3{2(s-\tfrac1{24})},\\
		\frac{P''(s)}{P'(s)}
		&<
		\frac{\pi}{\sqrt{6(s-\tfrac1{24})}}
		<
		\frac{257}{200\sqrt{s-\tfrac1{24}}}.
	\end{align*}
	This proves \eqref{eq:P-ratio}.
	
	We now prove \eqref{eq:P-ratio-derivative}.
	By  \eqref{eq:Pderiv}, we have
	\begin{align}\label{eq:P-ratio-identities}
		\frac{\mathrm d}{\mathrm ds}
		\frac{P''(s)}{P'(s)}
		=
		\frac{P'''(s)}{P'(s)}
		-\left(\frac{P''(s)}{P'(s)}\right)^2=
		\frac{\pi^2}{6(s-\tfrac1{24})}
		\frac{H_2(v)H_0(v)-H_1(v)^2}{H_0(v)^2}.
	\end{align}
	A direct calculation gives
	\begin{align}\label{eq:P-ratio-polynomial}
		H_2(v)H_0(v)-H_1(v)^2
		=
		-v\bigl(
		1-12v+60v^2-150v^3+180v^4-90v^5
		\bigr).
	\end{align}
For $0<v<1/90$, we have
\begin{align}\label{eq:P-ratio-polynomial-lower}
	1-12v+60v^2-150v^3+180v^4-90v^5
	\ge
	1-\frac{12}{90}-\frac{150}{90^3}
	-\frac{90}{90^5}
	>\frac45.
\end{align}
Since $0<H_0(v)^2<1$, dividing the left side of
\eqref{eq:P-ratio-polynomial-lower} by $H_0(v)^2$
gives a value greater than $4/5$.

For the upper bound, subtracting $H_0(v)^2$ from
the left side of \eqref{eq:P-ratio-polynomial-lower},
we obtain
	\begin{align*}
		v(-6+45v-132v^2+171v^3-90v^4)
		<
		v\left(
		-6+\frac{45}{90}+\frac{171}{90^3}
		\right)<0.
	\end{align*}
	Therefore,
	\begin{align}\label{eq:P-ratio-polynomial-bound}
		\frac45
		<
		\frac{
			1-12v+60v^2-150v^3+180v^4-90v^5
		}{H_0(v)^2}
		<1.
	\end{align}
	Combining \eqref{eq:P-ratio-identities},
	\eqref{eq:P-ratio-polynomial}, and
	\eqref{eq:P-ratio-polynomial-bound}, we get
	\begin{align*}
		-\frac{\pi}{2\sqrt6(s-\tfrac1{24})^{3/2}}
		<
		\frac{\mathrm d}{\mathrm ds}
		\frac{P''(s)}{P'(s)}
		<
		-\frac45
		\frac{\pi}{2\sqrt6(s-\tfrac1{24})^{3/2}}.
	\end{align*}
	Finally, $5/8<\pi/(2\sqrt6)<257/400$
	gives \eqref{eq:P-ratio-derivative}.
	This completes the proof.
\end{proof}

\begin{lemma}\label{lem:P-difference}
	Define $\Delta(d)=P(d)-P(d-1)$ for real $d>25/24$.
	For $d\ge1252$, we have
	\begin{align}\label{eq:Delta-base}
		\Delta(d)
		>
		\frac4{25d^{3/2}}
		\exp\!\left(
		\pi\sqrt{\frac23(d-\tfrac{25}{24})}
		\right),
\qquad
		(\log\Delta)''(d)
		<-\frac1{3d^{3/2}}.
	\end{align}
\end{lemma}

\begin{proof}
	By the fundamental theorem of calculus, we have 
	\begin{align}\label{eq:Delta-integral}
		\Delta(d)=\int_0^1P'(d-u)\dd u.
	\end{align}
For $0\le u\le1$, we have
\begin{align*}
	d-\tfrac{25}{24}\le d-u-\tfrac1{24}<d,
	\qquad
	0<\frac{\sqrt6}{2\pi\sqrt{d-u-\tfrac1{24}}}<\frac1{90}.
\end{align*}
Since $H_0(v)$ decreases for $0\le v\le1/90$ and
$\pi H_0(1/90)/(12\sqrt2)>4/25$, we obtain
	\begin{align*}
		P'(d-u)>
		\frac4{25d^{3/2}}
		\exp\!\left(
		\pi\sqrt{\frac23(d-\tfrac{25}{24})}
		\right).
	\end{align*}
	Integration proves the first inequality
	in \eqref{eq:Delta-base}.
	
	Differentiating the integral twice gives
	\begin{align}\label{eq:Delta-derivatives}
		\Delta'(d)
		&=\int_0^1P''(d-u)\dd u,
		\qquad
		\Delta''(d)=\int_0^1P'''(d-u)\dd u.
	\end{align}
	By \eqref{eq:P-ratio-identities}
	and \eqref{eq:Delta-derivatives}, we obtain
	\begin{align}\label{eq:Delta-variance}
		(\log\Delta)''(d)
		&=
		\frac1{\Delta(d)}
		\int_0^1P'(d-u)
		\left(
		\frac{P'''(d-u)}{P'(d-u)}
		-\frac{P''(d-u)^2}{P'(d-u)^2}
		\right)\dd u
		\notag\\
		&\quad+
		\frac1{\Delta(d)}
		\int_0^1P'(d-u)
		\left(
		\frac{P''(d-u)}{P'(d-u)}
		-\frac{\Delta'(d)}{\Delta(d)}
		\right)^2\dd u.
	\end{align}
We first estimate the first integral in
\eqref{eq:Delta-variance}.
For $0\le u\le1$, \eqref{eq:P-ratio-derivative} gives
\begin{align*}
	\frac{P'''(d-u)}{P'(d-u)}
	-\left(\frac{P''(d-u)}{P'(d-u)}\right)^2
	<
	-\frac1{2(d-u-\tfrac1{24})^{3/2}}
	<
	-\frac1{2d^{3/2}}.
\end{align*}
Since $P'(d-u)>0$ and
$\int_0^1P'(d-u)\,\mathrm du=\Delta(d)$,
the first integral is less than $-1/(2d^{3/2})$.

We now estimate the second integral in
\eqref{eq:Delta-variance}.
By \eqref{eq:Delta-derivatives}, we have
\begin{align}\label{eq:Delta-ratio-mean}
	\frac{\Delta'(d)}{\Delta(d)}
	=
	\frac1{\Delta(d)}
	\int_0^1P'(d-u)
	\frac{P''(d-u)}{P'(d-u)}\,\mathrm du.
\end{align}
For $0\le u\le1$, we have $P'(d-u)/\Delta(d)>0$, and
$
	\int_0^1\frac{P'(d-u)}{\Delta(d)}\,\mathrm du=1.
$
Also, \eqref{eq:P-ratio-derivative} shows that
$P''(s)/P'(s)$ decreases with $s$.
Consequently, both
$P''(d-u)/P'(d-u)$ and
$\Delta'(d)/\Delta(d)$ lie between
$P''(d)/P'(d)$ and $P''(d-1)/P'(d-1)$.

By \eqref{eq:P-ratio-derivative}, the difference
between these two bounds satisfies
\begin{align}\label{eq:Delta-ratio-range}
	0<
	\frac{P''(d-1)}{P'(d-1)}
	-\frac{P''(d)}{P'(d)}
	&=
	-\int_{d-1}^{d}
	\frac{\mathrm d}{\mathrm ds}
	\frac{P''(s)}{P'(s)}\,\mathrm ds<
	\frac{257}{400(d-\tfrac{25}{24})^{3/2}}
	<
	\frac1{(d-2)^{3/2}}.
\end{align}
Thus, for $0\le u\le1$, we obtain
\begin{align*}
	\left|
	\frac{P''(d-u)}{P'(d-u)}
	-\frac{\Delta'(d)}{\Delta(d)}
	\right|
	<
	\frac1{(d-2)^{3/2}}.
\end{align*}
It follows that 
\begin{align*}
	\frac1{\Delta(d)}
	\int_0^1P'(d-u)
	\left(
	\frac{P''(d-u)}{P'(d-u)}
	-\frac{\Delta'(d)}{\Delta(d)}
	\right)^2\mathrm du
	<
	\frac1{(d-2)^3}.
\end{align*}

Combining the bounds for the two integrals in
\eqref{eq:Delta-variance}, we obtain
\begin{align*}
	(\log\Delta)''(d)
	<
	-\frac1{2d^{3/2}}+\frac1{(d-2)^3}
	\le
	-\frac1{2d^{3/2}}+\frac8{d^3}
	<
	-\frac1{3d^{3/2}},
\end{align*}
where the last two inequalities use
$d-2\ge d/2$ and $d^{3/2}>48$, respectively.
Both follow from $d\ge1252$.
This proves the second inequality in \eqref{eq:Delta-base}.
\end{proof}

\begin{lemma}\label{lem:P-growth}
	For $n\ge5000$ and real $s$ satisfying
	$n/4\le s-\tfrac1{24}\le n$, we have
	\begin{align}\label{eq:P-growth}
		\frac{\mathrm d}{\mathrm ds}
		\log\frac{P'(s)}{s-\tfrac1{24}}
		=
		\frac{P''(s)}{P'(s)}
		-\frac1{s-\tfrac1{24}}
		>
		\frac{31}{25\sqrt n}.
	\end{align}
\end{lemma}

\begin{proof}
	By \eqref{eq:P-ratio}, we have 
	\begin{align*}
		\frac{P''(s)}{P'(s)}
		-\frac1{s-\tfrac1{24}}
		\ge
		\frac{\pi}{\sqrt{6(s-\tfrac1{24})}}
		-\frac5{2(s-\tfrac1{24})}.
	\end{align*}
	The last expression decreases as $s-\tfrac1{24}$
	increases for $s-\tfrac1{24}\ge1250$.
	Since $s-\tfrac1{24}\le n$, we obtain
	\begin{align*}
		\frac{P''(s)}{P'(s)}
		-\frac1{s-\tfrac1{24}}
		\ge
		\frac1{\sqrt n}
		\left(
		\frac{\pi}{\sqrt6}-\frac5{2\sqrt n}
		\right)
		>
		\frac1{\sqrt n}
		\left(
		\frac{314}{245}-\frac1{28}
		\right)
		=
		\frac{1221}{980\sqrt n}
		>
		\frac{31}{25\sqrt n}.
	\end{align*}
This completes the proof.
\end{proof}

\subsection{Proof of Theorem~\ref{thm:strict}
	for $2\sqrt n\le m\le n/2+1$}

\begin{theorem}\label{prop:outer}
	Fix $k\in\{1,2\}$ and an integer $n\ge5000$.
	For every integer $m$ with $2\sqrt n\le m\le n/2+1$,
	we have $A_k(m+j,n)>0$ for $j=-1,0,1$, and
	\begin{align}\label{eq:outermargin}
		\Lambda_k(m,n)
		=\log\frac{A_k(m,n)^2}
		{A_k(m-1,n)A_k(m+1,n)}
		>\frac1{48n^{3/2}}.
	\end{align}
\end{theorem}

\begin{proof}
We construct a positive, strictly log-concave finite sum
on $[m-1,m+1]$, compare its values at $m+j$ with
$A_k(m+j,n)$ for $j=-1,0,1$, and apply
\eqref{eq:kernel} and \eqref{eq:transfer}
to prove \eqref{eq:outermargin}.
	
	\medskip\noindent
	\textbf{Step 1: a fixed finite sum.}
	For real $y$, define
	\begin{align}\label{eq:outer-sj}
		s_j(y)
		=n-\left(k-\tfrac12\right)j^2-yj+\frac j2,
		\qquad j\ge1.
	\end{align}
	Let $L$ be the largest positive integer satisfying
	\begin{align}\label{eq:outer-cutoff}
		s_L(m)\ge\frac n4+\sqrt n.
	\end{align}
	Such an integer exists because
	$s_1(m)\ge n/2-2>n/4+\sqrt n$,
	whereas $s_j(m)$ is negative
	for all sufficiently large $j$.
	Moreover,
	$s_L(m)\le n-L^2/2-(2\sqrt n-\tfrac12)L$.
	If $L\ge3\sqrt n/8$, this would give
	\begin{align*}
		s_L(m)
		\le
		\frac{23n}{128}+\frac{3\sqrt n}{16}
		<
		\frac n4+\sqrt n,
	\end{align*}
	contrary to \eqref{eq:outer-cutoff}.
	Since $\sqrt n>8$, we therefore have
	\begin{align}\label{eq:outer-L}
		1\le L<\frac{3\sqrt n}{8},
		\qquad
		L+1<\frac{\sqrt n}{2}.
	\end{align}
	
	For $y\in[m-1,m+1]$, define
	\begin{align}\label{eq:outer-S}
		S(y)=
		\sum_{j=1}^{L}(-1)^{j-1}
		\bigl(P(s_j(y))-P(s_j(y)-j)\bigr).
	\end{align}
	For $1\le j\le L$ and $0\le u\le j$, we have
	\begin{align}\label{eq:outer-retained-domain}
		s_j(y)-u-\frac1{24}
		\ge
		s_L(m)-2L-\frac1{24}
		>
		\frac n4+\frac{\sqrt n}{4}-\frac1{24}
		>
		\frac n4\ge1250.
	\end{align}
	Thus $S(y)$ and the retained terms have continuous
	derivatives of every order on $[m-1,m+1]$.
	
\medskip\noindent
\textbf{Step 2: proving $T''(y)\le0$ by pairing consecutive terms.}
Throughout this step, $m-1\le y\le m+1$.
By \eqref{eq:outer-sj}, we have $s_1(y)=n-y-k+1$ and
$n/2-3\le s_1(y)<n$.
For $1\le j\le L$ and $0\le u\le j$,
\eqref{eq:outer-retained-domain} and $s_j(y)\le s_1(y)<n$ give
\begin{align}\label{eq:outer-step2-domain}
	\frac n4<s_j(y)-u-\tfrac1{24}<n.
\end{align}
For $2\le j\le L$, define
\begin{align}\label{eq:outer-normalization}
	r_j(y)
	&=\frac{P(s_j(y))-P(s_j(y)-j)}{\Delta(s_1(y))}
	=\frac1{\Delta(s_1(y))}
	\int_0^jP'(s_j(y)-u)\dd u,
\end{align}
and
\begin{align}\label{eq:outer-normalization2}
	T(y)&=\sum_{j=2}^{L}(-1)^{j-1}r_j(y).
\end{align}
Then $S(y)=\Delta(s_1(y))(1+T(y))$.
We prove that $r_j''(y)$ is positive and decreases with $j$;
pairing consecutive terms will then give $T''(y)\le0$.

\begin{claim}\label{cl:outer-rj-bounds}
	For $2\le j\le L$, we have
	\begin{align}\label{eq:outer-rj-bounds}
		\frac{29(j-1)^2}{20}
		<
		\frac{\Delta(s_1(y))r_j''(y)}
		{\displaystyle
			\int_0^j
			\frac{P'(s_j(y)-u)}{s_j(y)-u-\tfrac1{24}}\dd u}
		<
		\frac{(257j-123)^2+520}{40000}.
	\end{align}
	In particular, $r_j''(y)>0$.
\end{claim}

\begin{proof}
	For $s-\tfrac1{24}\ge1250$,
	\eqref{eq:P-ratio} implies
	\begin{align}\label{eq:outer-ratio-constant}
		\frac{247}{200\sqrt{s-\tfrac1{24}}}
		<
		\frac{P''(s)}{P'(s)}
		<
		\frac{257}{200\sqrt{s-\tfrac1{24}}},
	\end{align}
	since $\pi/\sqrt6-3/(2\sqrt{1250})>247/200$.
For $2\le j\le L$ and $0\le u\le j$,
\eqref{eq:outer-sj} gives $s_1(y)-s_j(y)\ge y+1>2$.
Thus, for $0\le v\le1$, we have
$s_j(y)-u<s_1(y)-2<s_1(y)-v$.
Since $P''(s)/P'(s)$ decreases by
\eqref{eq:P-ratio-derivative},
\begin{align*}
	\frac{P''(s_1(y)-v)}{P'(s_1(y)-v)}
	\le
	\frac{P''(s_j(y)-u)}{P'(s_j(y)-u)}.
\end{align*}
Also, \eqref{eq:outer-ratio-constant} and
$s_1(y)-v-\tfrac1{24}<n$ give
\begin{align*}
	\frac{P''(s_1(y)-v)}{P'(s_1(y)-v)}
	>
	\frac{247}{200\sqrt{s_1(y)-v-\tfrac1{24}}}
	>
	\frac{247}{200\sqrt n}.
\end{align*}
Multiplying these inequalities by $P'(s_1(y)-v)>0$
and integrating from $v=0$ to $v=1$, using $\Delta(s_1(y))=P(s_1(y))-P(s_1(y)-1)$
and \eqref{eq:Delta-derivatives}, we obtain
\begin{align*}
	\frac{247}{200\sqrt n}\Delta(s_1(y))
	&=
	\frac{247}{200\sqrt n}
	\int_0^1P'(s_1(y)-v)\dd v\\
	&<
	\int_0^1P''(s_1(y)-v)\dd v\\
	&\le
	\frac{P''(s_j(y)-u)}{P'(s_j(y)-u)}
	\int_0^1P'(s_1(y)-v)\dd v\\
	&=
	\frac{P''(s_j(y)-u)}{P'(s_j(y)-u)}
	\Delta(s_1(y)).
\end{align*}
Dividing by $\Delta(s_1(y))>0$ and using $\int_0^1P''(s_1(y)-v)\dd v
=\Delta'(s_1(y))$, we obtain
\begin{align}\label{eq:outer-Delta-ratio}
	\frac{247}{200\sqrt n}
	<
	\frac{\Delta'(s_1(y))}{\Delta(s_1(y))}
	\le
	\frac{P''(s_j(y)-u)}{P'(s_j(y)-u)}.
\end{align}
By \eqref{eq:outer-step2-domain}, we have 
$\sqrt{s_j(y)-u-\tfrac1{24}}>\sqrt n/2$. Then 
\begin{align}\label{eq:outer-Delta-ratio-scaled}
	\frac{\Delta'(s_1(y))}{\Delta(s_1(y))}
	\sqrt{s_j(y)-u-\tfrac1{24}}
	>
	\frac{247}{200\sqrt n}\frac{\sqrt n}{2}
	=
	\frac{247}{400}
	>
	\frac{123}{200}.
\end{align}
	
	By \eqref{eq:Delta-base}, we have 
	$-(\log\Delta)''(s_1(y))>0$.
	Since the second integral in \eqref{eq:Delta-variance}
	is nonnegative, \eqref{eq:P-ratio-derivative} implies that
	$-(\log\Delta)''(s_1(y))
	<257/(400(s_1(y)-2)^{3/2})$.
	Since $s_j(y)-u<s_1(y)-2$ and $s_1(y)\ge2497$,
	we obtain
	\begin{align}\label{eq:outer-Delta-curvature-upper}
		0<
		-(\log\Delta)''(s_1(y))
		\bigl(s_j(y)-u-\tfrac1{24}\bigr)
		<
		\frac{257}{400\sqrt{2495}}
		<\frac{13}{1000}.
	\end{align}
	
Since $s_j'(y)=-j$ and $s_1'(y)=-1$,
the product and chain rules give
\begin{align}\label{eq:outer-rj-first}
	r_j'(y)
	&=
	-\frac{j}{\Delta(s_1(y))}
	\int_0^jP''(s_j(y)-u)\dd u
	+
	\frac{\Delta'(s_1(y))}{\Delta(s_1(y))^2}
	\int_0^jP'(s_j(y)-u)\dd u
	\notag\\
	&=
	\frac1{\Delta(s_1(y))}
	\int_0^jP'(s_j(y)-u)
	\left(
	\frac{\Delta'(s_1(y))}{\Delta(s_1(y))}
	-j\frac{P''(s_j(y)-u)}{P'(s_j(y)-u)}
	\right)\dd u.
\end{align}
By \eqref{eq:outer-rj-first}, we get 
	\begin{align}\label{eq:rsecond}
		r_j''(y)
		={}&
		\frac1{\Delta(s_1(y))}
		\int_0^jP'(s_j(y)-u)
		\Biggl(
		\left(
		j\frac{P''(s_j(y)-u)}{P'(s_j(y)-u)}
		-\frac{\Delta'(s_1(y))}{\Delta(s_1(y))}
		\right)^2
		\notag\\
		&\quad+
		j^2\left(
		\frac{P'''(s_j(y)-u)}{P'(s_j(y)-u)}
		-\frac{P''(s_j(y)-u)^2}{P'(s_j(y)-u)^2}
		\right)
		-(\log\Delta)''(s_1(y))
		\Biggr)\dd u.
	\end{align}
	Here we used
	$\frac{\mathrm d}{\mathrm dy}
	(\Delta'(s_1(y))/\Delta(s_1(y)))
	=-(\log\Delta)''(s_1(y))$,
	which follows from $s_1'(y)=-1$.
	
By \eqref{eq:outer-Delta-ratio}, we have
\begin{align*}
	j\frac{P''(s_j(y)-u)}{P'(s_j(y)-u)}
	-\frac{\Delta'(s_1(y))}{\Delta(s_1(y))}
	\ge
	(j-1)\frac{P''(s_j(y)-u)}{P'(s_j(y)-u)}
	>0.
\end{align*}
Thus \eqref{eq:outer-ratio-constant} and
\eqref{eq:outer-Delta-ratio-scaled} give
\begin{align*}
	\frac{(247/200)^2(j-1)^2}{s_j(y)-u-\tfrac1{24}}
	&<
	\left(
	j\frac{P''(s_j(y)-u)}{P'(s_j(y)-u)}
	-\frac{\Delta'(s_1(y))}{\Delta(s_1(y))}
	\right)^2
	<
	\frac{(257j-123)^2}
	{40000(s_j(y)-u-\tfrac1{24})}.
\end{align*}
Moreover, by \eqref{eq:P-ratio-derivative}, we obtain 
\begin{align*}
	j^2\left(
	\frac{P'''(s_j(y)-u)}{P'(s_j(y)-u)}
	-\frac{P''(s_j(y)-u)^2}{P'(s_j(y)-u)^2}
	\right)<0.
\end{align*}
	By \eqref{eq:P-ratio-derivative}, it is greater than
	$-257(j-1)^2/(3500(s_j(y)-u-\tfrac1{24}))$,
	because $j^2\le4(j-1)^2$ and
	$\sqrt{s_j(y)-u-\tfrac1{24}}>35$.
	Finally, \eqref{eq:outer-Delta-curvature-upper}
	bounds the last term.
	Using
	\begin{align*}
		\left(\frac{247}{200}\right)^2-\frac{257}{3500}
		>\frac{29}{20},
		\qquad
		\frac{13}{1000}=\frac{520}{40000},
	\end{align*}
	we integrate the resulting bounds in \eqref{eq:rsecond}
	to obtain \eqref{eq:outer-rj-bounds}.
	Its denominator is positive by \eqref{eq:Pprime-form},
	so $r_j''(y)>0$.
\end{proof}

\begin{claim}\label{cl:outer-adjacent-second}
	For $2\le j<L$, we have $r_{j+1}''(y)<r_j''(y)$.
\end{claim}

\begin{proof}
	By \eqref{eq:outer-sj}, we get 
	\begin{align}\label{eq:outer-adjacent-separation}
		s_j(y)-j-s_{j+1}(y)
		=y+(k-1)(2j+1)\ge y.
	\end{align}
	By \eqref{eq:P-growth}, we have 
	$P'(s)/(s-\tfrac1{24})$ increases.
For $2\le j<L$, the monotonicity in \eqref{eq:P-growth} gives
\begin{align*}
	\int_0^j
	\frac{P'(s_j(y)-u)}{s_j(y)-u-\tfrac1{24}}\dd u
	&\ge
	j\frac{P'(s_j(y)-j)}{s_j(y)-j-\tfrac1{24}},\\
	\int_0^{j+1}
	\frac{P'(s_{j+1}(y)-u)}{s_{j+1}(y)-u-\tfrac1{24}}\dd u
	&\le
	(j+1)\frac{P'(s_{j+1}(y))}
	{s_{j+1}(y)-\tfrac1{24}}.
\end{align*}
Integrating both sides of \eqref{eq:P-growth}
from $s_{j+1}(y)$ to $s_j(y)-j$ and using
\eqref{eq:outer-adjacent-separation}, we obtain
\begin{align*}
	\log\frac{
		P'(s_j(y)-j)/(s_j(y)-j-\tfrac1{24})
	}{
		P'(s_{j+1}(y))/(s_{j+1}(y)-\tfrac1{24})
	}
	>
	\frac{31(s_j(y)-j-s_{j+1}(y))}{25\sqrt n}
	\ge\frac{31y}{25\sqrt n}.
\end{align*}
Therefore, 
\begin{align*}
	\frac{
		P'(s_{j+1}(y))/(s_{j+1}(y)-\tfrac1{24})
	}{
		P'(s_j(y)-j)/(s_j(y)-j-\tfrac1{24})
	}
	<
	\exp\!\left(-\frac{31y}{25\sqrt n}\right).
\end{align*}
	Dividing the upper bound for $r_{j+1}''(y)$
	by the lower bound for $r_j''(y)$ in
	\eqref{eq:outer-rj-bounds} therefore gives
	\begin{align}\label{eq:outer-adjacent-ratio}
		\frac{r_{j+1}''(y)}{r_j''(y)}
		<
		\frac{(j+1)((257j+134)^2+520)}
		{58000j(j-1)^2}
		\exp\!\left(-\frac{31y}{25\sqrt n}\right).
	\end{align}
	For $j\ge2$, we have
	$(j+1)/j\le3/2$ and
	$(257j+134)/(j-1)=257+391/(j-1)\le648$.
	Thus the fraction on the right of
	\eqref{eq:outer-adjacent-ratio}
	is at most $3(648^2+520)/116000$.
	Also, $y\ge2\sqrt n-1$ and $\sqrt n>70$ give
	\begin{align}\label{eq:outer-geometric-rate}
		\frac{31y}{25\sqrt n}
		>\frac{4309}{1750}>\frac{49}{20},
		\qquad e^{49/20}>\frac{80}{7}.
	\end{align}
	Substituting into \eqref{eq:outer-adjacent-ratio},
	we get
	\begin{align*}
		\frac{r_{j+1}''(y)}{r_j''(y)}
		<
		\frac{3(648^2+520)}{116000}\frac7{80}
		<1.
	\end{align*}
	This proves the claim.
\end{proof}

For $m-1\le y\le m+1$, we have $T''(y)=0$ if $L=1$.
For $L\ge2$, Claim~\ref{cl:outer-adjacent-second} gives
$-r_{2\ell}''(y)+r_{2\ell+1}''(y)<0$
for every positive integer $\ell$ with $2\ell+1\le L$.
If $L$ is even, the additional term $-r_L''(y)$
is negative by Claim~\ref{cl:outer-rj-bounds}.
Therefore,
\begin{align}\label{eq:outer-T-second}
	T''(y)\le0.
\end{align}
	
	\medskip\noindent
	\textbf{Step 3: positivity and strict log-concavity of $S(y)$.} 
	For $2\le j\le L$ and $s_j(y)-j\le s\le s_1(y)$,
	\eqref{eq:outer-step2-domain} and \eqref{eq:P-growth} give
	$P''(s)/P'(s)>31/(25\sqrt n)>0$.
	Thus $P'(s)$ increases with $s$ on this interval.
	
	By \eqref{eq:outer-normalization} and
	\eqref{eq:Delta-integral}, we obtain
	\begin{align}\label{eq:outer-rj-size}
		0<r_j(y)
		=
		\frac{
			\int_0^jP'(s_j(y)-u)\dd u
		}{
			\int_0^1P'(s_1(y)-v)\dd v
		}
		\le
		\frac{jP'(s_j(y))}{P'(s_1(y)-1)}
		\le
		j\exp\!\left(
		-\frac{31((j-1)y-1)}{25\sqrt n}
		\right).
	\end{align}
	The last inequality in \eqref{eq:outer-rj-size}
	follows by integrating
	$P''(s)/P'(s)>31/(25\sqrt n)$
	from $s_j(y)$ to $s_1(y)-1$, since
	$s_1(y)-1-s_j(y)\ge(j-1)y-1$.
Using \eqref{eq:outer-rj-size},
\eqref{eq:outer-geometric-rate}, and
$\exp(31/(25\sqrt n))<25/24$, we obtain
\begin{align*}
	|T(y)|
	&=
	\left|
	\sum_{j=2}^{L}(-1)^{j-1}r_j(y)
	\right|
	\le
	\sum_{j=2}^{L}r_j(y)
	<
	\frac{25}{24}
	\sum_{j=2}^{\infty}
	j\left(\frac7{80}\right)^{j-1}
	=
	\frac{8925}{42632}<\frac14.
\end{align*}
	Since $31/(25\sqrt n)<1/25$ and
	$e^x\le(1-x)^{-1}$ for $0\le x<1$, we have
	\begin{align}\label{eq:outer-exponential-bound}
		\exp\!\left(\frac{31}{25\sqrt n}\right)
		\le
		\frac1{1-31/(25\sqrt n)}
		<\frac{25}{24}.
	\end{align}
	
	The estimate $|T(y)|<1/4$ and the positivity of
	$\Delta(s_1(y))$ give
	\begin{align}\label{eq:S-positive}
		S(y)=\Delta(s_1(y))(1+T(y))
		>\frac34\Delta(s_1(y))>0.
	\end{align}
	
	Next, since $s_1'(y)=-1$ and $s_1''(y)=0$,
	differentiating $\log S(y)$ twice gives
	\begin{align}\label{eq:Scurvature}
		(\log S)''(y)
		&=
		(\log\Delta)''(s_1(y))
		+\frac{T''(y)}{1+T(y)}
		-\frac{T'(y)^2}{(1+T(y))^2}
		<
		-\frac1{3s_1(y)^{3/2}}
		\le
		-\frac1{3n^{3/2}}.
	\end{align}
	Here the first inequality follows from
	\eqref{eq:Delta-base}, \eqref{eq:outer-T-second},
	and $1+T(y)>0$; the second uses $s_1(y)<n$.
	
	Finally, substitute the first bound in
	\eqref{eq:Delta-base} into \eqref{eq:S-positive}.
	Since $s_1(y)-25/24>s_1(y)-2$, we obtain
	\begin{align}\label{eq:Slower}
		S(y)>
		\frac3{25s_1(y)^{3/2}}
		\exp\!\left(
		\pi\sqrt{\frac23(s_1(y)-2)}
		\right).
	\end{align}
	
	\medskip\noindent
	\textbf{Step 4: comparison with the integer coefficients.}
	
Now let $y\in\{m-1,m,m+1\}$.
By \eqref{eq:coeff} and \eqref{eq:outer-S},
we estimate $A_k(y,n)-S(y)$ by bounding the terms
with $j\ge L+1$ and the differences $p(s)-P(s)$
at $s=s_j(y)$ and $s=s_j(y)-j$ for $1\le j\le L$.

If $s_{L+1}(y)<0$, all terms with $j\ge L+1$ vanish.
Otherwise, Lemma~\ref{lem:partition-tail} gives
\begin{align}\label{eq:outer-partition-tail}
	\left|
	\sum_{j=L+1}^{\infty}(-1)^{j-1}
	\bigl(p(s_j(y))-p(s_j(y)-j)\bigr)
	\right|
	\le p(s_{L+1}(y)).
\end{align}
By maximality of $L$ in \eqref{eq:outer-cutoff},
$s_{L+1}(m)<n/4+\sqrt n$.
Since $|y-m|\le1$, \eqref{eq:outer-L} gives
\begin{align}\label{eq:outer-tail-range}
	s_{L+1}(y)
	=s_{L+1}(m)-(L+1)(y-m)
	\le s_{L+1}(m)+L+1
	<\frac n4+\frac32\sqrt n
	<\frac{3n}{10}.
\end{align}
Since $s_1(y)-2\ge n/2-5>49n/100$, we have
\begin{align}\label{eq:outer-square-root-bounds}
	\sqrt{s_1(y)-2}>\frac7{10}\sqrt n.
\end{align}

Also, \eqref{eq:outer-tail-range} gives
$\max(s_{L+1}(y),n/8)<3n/10$.
Consequently,
\begin{align}\label{eq:outer-tail-root-bound}
	\sqrt{\max(s_{L+1}(y),n/8)}
	<
	\sqrt{\frac{3n}{10}}
	<
	\frac{11}{20}\sqrt n.
\end{align}

We now divide the upper bound for $p(s_{L+1}(y))$
in \eqref{eq:envelopes} by the lower bound for $S(y)$
in \eqref{eq:Slower}, obtaining
\begin{align}\label{eq:outer-omitted-bound}
	\frac{p(s_{L+1}(y))}{S(y)}
	&<
	\frac{10s_1(y)^{3/2}}n
	\exp\!\left(
	-\pi\sqrt{\frac23}
	\left(
	\sqrt{s_1(y)-2}
	-\sqrt{\max(s_{L+1}(y),n/8)}
	\right)
	\right)
	\notag\\
	&<
	10\sqrt n\,e^{-3\sqrt n/8}.
\end{align}
For the last inequality, \eqref{eq:outer-square-root-bounds}
and \eqref{eq:outer-tail-root-bound} give a difference
greater than $(7/10-11/20)\sqrt n=3\sqrt n/20$.
We also used $s_1(y)<n$ and $\pi\sqrt{2/3}>5/2$.
	
For $1\le j\le L$ and
$s\in\{s_j(y),s_j(y)-j\}$,
\eqref{eq:outer-retained-domain} gives
$s-\tfrac1{24}>n/4$ and $s\le s_1(y)$.
By \eqref{eq:p-error}, we have
\begin{align}\label{eq:outer-single-difference}
	|p(s)-P(s)|
	\le
	\frac{\exp\!\left(\pi\sqrt{(s-\tfrac1{24})/6}\right)}
	{8(s-\tfrac1{24})}
	<
	\frac1{2n}
	\exp\!\left(\pi\sqrt{s_1(y)/6}\right).
\end{align}

Since each index $j$ gives two differences,
there are $2L<3\sqrt n/4$ differences in total.
By \eqref{eq:outer-single-difference} and
the lower bound for $S(y)$ in \eqref{eq:Slower}, we obtain
\begin{align}\label{eq:outer-retained-bound}
	&\frac1{S(y)}
	\sum_{j=1}^{L}
	\Bigl(
	|p(s_j(y))-P(s_j(y))|
	+
	|p(s_j(y)-j)-P(s_j(y)-j)|
	\Bigr)
	\notag\\
	&\quad<
	\frac{25s_1(y)^{3/2}}{8\sqrt n}
	\exp\!\left(
	-\pi\sqrt{\frac23}
	\left(
	\sqrt{s_1(y)-2}-\frac{\sqrt{s_1(y)}}2
	\right)
	\right)
	\notag\\
	&\quad<
	4n e^{-\sqrt n/2}.
\end{align}
For the last inequality,
\eqref{eq:outer-square-root-bounds} and $s_1(y)<n$ give
$\sqrt{s_1(y)-2}-\sqrt{s_1(y)}/2>\sqrt n/5$.
Since $\pi\sqrt{2/3}>5/2$ and $25/8<4$,
the stated bound follows.

By \eqref{eq:monotone}, we get 
$\sqrt n\,e^{-\sqrt n/8}$ decreases for $n\ge5000$.
Since $e^3>20$, we have
\begin{align}\label{eq:outer-first-decay}
	4\sqrt n\,e^{-\sqrt n/8}
	\le
	4\sqrt{5000}\,e^{-\sqrt{5000}/8}
	<1.
\end{align}
Also, by \eqref{eq:monotone},
$n^2e^{-3\sqrt n/8}$ decreases for $n\ge5000$.
It follows that
\begin{align}\label{eq:outer-second-decay}
	176n^2e^{-3\sqrt n/8}
	&\le
	176\cdot5000^2e^{-3\sqrt{5000}/8}
	<1.
\end{align}

By \eqref{eq:outer-omitted-bound} and
\eqref{eq:outer-retained-bound}, the triangle inequality gives
\begin{align}\label{eq:outer-relative-bound}
	\left|\frac{A_k(y,n)}{S(y)}-1\right|
	<
	10\sqrt n\,e^{-3\sqrt n/8}
	+4n e^{-\sqrt n/2}
	<
	11\sqrt n\,e^{-3\sqrt n/8}
	<
	\frac1{16n^{3/2}}.
\end{align}
The second inequality uses \eqref{eq:outer-first-decay},
and the last uses \eqref{eq:outer-second-decay}.

\medskip\noindent
\textbf{Step 5: strict log-concavity of $A_k(m,n)$.}
	By \eqref{eq:S-positive}, we have $S(y)>0$ for
	$m-1\le y\le m+1$.
	The finite sum in \eqref{eq:outer-S} has continuous
	first and second derivatives on this interval.
	Since $\int_{-1}^{1}(1-|t|)\dd t=1$,
	\eqref{eq:kernel} and \eqref{eq:Scurvature} give
	\begin{align}\label{eq:outer-S-discrete-bound}
		\log\frac{S(m)^2}{S(m-1)S(m+1)}
		&=
		-\int_{-1}^{1}(1-|t|)
		(\log S)''(m+t)\dd t>\frac1{3n^{3/2}}.
	\end{align}
	
	For $j=-1,0,1$, take
	$x_j=S(m+j)$, $b_j=A_k(m+j,n)$, and
	$\eta=1/(16n^{3/2})$ in Lemma~\ref{lem:transfer}.
	By \eqref{eq:S-positive} and
	\eqref{eq:outer-relative-bound}, we have
	\begin{align}\label{eq:outer-transfer-conditions}
		x_j>0,\qquad
		\left|\frac{b_j}{x_j}-1\right|<\eta,
		\qquad 0<\eta<\frac15.
	\end{align}
	Therefore Lemma~\ref{lem:transfer} applies and gives
	$A_k(m+j,n)>0$ for $j=-1,0,1$.
	
	Since $\eta<1/5$, we have $4\eta/(1-\eta)<5\eta$.
	By \eqref{eq:transfer} and
	\eqref{eq:outer-S-discrete-bound}, we obtain
	\begin{align}\label{eq:outer-final-bound}
		\Lambda_k(m,n)
		\ge
		\log\frac{S(m)^2}{S(m-1)S(m+1)}
		-\frac{4\eta}{1-\eta}>
		\frac1{3n^{3/2}}
		-\frac5{16n^{3/2}}
		=\frac1{48n^{3/2}}.
	\end{align}
	This completes the proof.
\end{proof}

\section{Proof of the main theorem}
\label{sec:completion-conjectures}

We combine the estimates in
Sections~\ref{sec:central} and~\ref{sec:outer},
and Zhou's inequality \eqref{eq:zhou}
to prove Theorem~\ref{thm:strict}.

\begin{proof}[Proof of Theorem~\ref{thm:strict}]
	Fix $k\in\{1,2\}$ and integers $n\ge72$ and $m$
	satisfying $|m|\le n-k-71$.
	By \eqref{eq:Ak} and \eqref{eq:coeff},
	we have $A_k(-m,n)=A_k(m,n)$.
	Hence it suffices to consider
	$0\le m\le n-k-71$.

	For $72\le n<5000$, direct computation with
	Mathematica verifies \eqref{eq:strict-main}
	for $k=1,2$ and every integer
	$0\le m\le n-k-71$.
	We now assume $n\ge 5000$.

    For $0\le m\le2\sqrt n$,
	by Theorem~\ref{prop:central},
    we have 
	$A_k(m-1,n), A_k(m,n), A_k(m+1,n)>0$ and 
	\begin{align*}
	\Lambda_k(m,n)=
	\log\frac{A_k(m,n)^2}{A_k(m-1,n)A_k(m+1,n)}>\frac1{20000n}>0.
	\end{align*}
    For $2\sqrt n<m\le n/2+1$,
	by Theorem~\ref{prop:outer}, we also have 
	$A_k(m-1,n), A_k(m,n), A_k(m+1,n)>0$ and 
	\begin{align*}
	\Lambda_k(m,n)=
	\log\frac{A_k(m,n)^2}{A_k(m-1,n)A_k(m+1,n)}>\frac1{48n^{3/2}}>0.
	\end{align*}
	Thus for every $n\ge 5000,$
    we have established strict log-concavity of $A_k(m,n)$ for $k=1,2$ and $0\le m\le n/2+1.$
    Zhou~\cite[Proposition 1.6]{Zhou} proved that 
    strict log-concavity of $A_k(m,n)$ for $k=1,2$ and $n/2+1\le m\le n-k-71.$
	  This completes the proof.
\end{proof}

\end{document}